\documentclass[11pt]{amsart}

\usepackage{amsmath,amsthm, amssymb}
\usepackage[T1]{fontenc}
\usepackage{latexsym}
\usepackage{subfigure, epsfig, epsf}
\usepackage{color,fancyhdr,lastpage,graphicx}
\usepackage{xspace}
\usepackage{wasysym}
\usepackage{setspace}
\usepackage{bbm}
\usepackage{stmaryrd}
\usepackage{titlesec}

\newcommand{\NN}{\mathbb{N}}      
\newcommand{\RR}{\mathbb{R}}    \newcommand{\LL}{\mathbb{L}}    
   \newcommand{\VV}{\mathbb{V}}

 \newcommand{\mcC}{\mathcal{C}}  \newcommand{\mcF}{\mathcal{F}} 
    
     \newcommand{\mcD}{\mathcal{D}}   
  \newcommand{\mcL}{\mathcal{L}}

    \newcommand{\bfX}{{\bf X}}      
\newcommand{\bfV}{{\bf V}}

\newcommand{\bbX}{\mathbb{X}}   \newcommand{\bbM}{\mathbb{M}}

\newcommand{\EE}{\mathbb{E}} \newcommand{\PP}{\mathbb{P}}

\newcommand{\ep}{\epsilon}
\renewcommand{\leq}{\leqslant}
\renewcommand{\geq}{\geqslant}

\newcommand{\ssk}{\smallskip}

\titleformat{\section}[block]
  {\normalfont\sffamily}{\thesection}{.5em}{\titlerule\\[.8ex]\bfseries} 

\titleformat{\subsection}[wrap]{\sffamily\fontseries{b}\selectfont\filright}{\thesubsection.}{.5em}{}
                  \titlespacing{\subsection}{20pc}{1.5ex plus .1ex minus .2ex}{1pc}

\titleformat{\subsubsection}[wrap]{\sffamily\fontseries{b}\selectfont\filright}{\thesubsubsection.}{.5em}{}
                  \titlespacing{\subsubsection}{20pc}{1ex plus .1ex minus .2ex}{1pc}

\numberwithin{subsection}{section}
\numberwithin{subsubsection}{subsection}

\newtheoremstyle{mystyle}
{3pt}% ⟨Space above⟩
{3pt}% ⟨Space below⟩
{\sffamily}% ⟨Body font⟩
{0em}% ⟨Indent amount⟩
{\bfseries\sffamily}% ⟨Theorem head font⟩
{.}% ⟨Punctuation after theorem head⟩
{.5em}% ⟨Space after theorem head⟩
{}% ⟨Theorem head spec (can be left empty, meaning ‘normal’)⟩
\theoremstyle{mystyle}
\newtheorem{thm}{Theorem}

\newtheorem{cor}[thm]{\hspace{-0.15cm}  {Corollary} }% [section]
\newtheorem{lem}[thm]{\hspace{-0.15cm}  {Lemma} }%[section]
\newtheorem{prop}[thm]{\hspace{-0.15cm} {Proposition}}%[chapter]
\newtheorem{defn}[thm]{ \hspace{-0.3cm} {Definition}}%[chapter]
   \newtheorem*{definition}{ \hspace{-0.3cm} {Definition}}%[chapter]
\newtheorem{rem}[thm]{Remark}

   \newtheorem*{example}{Example}

\numberwithin{equation}{section} % Une équation p.q. est la q ème équation de la section p.

\newenvironment{Dem}{%
    \begin{list}{\hspace{0.5cm}{\textbf{\textsf{Proof --}}}}{%
        \setlength{\topsep}{0pt}%
        \setlength{\leftmargin}{0pt}%
        \setlength{\rightmargin}{0pt}%
        \setlength{\listparindent}{0pt}%
        \setlength{\itemindent}{0pt}%
        \setlength{\parsep}{0pt}%
        \addtolength{\leftmargin}{20pt}%
        \addtolength{\rightmargin}{0pt}%
    } \item }{\hfill{\space $\rhd$}\end{list}\smallskip}

\begin{document}
\title{Rough flows}

\author{I. Bailleul}
\address{IRMAR, 263 Avenue du General Leclerc, 35042 RENNES, France}
\email{ismael.bailleul@univ-rennes1.fr}
\thanks{}

\author{S. Riedel}
\address{Sebastian Riedel \\
Institut f\"ur Mathematik, Technische Universit\"at Berlin, Germany}
\email{riedel@math.tu-berlin.de}

\begin{abstract}
We introduce in this work a concept of rough driver that somehow provides a rough path-like analogue of an enriched object associated with time-dependent vector fields. We use the machinery of approximate flows to build the integration theory of rough drivers and prove well-posedness results for rough differential equations on flows and continuity of the solution flow as a function of the generating rough driver. We show that the theory of semimartingale stochastic flows developed in the 80's and early 90's fits nicely in this framework, and obtain as a consequence some strong approximation results for general semimartingale flows and provide a fresh look at large deviation theorems for 'Gaussian' stochastic flows.
\end{abstract}

\date{\today}

\vspace*{3ex minus 1ex}
\begin{center}
\huge\sffamily{Rough flows}
\end{center}
\vskip 5ex minus 1ex

\begin{center}
{\sf I. BAILLEUL\footnote{I.B. was partly supported by the ANR project ''Retour Post-doctorant'', no. 11-PDOC-0025; I.B. also thanks the U.B.O. for their hospitality, part of this work was written there.} and S. RIEDEL}
\end{center}

\vspace{1cm}

\begin{center}
\begin{minipage}{0.8\textwidth}
\renewcommand\baselinestretch{0.7} \rmfamily {\scriptsize {\bf \sc \noindent Abstract.} We introduce in this work a concept of rough driver that somehow provides a rough path-like analogue of an enriched object associated with time-dependent vector fields. We use the machinery of approximate flows to build the integration theory of rough drivers and prove well-posedness results for rough differential equations on flows and continuity of the solution flow as a function of the generating rough driver. We show that the theory of semimartingale stochastic flows developed in the 80's and early 90's fits nicely in this framework, and obtain as a consequence some new strong approximation results for general semimartingale flows and provide a fresh look at large deviation theorems for 'Gaussian' stochastic flows. }
\end{minipage}
\end{center}

\bigskip

\section*{Introduction}

An elementary construction recipe of flows was recently introduced in \cite{B15} and used there to get back the core results of Lyons' theory of rough differential equations \cite{Lyo98, FV10} in a very short and elementary way. This work emphasizes the fact that it may be worth considering flows of maps as the primary objects from which the individual trajectories can be built, as opposed to the classical point of view that constructs a flow from an uncountable collection of individual trajectories. Probabilists know how tricky it can be to deal with uncountably many null sets. As well-recognized now, the main success of Lyons' theory was to disentangle probability from pure dynamics in the study of stochastic differential equations by showing that the dynamics is a deterministic and continuous function of an enriched signal that is constructed from the noise in the equation by purely probabilistic means. This very clean picture led to new proofs and extensions of foundational results in the theory of 
stochastic differential equations, such as Stroock and Varadhan support theorem or the basics of Freidlin and Ventzel theory of large deviations for diffusions.

\medskip

It was realized in the late 70's that stochastic differential equations not only define individual trajectories, they also define flows of regular homeomorphisms, depending on the regularity of the vector fields involved in the dynamics \cite{El78, B81, IW81, Ku81}. This opened the door to the study of stochastic flows of maps for themselves \cite{Ha81, Bax80, LeJan82, LeJan85} and it did not take long time before Le Jan and Watanabe \cite{LeJanW82} clarified definitely the situation by showing that, in a semimartingale setting, there is a one-to-one correspondence between flows of diffeomorphisms and time-varying stochastic velocity fields, under proper regularity conditions on the objects involved. We offer in the present work an embedding of the theory of semimartingale stochastic flows into the theory of rough flows similar to the embedding of the theory of stochastic differential equations into the theory of rough differential equations. While the acquainted reader will have noticed that the latter 
framework can be used to deal with Brownian flows by seeing them as solution flows to some Banach space-valued rough differential equations driven by a Brownian rough path in some space of vector fields, such as done in \cite{Der10, DD12}, the situation is not so clear for more general random velocity fields and stochastic flows of maps. Our approach provides a simple setting for dealing with the general case, and we provide in this work an elementary direct approach to the construction of stochastic flows whose scope goes beyond the realm of semimartingale calculus.

\medskip

Our theory of rough flows is based on the deterministic "Approximate flow-to-flow" machinery introduced in \cite{B15}, which gives body to the following fact. To a 2-index family $(\mu_{ts})_{0\leq s\leq t\leq T}$ of maps which falls short from being a flow one can associate a unique flow $(\varphi_{ts})_{0\leq s\leq t\leq T}$ close to $(\mu_{ts})_{0\leq s\leq t\leq T}$; moreover the flow $\varphi$ depends continuously on the approximate flow $\mu$. We introduce a notion of rough driver, that is an enriched version of a time-dependent vector field, that is given by the additional datum of a time-dependent second order differential operator satisfying some algebraic and analytic conditions. A notion of solution to a differential equation driven by a rough driver will be given, in the line of what was done in \cite{B15} for rough differential equations, and the approximate flow-to-flow machinery will be seen to lead to a clean and simple well-posedness result for such equations. Importantly, the It\^o map, 
that associates to a rough driver the solution flow to its associated equation, is continuous. This continuity result is the key to deep results in the theory of stochastic flows. We shall indeed prove that reasonable semimartingale velocity fields can be lifted to rough drivers under some mild boundedness and regularity conditions, and that the solution flow associated to the semimartingale rough driver coincides almost surely with the solution flow to the Kunita-type Stratonovich differential equation driven by the velocity field. In this sense, our theory of rough flows encompasses the theory of stochastic flows. A Wong-Zaka\" i theorem will be proved for a general class of semimartingale velocity fields and a sharp large deviation principle for Brownian flows will be proved will be proved as a consequence of the continuity of the It\^o map. No Wong-Zaka\" i result was available so far for semimartingale stochastic flows.

\medskip

The setting of rough drivers and rough flows is presented in Section \ref{SectionRoughFlows}, together with the approximate flow-to-flow machinery. This is the core of the deterministic machinery and everything that follows elaborates on this material, in a probabilistic setting. Some additional material on random rough drivers is in particular given in Section \ref{SubsectionKolmogorovTheorem}, where we provide some new variations around the Kolmogorov regularity theorem needed along the way and some sufficient conditions for a process to be bounded; this material may be of independent interest. We show in Section \ref{SectionStochasticRoughFlows} that reasonable semimartingale velocity fields can be lifted to rough drivers under appropriate mild boundedness and regularity assumptions and prove that the theory of semimartingale stochastic flows of maps is naturally embedded in the theory of rough flows. As an illustration of the continuity of the It\^o map, we prove in Section \ref{SubsectionStrongApproximations} a Wong-Zaka\" i theorem for stochastic flows of maps, and provide in Section \ref{SectionLDPSupp} a fresh look at large deviation theorems for Brownian stochastic flows. We proved in the follow up work \cite{BailleulRiedelScheutzow} that a large class of Gaussian vector fields can be lifted into rough drivers; this shows explicitly that the setting of rough drivers and rough flows goes beyond the semimartingale horizon. The point of view of rough flows presented here was also used in \cite{BailleulCatellier} to investigate the problem of stochastic turbulence.

\smallskip

The size of this work is related to the fact that it is intended to be as self-contained as possible; no a priori knowledge of the theory of stochastic flows is required for its understanding in particular. We have included as a consequence some material that is well-known from experts in stochastic flows or large deviation theory for instance. A reader interested only in the machinery of rough flows will have a complete picture by reading sections \ref{SubsectionApproximateFlows} to \ref{SubsectionItoFormula} and Theorem \ref{thm:kolmogorov_rough_drivers} in Section \ref{SubsectionKolmogorovTheorem}.

\bigskip

\noindent \textbf{\textsf{Notations.}} We gather here for reference a number of notations that will be used throughout the text. 

\ssk

\begin{itemize}
   \item We shall exclusively use the letter E to denote a Banach space; we shall denote by $\textrm{L}(\textrm{E})$ the set of continuous linear maps from E to itself, and for $M\in\textrm{L}(\textrm{E})$, we shall write $|M|$ for its operator norm. In this setting, differentiability and regularity notions are understood in the sense of Fr\'echet.\vspace{0.1cm}
   
   \item For functions $x \colon [0,T] \to E$, we will use the notation $x_{ts} = x_t -x_s$ for the increments of $x$.   
   
   \item Whenever useful, vector fields are identified with the first order differential operator they define in a canonical way. 
   
   \item As is common, we shall use Einstein's summation convention that $a^ib_i := \sum_i a^ib_i$.  
   
   \item Last, recall that a flow on E is a family $(\varphi_{ts})_{0\leq s\leq t\leq T}$ of maps from E to itself such that $\varphi_{tt}= \textrm{Id}$, for all $0\leq t\leq T$, and $\varphi_{tu}\circ\varphi_{us} = \varphi_{ts}$, for all $0\leq s\leq u\leq t\leq T$. The letter $T$, here and below, will always stand for a finite time horizon.  
\end{itemize}

%------------------------%
\section{Rough flows}
\label{SectionRoughFlows}
%------------------------%

%%---------------------------------------------%%
\subsection{Flows and approximate flows}
\label{SubsectionApproximateFlows}
%%---------------------------------------------%%

We introduced in \cite{B15} a simple machinery for constructing flows on E from approximate flows that can be understood as a generalization of Lyons' workhorse \cite{Lyo98} for constructing a rough path from an almost rough path; this is the core tool for the construction of the rough integral. Roughly speaking, the "Approximate flow-to-flow" machinery says that if we are given a family of maps $(\mu_{ts})_{0\leq s\leq t\leq T}$ from E to itself, and if the maps $\mu$ are close to defining a flow, in the sense that $\mu_{tu}\circ\mu_{us}-\mu_{ts}$ is small in a quantitative way, for $s\leq u\leq t$ with $t-s$ small, then there exists a unique flow close to $\mu$. In the rough paths setting, Lyons almost multiplicative functionals involve a family $a = (a_{ts})_{0\leq s\leq t\leq T}$ of elements of a tensor algebra such that $a_{tu}\,a_{us}-a_{ts}$ has some given size whenever $s\leq u\leq t$ with $t-s$ small, with the product on the tensor space used here. Despite their similarity, Lyons' setting differs 
from the present setting in that multiplication in a tensor algebra satisfies the distributivity property $ab-ac = a(b-c)$, which obviously does not hold if $a,b,c$ are maps and the product stands for composition. This seemingly minor point makes a real difference though, so it is fortunate that one can still get an analogue of Lyons' theorem in a function space setting. This comes at a little price on the regularity of the set of maps $\mu$ that one can consider. As usual, for $0<r\leq 1$, we denote by $\mcC^r$ the space of $r$-H\"older functions, with the understanding that they are Lipschitz continuous for $r=1$.

\medskip

\begin{definition}  {\sf 
\label{DefnC1ApproximateFlow}
Let $0<r\leq 1$ be given. A \emph{{\bf $\mcC^r$-approximate flow}} on $\textrm{E}$ is a family $\big(\mu_{ts}\big)_{0\leq s\leq t\leq T}$ of $(1+\rho)$-Lipschitz maps from $\textrm{E}$ to itself, for some $0<\rho\leq 1$, depending continuously on $(s,t)$ in the topology of uniform convergence and enjoying the following two properties.
\begin{itemize}
   \item {\bf Perturbation of the identity --} There exists a constant $\alpha$ with 
   $$
   0<1-\rho<\alpha<1,
   $$ 
   such that one has for any $0\leq s\leq t\leq T$, and any $x\in \textrm{E}$, the decomposition
 \begin{equation}
   \label{EqRegularityBounds}
   \quad D_x\mu_{ts} = \textrm{Id} + A^{ts}_x + B^{ts}_x,
   \end{equation} 
for some $\textrm{L}(\textrm{E})$-valued $\rho$-Lipschitz maps $A^{ts}$ with $\rho$-Lipschitz norm bounded above by $c|t-s|^\alpha$, and some $\textrm{L}(\textrm{E})$-valued $\mcC^1$ bounded maps $B^{ts}$, with $\mcC^1$-norm bounded above by $o_{t-s}(1)$. \vspace{0.1cm}

   \item {\bf $\mcC^r$-approximate flow property --} There exists some positive constants $c_1$ and $a>1$, such that one has
\begin{equation}
\label{EqMuMu}
\big\|\mu_{tu}\circ\mu_{us}-\mu_{ts}\big\|_{\mcC^r} \leq c_1\,|t-s|^a
\end{equation}
for all $0\leq s\leq u\leq t\leq T$.
\end{itemize}   }
\end{definition}

\medskip

Note that one requires a quantitative bound on $A$ while we only require a qualitative information on $B$, at the price of some more regularity for the latter. This fine decomposition of the differential of $\mu_{ts}$, as opposed to assuming only $D\mu_{ts} = \textrm{Id} + A^{ts}$, makes the setting more flexible. The introduction of the notion of approximate flow is justified by the following result proved in \cite{B15}. Given a partition $\pi_{ts}=\{s=s_0<s_1<\cdots<s_{n-1}<s_n=t\}$ of $(s,t)\subset [0,T]$, set 
$$
\mu_{\pi_{ts}} = \mu_{s_ns_{n-1}}\circ\cdots\circ\mu_{s_1s_0}.
$$

\medskip

\begin{thm}[Constructing flows on $\textrm{E}$]  {\sf 
\label{ThmConstructingFlows}  
A $\mcC^r$-approximate flow, with $\frac{1}{a}<r$, defines a unique flow $(\varphi_{ts})_{0\leq s\leq t\leq T}$ on $\textrm{E}$ to which one can associate a positive constant $\delta$ such that the inequality
$$
\big\|\varphi_{ts}-\mu_{ts}\big\|_\infty \leq |t-s|^a
$$
holds for all $0\leq s\leq t\leq T$ with $t-s\leq \delta$; this flow satisfies the inequality 
\begin{equation}
\label{EqApproxVarphiMu}
\big\|\varphi_{ts}-\mu_{\pi_{ts}}\big\|_\infty\leq 2\,c_1\,T\,\big|\pi_{ts}\big|^{ar - 1}
\end{equation} 
for any partition $\pi_{ts}$ of any interval $(s,t)\subset [0,T]$, with mesh $\big|\pi_{ts}\big|\leq \delta$. Moreover, the $\mcC^r$ norm of the maps $\varphi_{ts}$ is uniformly bounded by a function of the constant $c_1$ that appears in \eqref{EqMuMu}, for all $0\leq s\leq t\leq T$.   }
\end{thm}

\medskip

This statement generalises Gubinelli' sewing lemma \cite{Gub04}, such as reshaped by Feyel and de la Pradelle in \cite{FDlP06}, to the non-commutative, non-associative setting of maps on E. (The non-commutative sewing lemma of Feyel, de la Pradelle and Mokobodski \cite{FDlPM} requires associativity and cannot be used in the present setting.)

\medskip

Theorem \ref{ThmConstructingFlows} is stated in \cite{B15} for $\mcC^1$-approximate flows; the proof given there works verbatim for $\mcC^r$-approximate flows provided $\frac{1}{a}<r$; a $\mcC^1$ map is then understood in that setting as a Lipschitz map. The crucial point in the above statement is the fact that if $\mu$ depends continuously in $\mcC^r$ on some parameter then $\varphi$ also happens to depend continuously on that parameter, in $\mcC^0$, as a direct consequence of estimate \eqref{EqApproxVarphiMu}. As made clear in \cite{B15}, Theorem \ref{ThmConstructingFlows} can be seen as the cornerstone of the theory of rough differential equations, with the continuity of the It\^o-Lyons solution map given as a consequence of the aforementioned continuity of $\varphi$ on a parameter. We shall see in the present work that Theorem \ref{ThmConstructingFlows} is all we need to get back and extend the core results of the theory of stochastic flows intensively developed in the 80's and early 90's. We shall need 
for that purpose to introduce a notion of enriched velocity field that will somehow play the role in our setting of weak geometric H\"older $p$-rough paths, with $2\leq p<3$, in rough paths theory.

\medskip

\textit{We shall thus pick a regularity exponent $2\leq p<3$ here, once and for all}. Let us insist here that like for the theory of rough paths, the technical shape of the theory of rough drivers depends on that regularity exponent. Only two objects are needed in the definition of a rough driver when $2\leq p<3$; for $3\leq p<4$, we would need to introduce an additional object in the definition of a rough driver, that would thus consist of three objects; and so on. There is no real difficulty other than notational in giving a general theory, but as our applications of rough flows to the study of stochastic flows only require to develop the theory in the case where $2\leq p<3$, we stick to that setting and invite the reader to make up herself her mind about what the general theory looks like.

\bigskip

%%-------------------------%%
\subsection{Rough drivers}
\label{SubsectionRoughDrivers}
%%-------------------------%%

Let $\big(V(\cdot,t)\big)_{0\leq t\leq T}$ be a time-dependent vector field on E, with time increments 
$$
V_{ts}(\cdot) := V(\cdot,t)-V(\cdot,s).
$$ 
To get a hand on the definition of a weak geometric $p$-rough driver given below, think of $V_{ts}$ as given by the formula
\begin{equation}
\label{EqModelRoughDriver}
V_{ts} = V X_{ts},
\end{equation}
where $V(x) \in \textrm{L}(\RR^\ell,\RR^d)$, for all $x\in\RR^d$, and ${\bfX} = (X,\bbX)$ is a $p$-rough path over $\RR^\ell$. Write $V_i$ for the image by $V$ of the $i^{\textrm{th}}$ vector in the canonical basis of $\RR^\ell$. A solution path $x_\bullet$ to the rough differential equation 
$$
dx_t = V(x_t)\,{\bfX}(dt)
$$
can be characterized as a path satisfying some uniform Euler-Taylor expansion of the form
$$
f(x_t) = f(x_s) + X^i_{ts}\big(V_if\big)(x_s) + \bbX^{jk}_{ts}\big(V_jV_kf\big)(x_s) + O\Big(|t-s|^\frac{3}{p}\Big)
$$
for all sufficiently regular real-valued functions $f$ on $\RR^d$. The present Section will make it clear that the operators $X^i_{ts}V_i = VX_{ts}$ and $\bbX^{jk}_{ts}V_jV_k = (DV)V \,\bbX_{ts}$ are all we need in this formula to run the theory, with no need to separate their space part, given by $V$ and $(DV)V$, from their time part ${\bfX}_{ts}$.

\medskip

\begin{definition}  {\sf
Let $2\leq p<3$, and $p-2<\rho\leq 1$ be given. A weak geometric $(p,\rho)$-{\bf rough driver} is a family $\big({\bfV}_{ts}\big)_{0\leq s\leq t\leq T}$, with 
$$
{\bfV}_{ts} := \big(V_{ts},\VV_{ts}\big),
$$ 
and $\VV_{ts}$ a second order differential operator, such that \vspace{0.1cm}
\begin{itemize}
   \item[{\bf (i)}] the vector fields $V_{ts}$ are additive as functions of time
   $$
   V_{ts} = V_{tu} + V_{us}   
   $$
	for all $s<u<t$, and each $V_{ts}$ is of class $\mcC^{2+\rho}$ on E, with
$$
\underset{0\leq s<t\leq T}{\sup}\;\frac{\big\|V_{ts}\big\|_{\mcC^{2+\rho}}}{|t-s|^{\frac{1}{p}}} < \infty,
$$ \vspace{0.1cm}

   \item[{\bf (ii)}] the second order differential operators 
   $$
   W_{ts} := \VV_{ts} - \frac{1}{2}V_{ts}V_{ts}, 
   $$
   are actually vector fields, and
   $$
   \underset{0\leq s<t\leq T}{\sup}\;\frac{\big\|W_{ts}\big\|_{\mcC^{1+\rho}}}{|t-s|^{\frac{2}{p}}} < \infty,   
   $$ \vspace{0.1cm}

   \item[{\bf (iii)}] we have 
   $$
   \VV_{ts} = \VV_{tu} + V_{us}V_{tu} + \VV_{us}, 
   $$
   for any $0\leq s\leq u\leq t\leq T$.   
\end{itemize}    }
\end{definition}

\ssk

With in mind the model weak geometric $p$-rough driver given by formula \eqref{EqModelRoughDriver}, the requirement $p-2<\rho$ appears as a natural assumption to impose, given  known well-posedness results on rough differential equations \cite{D07}; the first order condition on the operators $W_{ts}$ justifies that we call $\bfV$ a \textit{weak geometric} $p$-rough driver, and condition {\bf (iii)} stands for an analogue of \textit{Chen's relation.} We shall freely talk about rough drivers rather than weak geometric $(p,\rho)$-rough drivers in the sequel. 

\ssk

\begin{definition}
We define the (pseudo-)norm of $\bfV$ to be
\begin{equation}
\label{EqDefnNorm}
\|{\bfV}\|_{p,\rho} := \underset{0\leq s< t\leq T}{\sup}\;\left\{\frac{\big\|V_{ts}\big\|_{\mcC^{2+\rho}}}{|t-s|^{\frac{1}{p}}} \vee \;\frac{\big\|W_{ts}\big\|_{\mcC^{1+\rho}}}{|t-s|^{\frac{2}{p}}}\right\}.
\end{equation}
and define an associated (pseudo-)metric on the set $\mcD_{p,\rho}$ of weak geometric $(p,\rho)$-rough drivers setting
$$
d_{p,\rho}(\bfV,\bfV') := \big\|\bfV-\bfV'\big\|_{p,\rho}.
$$
Like the space of rough paths the space of rough drivers is not a linear space. We will also need the homogeneous metric
\begin{align}\label{eqn:def_homog_metric}
 \mathfrak{d}_{p,\rho} (\bfV,\bfV') := \underset{0\leq s< t\leq T}{\sup}\;\left\{\frac{\big\|V_{ts} - V'_{ts} \big\|_{\mcC^{2+\rho}}}{|t-s|^{\frac{1}{p}}} \vee \; \sqrt{\frac{\big\|W_{ts} - W'_{ts}\big\|_{\mcC^{1+\rho}}}{|t-s|^{\frac{2}{p}}}}\right\}.
\end{align}
We will often drop the subindices $p,\rho$ when it is clear from the context in which space we are working in. Note that $d$ and $\mathfrak{d}$ induce the same topology in the space of rough drivers.

\end{definition}

\ssk

Note that one should add the $\mcC^{2+\rho}$-norm of $V(\cdot,0)$ in formula \eqref{EqDefnNorm} to define a proper norm on the space of rough drivers. This has no consequences as rough drivers are only used via their increments. Note also that given a rough driver $\bfV$ and $0<a\leq T$, one defines another rough driver ${\bfV}^a = \big(V^a,\VV^a\big)$, on the time interval $[0,a]$, setting
\begin{equation}
\label{EqReversedRoughVectorField}
\begin{split}
&V^a_{ts} = V_{a-s,a-t}, \\
&\VV^a_{ts} = -\VV_{a-s,a-t} + V_{a-s,a-t}V_{a-s,a-t},
\end{split}
\end{equation}
for all $0\leq s\leq t\leq a$. It is indeed elementary to check that these operators satisfy the algebraic conditions \textbf{(ii)} and \textbf{(iii)}, with 
\begin{equation}
\label{EqWr}
W^a_{ts} := \VV^a_{ts}-\frac{1}{2}V^a_{ts}V^a_{ts} = -W_{a-s,a-t};
\end{equation} 
that they satisfy the above analytic requirements is obvious. This rough driver is called the {\bf time reversal of the rough driver $\bfV$, from time $a$}. Note that $\big\|{\bfV}^a\big\|\leq\|{\bfV}\|$.

\bigskip

%%-----------------------%%
\subsection{Rough flows} 
\label{SubsectionRoughFlows}
%%-----------------------%%

We shall adopt below a definition of a solution flow to the equation
\begin{equation}
\label{EqRoughFlow}
d\varphi = {\bfV}(\varphi\,; dt)
\end{equation}
similar to the above definition of a solution path to a rough differential equation. A solution flow will be required to satisfy some uniform Euler-Taylor expansion of the form 
$$
f\circ\varphi_{ts} - \big\{ f + V_{ts}f + \VV_{ts}f \big\} = O\Big(|t-s|^\frac{3}{p}\Big),
$$
for all sufficiently regular real-valued functions $f$ on $\RR^d$. It is actually elementary to construct a family of maps $\big(\mu_{ts}\big)_{0\leq s\leq t\leq T}$ which enjoys the above Euler-Taylor expansion property. The key point is that this family will turn out to be a $\mcC^\rho$-approximate flow, if $\rho$ is not too small, so we shall get the existence and uniqueness of a solution flow from its very definition and Theorem \ref{ThmConstructingFlows}. 

\medskip

Given $0\leq s\leq t\leq T$, consider the ordinary differential equation 
\begin{equation}
\label{ApproximateRDE}
\dot y_u = V_{ts}(y_u) + W_{ts}(y_u), \quad 0\leq u\leq 1,
\end{equation}
and denote by $\mu_{ts}$ its well-defined time $1$ map, associating to any $x\in \textrm{E}$ the value at time $1$ of the solution of equation \eqref{ApproximateRDE} started from $x$. Elementary results on ordinary differential equations imply that if one considers $y_u$ as a function of $x$ , for $0\leq u\leq 1$, then we have
\begin{equation}
\label{EqBasicEstimate}
\big\|y_u - \textrm{Id}\big\|_{\mcC^1} \leq c\,\|{\bfV}\|\,|t-s|^{\frac{1}{p}}
\end{equation}
for some universal positive constant $c$. Proposition \ref{PropFundamentalEstimate} below shows that the maps $\mu_{ts}$ have the awaited Euler-Taylor expansion expected from a solution flow to equation \eqref{EqRoughFlow}.

\medskip

\begin{prop}  {\sf
\label{PropFundamentalEstimate}
We have 
\begin{equation}
\label{EqFundamentalEstimate}
\Big\| f\circ\mu_{ts} - \big\{f+V_{ts}f+\VV_{ts}f\big\}\Big\|_\infty \leq c\,\|f\|_{\mcC^{2+\rho}}|t-s|^{\frac{3}{p}},
\end{equation}
for any $f\in\mcC^{2+\rho}$ and any $0\leq s\leq t\leq T$.   }
\end{prop}

\medskip

The proof of this statement is straightforward and relies on the the following formula. For all $x\in \textrm{E}$ and all $f\in\mcC^2$, we have
\begin{equation*}
\begin{split}
f\big(\mu_{ts}(x)\big) &= f(x) + \int_0^1 \big(V_{ts}f\big)(y_u)\,du + \int_0^1 \big(W_{ts}f\big)(y_u)\,du \\
                                     &= f(x) + \big(V_{ts}f\big)(x) + \big(\VV_{ts}f\big)(x) + \ep^f_{ts}(x)
\end{split}
\end{equation*}
where
\begin{equation*}
\begin{split}
\ep^f_{ts}(x) :=  \int_0^1\int_0^u \Big\{\big(V_{ts}V_{ts}f\big)(y_r) - \big(V_{ts}V_{ts}f\big)(x)\Big\}\,drdu &+ \int_0^1\int_0^u \big(W_{ts}V_{ts}f\big)(y_r)\,drdu \\ &+\int_0^1\big\{\big(W_{ts}f\big)(y_u)-\big(W_{ts}f\big)(x)\big\}\,du.
\end{split}
\end{equation*}
The inequality
\begin{equation}
\label{EqEstimateEpFTs}
\Big\|\ep^f_{ts}\Big\|_{\mcC^\rho} \leq c\,\big(1+\|{\bfV}\|^3\big)\|f\|_{\mcC^{2+\rho}}\,|t-s|^{\frac{3}{p}},
\end{equation}
justifies Proposition \ref{PropFundamentalEstimate}. 

\medskip

\begin{thm} {\sf 
\label{ThmC1ApproximateFlow}
The family of maps $\big(\mu_{ts}\big)_{0\leq s\leq t\leq T}$ is a $\mcC^\rho$-approximate flow which depends continuously on $\big((s,t),{\bfV}\big)$ in $\mcC^0$ topology.   }
\end{thm}

\medskip

\begin{Dem}
The family $\big(\mu_{ts}\big)_{0\leq s\leq t\leq T}$ satisfies the regularity assumptions \eqref{EqRegularityBounds} as a direct consequence of classical results on the dependence of solutions to ordinary differential equations with respect to parameters, including the initial condition for the equation. These results also imply the continuous dependence of $\mu_{ts}$ on $\big((s,t),{\bfV}\big)$ in $\mcC^\rho$ topology. To show that the family $\mu$ defines a $\mcC^\rho$-approximate flow, write, for $0\leq s\leq u\leq t\leq T$,
\begin{equation*}
\begin{split}
\mu_{tu}\big(\mu_{us}(x)\big) &= \mu_{us}(x) + V_{tu}\big(\mu_{us}(x)\big) + \big(\VV_{tu}\textrm{Id}\big)\big(\mu_{us}(x)\big) + \ep^{\textrm{Id}}_{tu}\big(\mu_{us}(x)\big) \\
  												&= x\,+V_{us}(x)+\big(\VV_{us}\textrm{Id}\big)(x) + \ep^{\textrm{Id}}_{us}(x) \\
  												&\quad\quad\, +V_{tu}(x) + \big(V_{us}V_{tu}\big)(x) + \big(\VV_{us}V_{tu}\big)(x) + \ep^{V_{tu}}_{us}(x) \\
  												&\quad\quad\, +\big(\VV_{tu}\textrm{Id}\big)(x) + \Big(\big(\VV_{tu}\textrm{Id}\big)\big(\mu_{us}(x)\big) - \big(\VV_{tu}\textrm{Id}\big)(x)\Big) + \ep^{\textrm{Id}}_{tu}\big(\mu_{us}(x)\big) \\
  												&= \mu_{ts}(x) + \Big\{\big(\VV_{us}V_{tu}\big)(x) + \Big(\big(\VV_{tu}\textrm{Id}\big)\big(\mu_{us}(x)\big) - \big(\VV_{tu}\textrm{Id}\big)(x)\Big) + \ep^{V_{tu}}_{us}(x) \\
  												&\quad\quad\, +\ep^{\textrm{Id}}_{us}(x) + \ep^{\textrm{Id}}_{tu}\big(\mu_{us}(x)\big)\Big\}.
\end{split}
\end{equation*}
The approximate flow property then follows from the regularity assumptions on $V_{ts}$ and $\VV_{ts}$, and estimate \eqref{EqEstimateEpFTs}.
\end{Dem}

\medskip

With the notations used in the definition of an approximate flow, the exponent $a$ that appears here in the approximate flow identity \eqref{EqMuMu} is $a=\frac{3}{p}$.

\medskip

\begin{defn}   {\sf
A \textbf{flow} $\big(\varphi_{ts}\big)_{0\leq s\leq t\leq T}$ is said to \textbf{solve the rough differential equation}
\begin{equation}
\label{RDE}
d\varphi = {\bfV}(\varphi \,; dt)
\end{equation}
if there exists a possibly $\bfV$-dependent positive constant $\delta$ such that the inequality
\begin{equation}
\label{DefnSolRDE}
\big\|\varphi_{ts}-\mu_{ts}\big\|_\infty \leq |t-s|^\frac{3}{p}
\end{equation}
holds for all $0\leq s\leq t\leq T$ with $t-s\leq\delta$. Flows solving a differential equation of the form \eqref{RDE} are called {\bf rough flows}. If equation \eqref{RDE} is well-posed, the map which associates to a rough driver $\bfV$ the solution flow to equation \eqref{RDE} is called the {\bf It\^o map}.  }
\end{defn}

\medskip

Following Cass and Weidner \cite{CassWeidner}, one can equivalently take the Taylor expansion property
$$
\varphi_{ts} = \textrm{Id} + V_{ts}\textrm{Id} + \VV_{ts}\textrm{Id} + O\big(|t-s|^\frac{3}{p}\big)
$$ 
as a defining property of a solution flow to the rough differential equation \eqref{RDE}. The following well-posedness result comes as a direct consequence of Theorem \ref{ThmConstructingFlows} and Theorem \ref{ThmC1ApproximateFlow}. A family of maps is said to be uniformly $\mcC^\rho$ is it has uniformly bounded $\mcC^\rho$-norm.

\medskip

\begin{thm}  {\sf 
\label{ThmMain}
Assume $\rho > \frac{p}{3}$. Then the differential equation on flows
$$
d\varphi = {\bfV}(\varphi \,; dt)
$$ 
has a unique solution flow; it takes values in the space of uniformly $\mcC^\rho$ homeomorphisms of $\textrm{E}$, with uniformly $\mcC^\rho$ inverses, and depends continuously on $\bfV$ in the topology of uniform convergence.   }
\end{thm}

\medskip

\begin{Dem}
It follows from the proof of Theorem \ref{ThmC1ApproximateFlow} that one can choose as a constant $c_1$ in inequality \eqref{EqMuMu} a multiple of $1+\|\bfV\|^4$, so we have from Theorem \ref{ThmConstructingFlows} the estimate
$$
\big\|\varphi_{ts}-\mu_{\pi_{ts}}\big\|_\infty \leq c\,\big(1+\|{\bf V}\|^4\big)\,T\,\big|\pi_{ts}\big|^{\rho\,\frac{3}{p}-1},
$$ 
for any partition $\pi_{ts}$ of $(s,t)\subset [0,T]$ with mesh $\big|\pi_{ts}\big|$ small enough, say no greater than $\delta$. Note that the exponent $\rho\,\frac{3}{p}-1$ is positive. As these bounds are uniform in $(s,t)$, and for $\bfV$ in a bounded set of the space of rough drivers, and each $\mu_{\pi_{ts}}$ is a continuous function of $\bfV$, the flow $\varphi$ depends continuously on $\big((s,t),\bfV\big)$.

\ssk

To prove that $\varphi$ is a homeomorphism, note that it follows from \eqref{EqWr} that, for $0\leq a\leq b\leq t$, each $\mu_{ba}$ is a diffeomorphism with inverse given by the time one map $\mu^t_{t-a,t-b}$ of the ordinary differential equation
$$
\dot y_u = -V_{ba}(y_u) - W_{ba}(y_u) = V^t_{t-a,t-b}(y_u) + W^t_{t-a,t-b}(y_u),\quad 0\leq u\leq 1,
$$
associated with the time reversed rough driver ${\bfV}^t$. As $\mu^t$ has the same properties as $\mu$, the maps 
$$
\big(\mu_{\pi_{ts}}\big)^{-1} = \mu_{s_1s_0}^{-1}\circ\cdots\circ\mu_{s_n s_{n-1}}^{-1} = \mu^t_{s_n s_{n-1}}\circ\cdots\circ \mu^t_{s_1s_0}
$$ 
converge uniformly to some continuous map $\varphi_{ts}^{-1}$,as the mesh of the partition $\pi_{ts}$ tends to $0$; these limit maps $\varphi_{ts}^{-1}$ satisfy by construction $\varphi_{ts}\circ\varphi_{ts}^{-1} = \textrm{Id}$.

\ssk

As Theorem \ref{ThmConstructingFlows} provides a uniform control of the $\mcC^\rho$ norm of the maps $\varphi_{ts}$, the same control holds for their inverses since $\big\|{\bfV}^t\big\|\leq \|{\bfV}\|$. We propagate this control from the set $\big\{(s,t)\in [0,T]^2\,;\,s\leq t, \; t-s\leq \delta\big\}$ to the whole of $\big\{(s,t)\in [0,T]^2\,;\,s\leq t\big\}$ using the flow property of $\varphi$.
\end{Dem}

\medskip
 
Note that the solution flow to the rough differential equation 
$$
d\psi = {\bfV}^T(\psi \,; dt),
$$ 
driven by the time reversal of the rough driver $\bfV$, from time $T$, provides the inverse flow of $\varphi$, in the sense that 
$$
\varphi_{ts}^{-1} = \psi_{T-s.T-t},
$$
for all $0\leq s\leq t\leq T$. Last, note that it is elementary to adapt the above results to add a globally Lipschitz drift in the dynamics; the above results hold in that setting as well. 

\bigskip

\noindent {\large{\sc Remark{\bf .}}} \textit{The rough drivers introduced here are somewhat a dual version of similar objects that were introduced very recently in \cite{BailleulGubinelli} by one of the authors and Gubinelli in the study of the well-posedness of a general family of linear hyperbolic symmetric systems of equations driven by time-dependent vector fields that are only distributions in the time direction. The latter work deals with evolutions in function spaces and uses functional analytic tools in the setting of controlled paths to make a first step towards a general theory of rough PDEs, in the lines of the classical PDE approach based on duality, a priori estimates and compactness results. The present work does not overlap with the latter.
}  

\bigskip

\textsf{\textbf{How the story goes on.}} The entire technical core of the theory of rough flows is contained in Section \ref{SubsectionApproximateFlows} and Section \ref{SubsectionRoughFlows}. The remainder of this work is dedicated to   \vspace{0.2cm}

\begin{itemize}
   \item showing how one can lift semimartingale velocity fields into rough drivers -- Section \ref{SubsectionKolmogorovTheorem} and \ref{SubsectionMartingaleDrivers},   \vspace{0.15cm}
   
   \item showing that stochastic and rough flows concide for semimartingale velocity fields -- Section \ref{SubsectionStochRoughFlows},   \vspace{0.15cm}
   
   \item probing a Wong-Zaka\" i-type support theorem for semimartingale stochastic flows -- Section \ref{SubsectionStrongApproximations},   \vspace{0.15cm}
   
   \item proving some sharp Schilder and large deviation theorems for flows generated by Gaussian rough drivers -- Section \ref{SubsectionLDP}.
\end{itemize}

\medskip

We emphasize here that we proved in the subsequent work \cite{BailleulRiedelScheutzow} that one can lift into rough drivers a whole class of Gaussian velocity fields, showing that the setting of rough flows goes beyond the setting of semimartingale calculus.

\bigskip

%%--------------------------------------------------%%
\subsection{An It\^o formula for rough flows}
\label{SubsectionItoFormula}
%%--------------------------------------------------%%

With a view to identifying stochastic and rough flows in Section \ref{SubsectionStochRoughFlows}, we prove here an elementary It\^o formula analogue to Friz and Hairer's It\^o formula in \cite{FH14}. As a matter of fact, Theorem \ref{ThmItoFormula} below states that any $\frac{1}{p}$-H\"older path in a Banach space satisfies an It\^o formula, outside the setting of rough or controlled paths. To state and prove it recall Feyel and de la Pradelle sewing lemma \cite{FDlP06}, that can be seen as a precursor of the construction theorem for flows, Theorem \ref{ThmConstructingFlows}. Given a partition $\pi_{ts} = \{s=s_0<s_1<\cdots<s_{n-1}<s_n=t\}$ of an interval $[s,t]$, and an $\textrm{E}$-valued $2$-index map $z = (z_{ts})_{0\leq s\leq t\leq T}$, set 
$$
z_{\pi_{ts}} := z_{s_ns_{n-1}} + \cdots + z_{s_1s_0}.
$$

\medskip

\begin{thm}[\cite{FDlP06}]  {\sf
\label{ThmSewingLemmaFdlP}
Let $\big(z_{ts}\big)_{0\leq s\leq t\leq T}$ be an $\textrm{E}$-valued $2$-index continuous map to which one can associate some positive constants $c_1$ and $a>1$ such that 
\begin{equation}
\label{EqAlmostAdditivity}
\big|\big(z_{tu}+z_{us}\big)-z_{ts}\big| \leq c_1\,|t-s|^a
\end{equation}
holds for all $0\leq s\leq u\leq t\leq T$. Then there exists a unique continuous function $Z : [0,T]\rightarrow \RR$, with increments $Z_{ts} := Z_t-Z_s$,to which one can  associate a positive constant $\delta$ such that the inequality
$$
\big|Z_{ts}-z_{ts}\big| \leq |t-s|^a,
$$
holds for all $0\leq s\leq u\leq t\leq T$, with $t-s\leq \delta$; this map $Z$ satisfies the inequality
$$
\big| Z_{ts} - z_{\pi_{ts}}\big| \leq 2\,c_1\,T\,\big|\pi_{ts}\big|^{a-1}
$$
for any partition $\pi_{ts}$ of any interval $[s,t]\subset [0,T]$, with mesh $\big|\pi_{ts}\big|\leq \delta$. It follows in particular that $Z$ depends continuously on any parameter in uniform topology if $z$ does.   }
\end{thm}

\medskip

\noindent A map $z$ satisfying condition \eqref{EqAlmostAdditivity} is said to be {\bf almost-additive}, and we write 
$$
Z_{ts} =: \int_s^t z_{du}.
$$
We equip the tensor product space $\textrm{E}\otimes \textrm{E}$ with a compatible tensor norm that makes the natural embedding $\textrm{L}\big(\textrm{E},\textrm{L}(\textrm{E},\RR)\big) \subset \textrm{L}\big(\textrm{E}\otimes \textrm{E},\RR\big)$ continuous. Given such an choice, one can identify the second differential of a $\mcC^2$ real-valued function on $\textrm{E}$ to an element of $\textrm{L}\big(\textrm{E}\otimes\textrm{E},\RR\big)$ that is symmetric; this is what we do below.

\medskip

\begin{thm}[It\^o formula]  {\sf
\label{ThmItoFormula}
Let $\textrm{F} : [0,T]\times \textrm{E}\rightarrow \RR$ be a $\mcC^1$-function of time with time derivative $\partial_t \textrm{F}(t,x)$ bounded and continuous, uniformly in $x\in\textrm{E}$. Assume also that F is of class $\mcC^3$ in the sense of Fr\'echet as a function of its $\textrm{E}$-component, with derivatives $\textrm{F}^{(1)}, \textrm{F}^{(2)}, \textrm{F}^{(3)}$ and $\partial_t\textrm{F}^{(1)}$, bounded uniformly in time. Let $(x_t)_{0\leq s\leq t\leq T}$ be $\frac{1}{p}$-H\"older $\textrm{E}$-valued map. Then the continuous $2$-index map 
$$
z_{ts} := \textrm{F}^{(1)}_{(s,x_s)}\,(x_t-x_s) + \frac{1}{2}\textrm{F}^{(2)}_{(s,x_s)}\,(x_t-x_s)^{\otimes 2}
$$
is almost-additive, and we have 
\begin{equation}
\label{EqItoFormula}
\textrm{F}\big(t,x_t\big) = \textrm{F}\big(s,x_s\big) + \int_s^t \big(\partial_r\textrm{F}\big)(r,x_r)\,dr + \int_s^t z_{du},
\end{equation}
for any $0\leq s\leq t\leq T$.  }
\end{thm}

\medskip

\begin{Dem}
The proof is a straightforward application of Feyel-de la Pradelle' sewing lemma, Theorem \ref{ThmSewingLemmaFdlP}. Given $0\leq s\leq u\leq t\leq T$, the algebraic identity
\begin{equation*}
\begin{split}
z_{tu}+z_{us} &= \textrm{F}^{(1)}_{(s,x_s)}\,(x_t-x_s) + \Big(\textrm{F}^{(1)}_{(u,x_u)}-\textrm{F}^{(1)}_{(s,x_s)}\Big)(x_t-x_u) \\
					     &\quad+ \frac{1}{2}\textrm{F}^{(2)}_{(u,x_u)}\,(x_t-x_u)^{\otimes 2} + \frac{1}{2} \textrm{F}^{(2)}_{(s,x_s)} \,(x_u-x_s)^{\otimes 2},
\end{split}
\end{equation*}
the regularity assumptions on F and the symmetric character of $\textrm{F}^{(2)}_{(s,x)}$, for any $x\in\textrm{E}$, we have
\begin{equation*}
\begin{split}
z_{tu}+z_{us} &= \textrm{F}^{(1)}_{(s,x_s)}\big(x_t-x_s\big)  + \textrm{F}^{(2)}_{(s,x_s)}\,(x_u-x_s)\otimes (x_t-x_u)   \\
& \quad+ O\Big(|t-s|^\frac{3}{p}\Big) + O\big(\|x_u-x_s\|^2\big)\,\big\|x_t-x_u\big\|   \\
&\quad+ \frac{1}{2}\textrm{F}^{(2)}_{(u,x_u)}\,(x_t-x_u)^{\otimes 2} + \frac{1}{2} \textrm{F}^{(2)}_{(s,x_s)}\,  (x_u-x_s)^{\otimes 2} \\
&= z_{ts} + O\Big(|t-s|^\frac{3}{p}\Big).
\end{split}
\end{equation*}
It\^o's formula \eqref{EqItoFormula} follows by noting that we have for all $n\geq 1$
\begin{equation*}
\begin{split}
\textrm{F}\big(t,x_t\big) &= \sum_{i=0}^{n-1} \Big\{\textrm{F}\big(s_{i+1},x_{s_{i+1}}\big) - F\big(s_i,x_{s_i}\big)\Big\} \\
							       &= o_n(1) + \sum_{i=0}^{n-1} (s_{i+1}-s_i)\,(\partial_s \textrm{F})\big(s_i,x_{s_i}\big) + \sum_{i=0}^{n-1} \Big\{\textrm{F}\big(s_i,x_{s_{i+1}}\big) - F\big(s_i,x_{s_i}\big)\Big\}, 
\end{split}
\end{equation*}
with
\begin{equation*}
\begin{split}
\textrm{F}\big(s_i, x_{s_{i+1}}\big) &- \textrm{F}\big(s_i, x_{s_i}\big)   \\
&= \textrm{F}^{(1)}_{(s_i,x_{s_i})}\big(x_{s_{i+1}}-x_{s_i}\big) + \frac{1}{2}\textrm{F}^{(2)}_{(s_ix_{s_i})}\,(x_{s_{i+1}}-x_{s_i})^{\otimes 2} + O\Big(\big| x_{s_{i+1}}-x_{s_i}\big|^3\Big) \\
&= z_{s_{i+1}s_i} + O\Big(|s_{i+1}-s_i|^\frac{3}{p}\Big).
\end{split}
\end{equation*}
\end{Dem}

\bigskip

As an example, consider the solution flow $\varphi$ to a rough differential equation on $\RR^d$
$$
d\varphi = {\bfV}(\varphi\,; dt).
$$
Write $\varphi_t$ for $\varphi_{t0}$, and consider it as an element of the space $\textrm{E}$ of continuous paths from $[0,T]$ to $\mcC\big(\RR^d,\RR^d\big)$, equipped with the norm of uniform convergence, with $\mcC\big(\RR^d,\RR^d\big)$ endowed with a norm inducing uniform convergence on compact sets. It satisfies by its very definition and Proposition \ref{PropFundamentalEstimate} the Euler-Taylor expansion
$$
\varphi_t = \varphi_s + \big(V_{ts}\textrm{Id}\big)\circ\varphi_s + \big(\VV_{ts}\textrm{Id}\big)\circ\varphi_s + O\big(|t-s|^\frac{3}{p}\big)
$$
so it is a $\frac{1}{p}$-H\"older path in that space. Now, given some points $y_1,\dots, y_k$ in $\RR^d$ and a $\mcC^3_b$ real-valued function $f$ on $(\RR^d)^k$, one can think of the function 
\begin{equation}
\label{EqExampleFunctioF}
\textrm{F}(\phi) = f\big(\phi(y_1),\dots,\phi(y_k)\big),
\end{equation}
for $\phi\in \textrm{E}$, as a typical time-independent example of function satisfying the conditions of Theorem \ref{ThmItoFormula}. One then has
\begin{equation*}
\begin{split}
\textrm{F}\big(\varphi_{s_{i+1}}\big) - \textrm{F}\big(\varphi_{s_i}\big) &= f\Big(\varphi_{s_{i+1}s_i}\big(\varphi_{s_i}(y_1)\big),\dots,\varphi_{s_{i+1}s_i}\big(\varphi_{s_i}(y_k)\big)\Big) - f\big(\varphi_{s_i}(y_1),\dots,\varphi_{s_i}(y_k)\big) \\
&= \sum_{m=1}^k \Big(\big(V^{\{m\}}_{s_{i+1}s_i}+\VV^{\{m\}}_{s_{i+1}s_i}\big)f\Big)\big(\varphi_{s_i}(y_1),\dots,\varphi_{s_i}(y_k)\big) + O_c\Big(|s_{i+1}-s_i|^\frac{3}{p}\Big),
\end{split}
\end{equation*}
where the upper index $\{m\}$ means that the operators act on the $m^\textrm{th}$ component of $f$. The above sum defines an almost-additive continuous function $z^f_{ts}$, taken here at time $\big(s_{i+1},s_i\big)$, so we have
$$
f\big(\varphi_t(y_1),\dots,\varphi_t(y_k)\big) = f\big(\varphi_s(y_1),\dots,\varphi_s(y_k)\big) + \int_s^t z^f_{du}
$$
for all times $0\leq s\leq t\leq T$.

\bigskip

%%----------------------------------------------------------%%
\subsection{A Kolmogorov-type regularity theorem}
\label{SubsectionKolmogorovTheorem}
%%----------------------------------------------------------%%

We shall use below the theory of rough drivers in a setting where the drivers are random. Like in the theory of rough paths, the primary object we are given is not the random rough driver itself, or the random rough path, but rather a genuine random vector field, or random path, which needs to be enhanced in a first step into a random rough driver, or random rough path. This first, purely probabilistic, step can typically be done using some Kolmogorov-type continuity arguments. We give in this Section some variations on this theme that will be used to enhance vector field-valued martingales into rough drivers in Section \ref{SubsectionMartingaleDrivers}; \textit{a reader interested only in these applications is advised to skip the technical details and only have a look at Theorem} \ref{thm:kolmogorov_rough_drivers}; for the other readers, it is our hope that this somewhat long section contains some material interesting in itself; it provides moment conditions under which one can get back uniform in time 
estimates on quantities of the form $(t-s)^{-\alpha}\|X_{ts}\|_{\mathcal{C}^a}$, such as required by the definition of a rough driver.

\medskip

The next Lemma gives sufficient conditions for a process defined on a possibly unbounded domain to be bounded. Recall the equivalence of having Gaussian tails to square-root growth of moments, cf. \cite[Lemma A.17]{FV10}.

\medskip

\begin{lem}\label{lem:cond_proc_bdd} {\sf
 Let $(E,d)$ be a complete, separable metric space. Let $D$ be an open subset of $\RR^d$, $X \colon D \to (E,d)$ a continuous stochastic process, $e \in E$ and $\kappa > 0$. Set
 \begin{align*}
  D_n := \big\{x \in D \, :\, n-1 \leq |x| < n \big\}
 \end{align*}
 and $\mathbf{N} := \{n \in \NN\, :\, D_n \neq \emptyset\}$.
Let $(a_n)_{n \in \mathbf{N}}$ be a sequence of non-negative real numbers and $(x_n)_{n \in \mathbf{N}}$ a sequence of elements in $D$ such that $x_n \in D_n$ for every $n \in \mathbf{N}$.
 
\begin{itemize}
  \item[(i)] Assume that there is a $q \in [1, \infty)$ and a $\gamma \in (\frac{d}{q},1]$ such that
 \begin{align*}
  \sup_{x,y \in D_n} \Big\| d\big(X(x),X(y)\big) \Big\|_{L^q} \leq \kappa a_n |x - y|^{\gamma}
 \end{align*} 
and that
 \begin{align*}
  \big\| d(X(x_n),e) \big\|_{L^q} \leq \kappa a_n
 \end{align*}
for every $n \in \mathbf{N}$. Set $(b_n) := (a_n n^{\gamma})$ and assume that $\|b \|_{l^q} \leq K < \infty$. Then there is a constant $C = C(q,\gamma)$ such that
 \begin{align*}
  \big\| \sup_{x \in D} d(X(x),e) \big\|_{L^{q}} \leq C K \kappa.
 \end{align*}
 
  \item[(ii)] Let $\gamma \in (0,1]$ and assume that for every $q \geq 1$ there is a $c_q$ such that for every $n \in \mathbf{N}$,
 \begin{align*}
  \sup_{x,y \in D_n} \Big\| d\big(X(x),X(y)\big) \Big\|_{L^q} \leq \kappa c_q a_n |x - y|^{\gamma}
 \end{align*} 
 and that
 \begin{align*}
  \big\| d(X(x_n),e) \big\|_{L^q} \leq \kappa c_q a_n
 \end{align*}
 where $c_q = \mathcal{O}(\sqrt{q})$ when $q \to \infty$. Assume that $a_n = \mathcal{O}\Big(n^{-\gamma} \big(1 + \log(n)\big)^{-\frac{1}{2}}\Big)$. Then for every $q \geq 1$ there is some constant $C = C(q,\gamma)$ such that
 \begin{align*}
  \Big\| \sup_{x \in D} d\big(X(x),e\big) \Big\|_{L^{q}} \leq C \kappa
 \end{align*}
  with $C = \mathcal{O}(\sqrt{q})$ when $q \to \infty$. In particular, the random variable $\sup_{x \in D} d(X(x),e)$ has Gaussian tails.
\end{itemize}   }
\end{lem}

\medskip

\begin{Dem}
Without loss of generality, one may assume $\kappa = 1$, otherwise we consider the metric $\tilde{d} = d/\kappa$ instead, and $\mathbf{N} = \NN$ -- otherwise we add small, disjoint balls to $D$ and define $X$ to be constant and equal to $e$ on these balls. We first prove claim (i). 

\ssk

Let $\alpha > \frac{d}{q}$ and set $p(u) = u^{\alpha + \frac{d}{q}}$. By the Garsia-Rodemich-Rumsey Lemma (cf. e.g. \cite[Lemma 2.4 (i)]{Sch09}), for every $x, y \in D_n$,
 \begin{align*}
  \frac{d\big(X(x),X(y)\big)}{|x - y|^{\alpha - \frac{d}{q}}} \leq C V_n^\frac{1}{q}
 \end{align*}
 where
 \begin{align*}
  V_n = \int_{D_n \times D_n} \frac{\Big| d\big(X(u),X(v)\big) \Big|^q}{|u-v|^{\alpha q + d}}\, du\, dv.
 \end{align*}
 Thus, by a change of variables,
 \begin{align*}
  \EE \left| \sup_{x,y \in D_n}  \frac{d\big(X(x),X(y)\big)}{|x - y|^{\alpha - \frac{d}{q}}}\right|^q &\leq C^q a^q_n \int_{D_n \times D_n}  |u-v|^{(\gamma - \alpha)q - d} \, du\, dv \\
  &\leq C^q a^q_n n^{d + (\gamma - \alpha)q} \int_{(0,1)^2}  |u-v|^{(\gamma - \alpha)q - d} \, du\, dv.
 \end{align*}
 Let $\beta \in (0, \gamma - \frac{d}{q})$ and set $\alpha = \frac{d}{q} + \beta < \gamma$. Then the integral is finite, and we have shown that
 \begin{align*}
  \left\| \sup_{x,y \in D_n}  \frac{d\big(X(x),X(y)\big)}{|x - y|^{\beta}}\right\|_{L^q} \leq C a_n n^{(\gamma - \beta)}.
 \end{align*}
 By the triangle inequality,
 \begin{align*}
  \left\| \sup_{x \in D_n} d\big(X(x),e\big) \right\|_{L^q} \leq C a_n (n^{\gamma} + 1) \leq 2C b_n.
 \end{align*}
 Thus we obtain
 \begin{align*}
  \EE \left( \sup_{x \in D} d\big(X(x),e\big)^q \right) &= \EE \left( \sup_n \sup_{x \in D_n} d\big(X(x),e\big)^q \right) \leq \sum_{n = 1}^{\infty} \EE \left( \sup_{x \in D_n} d\big(X(x),e\big)^q \right) \\
  &\leq 2^q C^q \sum_{n = 1}^{\infty} b_n^q < \infty
 \end{align*}
 and claim (i) is shown.

%  therefore
%  \begin{align*}
%   \PP \left( \sup_{x \in D_n} d\big(X(x),e\big) \geq t \right) \leq \frac{C^q b_n^q}{t^q} 
%  \end{align*}
%  and
%  \begin{align*}
%   \PP \left( \sup_{x \in D} d\big(X(x),e\big) \geq t \right) \leq \frac{C^q}{t^q} \sum_{n = 1}^{\infty} b_n^q = \frac{C^q \| b\|_{l^q}^q}{t^q} .
%  \end{align*}
%  Hence
%  \begin{align*}
%   \EE \left| \sup_{x \in D} d\big(X(x),e\big) \right|^{q'} \leq 1 + C^q \| b\|_{l^q}^q \int_1^{\infty} \frac{dt}{t^{q/q'}} \leq 1 + \frac{C^q \| b\|_{l^q}^q}{1 - \frac{q}{q'}}.
%  \end{align*}
 
\ssk

Now we prove claim (ii). Note that the constant in the Garsia-Rodemich-Rumsey Lemma may be chosen non-increasing in $q$. Therefore, we can argue similarly as before to see that for every $q \geq 1$ and $n \in \NN$,
 \begin{align*}
  \left\| \sup_{x \in D_n} d\big(X(x),e\big) \right\|_{L^q} \leq C_q b_n
 \end{align*}
 where $C_q = \mathcal{O}(\sqrt{q})$. This shows that the random variable has Gaussian tails, i.e. there is some constant $C$ such that for every $n \in \NN$
 \begin{align*}
   \PP \left( \sup_{x \in D_n} d\big(X(x),e\big) \geq t \right) \leq C \exp \left( - \frac{t^2}{C b_n^2} \right)
 \end{align*}
 for every $t \geq 0$. Hence
 \begin{align*}
  \PP \left( \sup_{x \in D} d\big(X(x),e\big) \geq t \right) &\leq C \sum_{n = 1}^{\infty} \exp \left( - \frac{t^2}{C b_n^2} \right) \leq C \sum_{n = 1}^{\infty} \exp \left( - \frac{t^2}{C} (1 + \log(n)) \right) \\
  &\leq C \exp \left( - \frac{t^2}{C} \right) \sum_{n = 1}^{\infty} n^{-\frac{t^2}{C}}
 \end{align*}
 and the sum is finite for $t$ large enough. This proves that $\sup_{x \in D} d(X(x),e)$ has Gaussian tails.
 % which is equivalent to say that its $L^q$ norm grows at most like $\sqrt{q}$ when $q \to \infty$.
\end{Dem}

\ssk

\begin{cor}\label{cor:cond_proc_bdd}	{\sf
 Let $D$ be an open subset of $\RR^d$, $(\textrm{E},\| \cdot \|)$ a separable Banach space, $X \colon D \to \textrm{E}$ a continuous stochastic process and $q > 1$.
 \begin{itemize}
  \item[(i)] Assume that there is a constant $\kappa > 0$ and $\gamma \in (0,1]$ such that for every $x,y \in D$ with $0 < |x - y| \leq 1$,
 \begin{align}\label{eqn:diff_est_stoch_proc}
    \big\|X(x) - X(y) \big\|_{L^q} \leq \kappa |x - y|^{\gamma}
 \end{align}
 % If $X$ has compact support, assume that $q > d/ \gamma$. If $X$ does not have compact support, assume 
 and that there is an $\eta \in (0,\infty)$ such that for every $x \in D$,
 \begin{align}\label{eqn:moment_growths_stoch_proc}
  \|X(x) \|_{L^q} \leq \frac{\kappa}{1 + |x|^{\eta}}
 \end{align}
  with $q$ sufficiently large satisfying
  \begin{align*}
   q > \frac{1}{\eta} + d\left( \frac{1}{\eta} + \frac{1}{\gamma} \right).
  \end{align*}
  Then the random variable $\sup_{x \in D} \| X(x) \|$ is almost surely finite. Moreover, %if $X$ has compact support, 
  there is a constant $C = C(\gamma,q,\eta)$ such that
  \begin{align}\label{eqn:moment_sup_norm_bound}
   \left\| \sup_{x \in D} \big\|X(x)\big\| \right\|_{L^{q}} \leq C \kappa.
  \end{align}
  % If $X$ does not have compact support, for every $q' \in [1,q)$ there is a constant $C' = C'(\gamma,\eta,q,q')$ such that \eqref{eqn:moment_sup_norm_bound} holds with $q$ and $C$ replaced by $q'$ and $C'$.  \vspace{0.1cm}
  
  \item[(ii)]  Assume that \eqref{eqn:diff_est_stoch_proc} and \eqref{eqn:moment_growths_stoch_proc} hold for every $q \geq 1$ with $\kappa = \kappa(q) \leq \sqrt{q} \hat{\kappa}$ and some $\eta \in (0,\infty)$.
  %. Moreover, assume that $X$ has compact support or that \eqref{eqn:moment_growths_stoch_proc} holds for every $q \geq 1$ with $\kappa = \kappa(q) \leq \sqrt{q} \hat{\kappa}$ and some $\eta \in (0,\infty)$. 
  
  Then $\sup_{x \in D} \|X(x)\|$ has Gaussian tails and there is a constant $C = C(\gamma,\eta)$ such that
  \begin{align*}
   \left\| \sup_{x \in D} \big\|X(x)\big\| \right\|_{L^q} \leq C \sqrt{q} \hat{\kappa}
  \end{align*}
  for every $q \geq 1$.  
\end{itemize}  }
\end{cor}

\ssk

\begin{Dem}
 % We first assume that $X$ does not have compact support. 
 We start with proving (i). Let $x,y \in D$ such that $|x - y|\geq 1$. Then, by \eqref{eqn:moment_growths_stoch_proc},
  \begin{align*}
  \big\|X(x) - X(y) \big\|_{L^q} \leq \big\|X(x)  \big\|_{L^q} + \big\| X(y) \big\|_{L^q} \leq 2 \kappa \leq 2 \kappa |x-y|^{\gamma}
 \end{align*}
 which shows that 
 \begin{align}\label{eqn:diff_est_stoch_proc_2} 
   \big \|X(x) - X(y) \big\|_{L^q} \leq 2\kappa |x - y|^{\gamma}
 \end{align}
 holds for every $x,y \in D$. Now let $x,y \in D$ such that $n-1 \leq |x|, |y| < n$. Interpolating between the inequality
 \begin{align*}
  \big\|X(x) - X(y) \big\|_{L^q} \leq \big\|X(x) \big\|_{L^q} + \big\| X(y)\big\|_{L^q} \leq \frac{2 \kappa}{1 + (n-1)^{\eta}}
 \end{align*}
 and inequality \eqref{eqn:diff_est_stoch_proc_2}, we see that for every $\lambda \in [0,1]$,
  \begin{align*}
 \big\|X(x) - X(y) \big\|_{L^q} \leq C \kappa n^{-(1-\lambda)\eta} |x - y|^{\gamma \lambda}
 \end{align*}
 and
 \begin{align*}
  \big\|X(x) \big\|_{L^q} \leq C \kappa n^{-(1 - \lambda)\eta}
 \end{align*}
 for every $x, y \in D_n$. Set $\gamma' := \gamma \lambda$ and $a_n := n^{-(1 - \lambda)\eta}$. In order to obtain $(a_n n^{\gamma'}) \in \ell^q(\NN)$, we must have $q(\lambda \gamma -(1 - \lambda)\eta)< -1$ which is equivalent to
  \begin{align*}
   \frac{1}{q} < \eta - \lambda(\eta + \gamma).
  \end{align*}
  The condition $\gamma' > \frac{d}{q}$ is equivalent to
  \begin{align*}
   \frac{1}{q} < \frac{\lambda \gamma }{d}.
  \end{align*}
  Choosing $\lambda^* = \frac{\eta}{\gamma/d + \eta + \gamma} \in (0,1)$, we have
  \begin{align*}
   \eta - \lambda^*(\eta + \gamma) = \frac{\lambda^* \gamma }{d} = \frac{\gamma \eta}{\gamma + d(\gamma + \eta)}
  \end{align*}
  which is indeed smaller than $\frac{1}{q}$ by assumption. Hence we may apply Lemma \ref{lem:cond_proc_bdd} to conlude (i). The claim (ii) follows by applying Lemma \ref{lem:cond_proc_bdd} (ii). 
%   In the case when $X$ has compact support, \eqref{eqn:diff_est_stoch_proc} implies that
%  \begin{align*}
%     \big\|X(x) - X(y) \big\|_{L^q} \leq C \kappa |x - y|^{\gamma}
%  \end{align*}
%   holds for every $x,y \in D$ where the constant $C$ depends on the diameter of the support of $X$. In (i) and (ii), the claim now follows from the usual Garsia-Rodemich-Rumsey Lemma.
\end{Dem}

\medskip

\noindent \textbf{\textsf{Example.}} \textit{  
 Consider the Gaussian process $X \colon (0,\infty) \to \RR$ where $X_t = \frac{B_t}{1 + t^{\alpha}}$, for a standard Brownian motion $B$, and $\alpha \in (\frac{1}{2},1]$. Then, if $t > 0$,
 \begin{align*}
  \| X_t \|_{L^q} \lesssim \sqrt{q} \| X_t \|_{L^2} \lesssim \frac{\sqrt{q}}{1 + t^{\alpha - \frac{1}{2}}}
 \end{align*}
 and for $s<t$,
 \begin{align*}
  \|X_t - X_s\|_{L^2} &= \frac{ \| B_t (1 + s^{\alpha}) - B_s (1 + t^{\alpha}) \|_{L^2}}{(1 + s^{\alpha}) (1 + t^{\alpha})}  \leq \frac{\| B_t - B_s \|_{L^2}}{1 + t^{\alpha}} + \frac{\|B_s\|_{L^2}|t^{\alpha} - s^{\alpha}|}{(1 + s^{\alpha})(1 + t^{\alpha})} \\
  &\leq \frac{|t-s|^{\frac{1}{2}}}{1 + t^{\alpha}} + \frac{2 s^{\frac{1}{2}} t^{\alpha - \frac{1}{2}}}{(1 + s^{\alpha})(1 + t^{\alpha})}|t-s|^{\frac{1}{2}} \\
  &\leq |t-s|^{\frac{1}{2}} 
 \end{align*}
 and
 \begin{align*}
  \big\|X_t - X_s\big\|_{L^q} \lesssim \sqrt{q} \big\|X_t - X_s\big\|_{L^2}.
 \end{align*}
  Applying part (ii) in Corollary \ref{cor:cond_proc_bdd} shows that the random variable 
 \begin{align*}
  \sup_{t \in (0,\infty)} \big|X_t\big|
 \end{align*}
 is finite and has Gaussian tails. Note that this is sharp in the sense that the law of the iterated logarithm for a Brownian motion implies that it is not possible to choose $\alpha = \frac{1}{2}$.  } 

\medskip 

Next, we apply the same ideas to give conditions for H\"older continuity. 
 
\ssk 
 
\begin{lem}\label{lem:cond_proc_unif_hoelder} {\sf
Let $(E,d)$ be a complete separable metric space, $D$ an open subset of $\RR^d$, $X \colon D \to (E,d)$ a continuous stochastic process and $\kappa > 0$. Set
 \begin{align*}
  D_n := \Big\{x \in D \, :\, n-1 \leq |x| < n \Big\}
 \end{align*}
and $\mathbf{N} := \{n \in \NN\, :\, D_n \neq \emptyset\}$. Let $(a_n)_{n \in \mathbf{N}}$ be a sequence of non-negative real numbers.

 \begin{itemize}
  \item[(i)] Assume that there is a $q > 1$ and a $\gamma \in (\frac{d}{q},1]$ such that for every $n \in \mathbf{N}$ and every $x,y \in D_n$ with $0 < |x - y| \leq 4$,
 \begin{align*}
   \left\| d(X(x),X(y)) \right\|_{L^q} \leq \kappa a_n |x - y|^{\gamma}.
 \end{align*} 
 Let $\beta \in \big(0, \gamma - \frac{d}{q}\big)$. Define the sequence $(b_n) := \big(a_n n^{\gamma - \beta}\big)$, and assume that $\|b \|_{l^q} \leq K < \infty$. Then there is a constant $C = C(q,\gamma)$ such that
 \begin{align*}
  \left\| \sup_{\stackrel{x, y \in D}{0 < |x-y| \leq 1}} \frac{d\big(X(x),X(y)\big)}{|x - y|^{\beta}} \right\|_{L^{q}} \leq C K \kappa.
 \end{align*} \vspace{0.1cm}
 
  \item[(ii)] Assume that there is some $\gamma \in (0,1]$ and that for every $q \geq 1$ there is a $c_q$ such that for every $n \in \mathbf{N}$ and every $x,y \in D_n$ with $0 < |x - y| \leq 4$,
 \begin{align*}
    \Big\| d\big(X(x),X(y)\big) \Big\|_{L^q} \leq \kappa c_q a_n |x - y|^{\gamma}
 \end{align*} 
 where $c_q = \mathcal{O}(\sqrt{q})$ when $q \to \infty$. Let $\beta \in (0, \gamma)$ and assume that $a_n = \mathcal{O}\Big(n^{-(\gamma - \beta)} (1 + \log(n))^{-\frac{1}{2}}\Big)$. Then for every $q \geq 1$ there is some constant $C = C(q,\gamma)$ such that
 \begin{align*}
  \left\| \sup_{\stackrel{x ,y \in D}{0 < |x-y| \leq 1}} \frac{d\big(X(x),X(y)\big)}{|x - y|^{\beta}} \right\|_{L^{q}} \leq C \kappa
 \end{align*}
  with $C = \mathcal{O}(\sqrt{q})$ when $q \to \infty$. In particular, the random variable 
  \begin{align*}
   \sup_{\stackrel{x, y \in D}{0 < |x-y| \leq 1}} \frac{d\big(X(x),X(y)\big)}{|x - y|^{\beta}}
  \end{align*}
  has Gaussian tails.
 \end{itemize}  }
\end{lem}

\ssk

\begin{Dem}
 Without loss of generality, one can choose $\kappa = 1$ and $\mathbf{N} = \NN$. For $n \in \NN$, set $\tilde{D}_n := \big\{D_n \cup D_{n+1} \cup D_{n +2}\big\}$. We first prove (i). Fix some $n \in \NN$ and some $k \in \NN$. Let $\alpha > \frac{d}{q}$ and define
\begin{align*}
   p_k(s) = \begin{cases}
	    s^{\alpha + \frac{d}{q}} &\text{if } s \in [0,4] \\
	    (4^{\alpha q + d} + k(s-4))^{\frac{1}{q}} &\text{if } s \geq 4.
          \end{cases}
\end{align*}
Fix $x, y \in \tilde{D}_n$ with $0 < |x- y| \leq 1$. From the Garsia-Rodemich-Rumsey Lemma,
\begin{align*}
  d\big(X(x),X(y)\big) \leq C  V_{n,k}^{\frac{1}{q}} |x-y|^{\alpha - \frac{d}{q}}
\end{align*}
where
\begin{align*} 
  V_{n,k} = \int_{\tilde{D}_n \times \tilde{D}_n} \frac{\Big| d\big(X(u),X(v)\big) \Big|^q}{p_k\big(|u - v|\big)^q}\, du\, dv.
\end{align*}
Thus
\begin{align*}
  \EE \left| \sup_{\stackrel{x,y \in \tilde{D}_n}{0 < |x - y| \leq 1}}  \frac{d\big(X(x),X(y)\big)}{|x - y|^{\alpha - \frac{d}{q}}} \right|^q &\leq C^q  \int_{\tilde{D}_n \times \tilde{D}_n}  \frac{\EE \Big[ d\big(X(u),X(v)\big)^q\Big]}{p_k(|u - v|)^q} \, du\, dv \\
  &\leq C^q \sum_{l = 0,1,2} \int_{D_{n + l} \times D_{n + l}}  \frac{\EE \Big[ d\big(X(u),X(v)\big)^q\Big]}{p_k(|u - v|)^q} \, du\, dv.
\end{align*} 
For every $m \in \NN$, we have
 \begin{align*}
  \int_{D_m \times D_m}  \frac{\EE \Big[ d\big(X(u),X(v)\big)^q\Big]}{p_k(|u - v|)^q} \, du\, dv &\leq a_m^q \int_{(D_m \times D_m) \cap \{|u - v| \leq 4\}} |u-v|^{(\gamma - \alpha)q - d} \, du\, dv \\
  &\quad +\int_{(D_m \times D_m) \cap \{|u - v| > 4\}} \frac{\EE \Big[ d\big(X(u),X(v)\big)^q\Big]}{4^{\alpha q + d} + k(|u - v| - 4)} \, du\, dv
 \end{align*}
Moreover, by a change of variables,
 \begin{align*}
  \int_{(D_m \times D_m) \cap \{|u - v| \leq 4 \}} |u-v|^{(\gamma - \alpha)q - d} \, du\, dv &\leq \int_{D_m \times D_m } |u-v|^{(\gamma - \alpha)q - d} \, du\, dv \\
  &= m^{d + (\gamma - \alpha)q} \int_{(0,1)^2}  |u-v|^{(\gamma - \alpha)q - d} \, du\, dv.
 \end{align*}
Set $\alpha = \frac{d}{q} + \beta < \gamma$. Then this integral is finite, and sending $k \to \infty$ shows that
 \begin{align*}
    \left\| \sup_{\stackrel{x,y \in \tilde{D}_n}{0 < |x - y| \leq 1}}  \frac{d\big(X(x),X(y)\big)}{|x - y|^{\beta}} \right\|_{L^q} &\leq C\Big(a_n n^{(\gamma - \beta)} + a_{n+1} (n+1)^{(\gamma - \beta)} + a_{n+2} (n+2)^{(\gamma - \beta)}\Big) \\
    &= C\big( b_n + b_{n+1} + b_{n+2} \big).
 \end{align*}
Now take $x,y \in D$ with $0 < |x - y| \leq 1$ and assume that
 \begin{align*}
    \frac{d\big(X(x),X(y)\big)}{|x - y|^{\beta}} \geq t.
 \end{align*}
Then there is an $n \in \NN$ such that $x \in D_n$ and since $|x - y|\leq 1$, $y \in \big\{D_{n-1} \cup D_n \cup D_{n+1}\big\}$, where we set $D_0 := D_1$. Thus we have shown that for every $t \geq 0$,
 \begin{align*}
  \left\{ \sup_{\stackrel{x,y \in D}{0 < |x-y| \leq 1}} \frac{d\big(X(x),X(y)\big)}{|x - y|^{\beta}} \geq t \right\} \subseteq \bigcup_{n \in \NN} \left\{ \sup_{\stackrel{x,y \in \tilde{D}_n}{0 < |x-y| \leq 1|}} \frac{d\big(X(x),X(y)\big)}{|x - y|^{\beta}} \geq t \right\}
 \end{align*}
and therefore
\begin{align*}
  \EE \left| \sup_{\stackrel{x,y \in D}{0 < |x-y| \leq 1}} \frac{d\big(X(x),X(y)\big)}{|x - y|^{\beta}} \right|^{q} \leq \sum_{n = 1}^{\infty} \EE \left| \sup_{\stackrel{x,y \in \tilde{D}_n}{0 < |x-y| \leq 1}} \frac{d\big(X(x),X(y)\big)}{|x - y|^{\beta}} \right|^{q} \leq 3 C^q K^q.
\end{align*}

%  \begin{align*}
%   \PP \left( \sup_{\stackrel{x,y \in D}{0 < |x-y| \leq 1}} \frac{d\big(X(x),X(y)\big)}{|x - y|^{\beta}} \geq t \right) &\leq \sum_{n = 1}^{\infty} \PP \left( \sup_{\stackrel{x,y \in \tilde{D}_n}{0 < |x-y| \leq 1}} \frac{d\big(X(x),X(y)\big)}{|x - y|^{\beta}} \geq t \right) \\
%   &\leq \frac{3^q C^q}{t^q} \sum_{n = 1}^{\infty} b_n^q = \frac{ 3^q C^q K^q}{t^q}.
%  \end{align*}
% Hence we can estimate
%  \begin{align*}
%   \EE \left| \sup_{\stackrel{x,y \in D}{0 < |x-y| \leq 1}} \frac{d\big(X(x),X(y)\big)}{|x - y|^{\beta}} \right|^{q'} \leq 1 + 3^q C^q K^q \int_1^{\infty} \frac{dt}{t^\frac{q}{q'}} \leq 1 + \frac{3^q C^q K^q}{1 - \frac{q}{q'}}.
%  \end{align*}
 
\ssk 
 
Now we prove (ii). Note that the constant in the Garsia-Rodemich-Rumsey Lemma may be chosen non-increasing in $q$. Therefore, we can argue similarly as before to see that for every $q \geq 1$ and $n \in \NN$,
 \begin{align*}
  \left\| \sup_{\stackrel{x,y \in \tilde{D}_n}{0 < |x - y| \leq 1}}  \frac{d\big(X(x),X(y)\big)}{|x - y|^{\beta}} \right\|_{L^q} \leq C_q \big(b_n + b_{n+1} + b_{n+2}\big)
 \end{align*}
where $C_q = \mathcal{O}(\sqrt{q})$. This shows that the random variable has Gaussian tails, i.e. there is some constant $C$ such that for every $n \in \NN$
 \begin{align*}
   \PP \left(\sup_{\stackrel{x,y \in \tilde{D}_n}{0 < |x - y| \leq 1}}  \frac{d\big(X(x),X(y)\big)}{|x - y|^{\beta}} \geq t \right) \leq C \exp \left( - \frac{t^2}{C (b_n + b_{n+1} + b_{n+2})^2} \right) 
 \end{align*}
for every $t \geq 0$. Hence
 \begin{align*}
  \PP \left( \sup_{\stackrel{x,y \in D}{0 < |x - y| \leq 1}}  \frac{d\big(X(x),X(y)\big)}{|x - y|^{\beta}} \geq t \right) &\leq C \sum_{n = 1}^{\infty} \exp \left( - \frac{t^2}{C (b_n + b_{n+1} + b_{n+2})^2} \right) \\ 
  &\leq C \sum_{n = 1}^{\infty} \exp \left( - \frac{t^2}{C} (1 + \log(n)) \right) \\
  &\leq C \exp \left( - \frac{t^2}{C} \right) \sum_{n = 1}^{\infty} n^{-t^2/C}
 \end{align*}
and the sum is finite for $t$ large enough. This proves that
 \begin{align*}
  \sup_{\stackrel{x,y \in D}{0 < |x - y| \leq 1}}  \frac{d\big(X(x),X(y)\big)}{|x - y|^{\beta}}
 \end{align*}
has Gaussian tails.
% which is equivalent to say that its $L^q$ norm grows like $\sqrt{q}$ when $q \to \infty$.
\end{Dem}

\medskip

\begin{example}
Let $X \colon (0,\infty) \to \RR$ be the Gaussian process defined as
 \begin{align*}
  X_t = \frac{B_t}{\sqrt{t \log(1 + t)}},
 \end{align*}
$B$ being a standard Brownian motion. Then
 \begin{align*}
  \|X_t - X_s \|_{L^2} \leq \frac{|t-s|^{\frac{1}{2}}}{\sqrt{t \log(1 + t)}} + \frac{\sqrt{t \log(1 + t)} - \sqrt{s \log(1 + s)}}{\sqrt{\log(1 + s)}\sqrt{s \log(1 + s)}}.
 \end{align*}
By the mean value theorem,
 \begin{align*}
  t \log(1 + t) - s \log(1 + s) \leq (\log(1 + t) + 1)(t-s)
 \end{align*}
and therefore
 \begin{align*}
  \sqrt{t \log(1 + t)} - \sqrt{s \log(1 + s)} \leq \sqrt{t \log(1 + t) - s \log(1 + s)} \leq \sqrt{(\log(1 + t) + 1)}|t-s|^{\frac{1}{2}}.
 \end{align*}
If $(n - 1) \leq s \leq t \leq n$, we have for any $q \geq 2$
 \begin{align*}
  \|X_t - X_s \|_{L^q} \lesssim \sqrt{q} \|X_t - X_s \|_{L^2} \lesssim \sqrt{q} a_n |t-s|^{\frac{1}{2}}
 \end{align*}
with $a_n = \mathcal{O}\Big(n^{-\frac{1}{2}} (1 + \log(n))^{-\frac{1}{2}}\Big)$. Part (ii) of Lemma \ref{lem:cond_proc_unif_hoelder} shows that for any $\beta \in \big(0,\frac{1}{2}\big)$, the random variable
 \begin{align*}
  \sup_{\stackrel{s,t \in (0,\infty)}{0 < |t - s| \leq 1}}  \frac{\big|X_t - X_s\big|}{|t - s|^{\beta}}
 \end{align*}
is finite and has Gaussian tails.
\end{example}

\medskip

\begin{cor}\label{cor:cond_proc_unif_hoelder} {\sf
 Let $D$ be an open subset of $\RR^d$, and $\big(E, \| \cdot \|\big)$ be a separable Banach space. Let $X \colon D \to E$ be a continuous stochastic process and $q > 1$.
 \begin{itemize}
  \item[(i)] Assume that there are constants $\kappa > 0$, $\gamma \in (0,1]$ and $\beta \in (0, \gamma)$ %, $\eta > 0$ 
  such that for every $x,y \in D$ with $0 < |x - y| \leq 1$,
 \begin{align}\label{eqn:diff_est_stoch_proc_unif_hol}
    \big\|X(x) - X(y) \big\|_{L^q} \leq \kappa |x - y|^{\gamma}
 \end{align}
 and that there is an $\eta \in (0,\infty)$ such that for every $x \in D$,
 \begin{align}\label{eqn:moment_growths_stoch_proc_unif_hol}
  \big\|X(x) \big\|_{L^q} \leq \frac{\kappa}{1 + |x|^{\eta}}
 \end{align}
  where $q > 1$ satisfies
  \begin{align*}
   q > \frac{\gamma}{\eta(\gamma - \beta)} + d\left(\frac{\gamma}{\eta(\gamma - \beta)} + \frac{1}{\gamma - \beta} \right).
  \end{align*}
  Then the random variable 
  \begin{align*}
   \sup_{\stackrel{x,y \in D}{0 < |x - y| \leq 1}} \frac{\big\|X(x) - X(y)\big\|}{|x - y|^{\beta}}
  \end{align*}
  is almost surely finite. Moreover, there is a constant $C = C(\gamma,\eta,\beta,q)$ such that
  \begin{align}\label{eqn:moment_bdd_hoelder_norm}
   \left\| \sup_{\stackrel{x,y \in D}{0 < |x - y| \leq 1}} \frac{\big\|X(x) - X(y)\big\|}{|x - y|^{\beta}} \right\|_{L^{q}} \leq C \kappa.
  \end{align}
  % If $X$ does not have compact support, for every $q' \in [1,q)$ there is a constant $C' = C'(\gamma,\eta,\beta,q,q')$ such that \eqref{eqn:moment_bdd_hoelder_norm} holds with $q$ and $C$ replaced by $q'$ resp. $C'$.  \vspace{0.1cm}
  
  \item[(ii)]  Assume that \eqref{eqn:diff_est_stoch_proc_unif_hol} and \eqref{eqn:moment_growths_stoch_proc_unif_hol} hold for every $q \geq 1$ with $\kappa = \kappa(q) \leq \sqrt{q} \hat{\kappa}$ and some $\eta \in (0,\infty)$. % Moreover, assume that $X$ has compact support or that there is some $\eta \in (0,\infty)$ such that \eqref{eqn:moment_growths_stoch_proc_unif_hol} holds for every $q \geq 1$ with $\kappa = \kappa(q) \leq \sqrt{q} \hat{\kappa}$.
  
  Then for every $\beta \in (0, \gamma)$,
  \begin{align*}
   \sup_{\stackrel{x,y \in D}{0 < |x - y| \leq 1}} \frac{\big\|X(x) - X(y)\big\|}{|x - y|^{\beta}}
  \end{align*}
  has Gaussian tails, and there is a constant $C = C(\gamma,\eta,\beta)$ such that
  \begin{align*}
   \left\| \sup_{\stackrel{x,y \in D}{0 < |x - y| \leq 1}} \frac{\big\|X(x) - X(y)\big\|}{|x - y|^{\beta}} \right\|_{L^q} \leq C \sqrt{q} \hat{\kappa}
  \end{align*}
  holds for every $q \geq 1$. 
\end{itemize}   }
\end{cor}

\ssk

\begin{Dem}
  The proof is similar to the proof of Corollary \ref{cor:cond_proc_bdd}, using Lemma \ref{lem:cond_proc_unif_hoelder} above. We leave the details to the reader.
% We will only prove the case of $X$ not having compact support, the other case is similar (and even easier). We start with (i). By \eqref{eqn:diff_est_stoch_proc_unif_hol} and \eqref{eqn:moment_growths_stoch_proc_unif_hol},
%  \begin{align}
%     \big\|X(x) - X(y) \big\|_{L^q} \leq 2\kappa |x - y|^{\gamma}
%  \end{align}
% holds for every $x,y \in D$. By interpolation, for every $\lambda \in [0,1]$,
%   \begin{align*}
%   \big\|X(x) - X(y) \big\|_{L^q} \leq C \kappa n^{-(1-\lambda)\eta} |x - y|^{\gamma \lambda}
%  \end{align*}
% and
%  \begin{align*} 
%   \big\|X(x) \big\|_{L^q} \leq C \kappa n^{-(1 - \lambda)\eta}
%  \end{align*} 
% for every $x, y \in D_n$. Set $\gamma' := \gamma \lambda$ and $a_n := n^{-(1 - \lambda)\eta}$. The condition $\big(a_n n^{\gamma' - \beta}\big) \in \ell^q(\NN)$ is satisfied when
%   \begin{align*}
%    \frac{1}{q} < \eta + \beta - \lambda(\eta + \gamma).
%   \end{align*}
% The condition $\beta < \gamma' - \frac{d}{q}$ is equivalent to
%   \begin{align*}
%    \frac{1}{q} < \frac{\lambda \gamma - \beta }{d}.
%   \end{align*}
% Choosing $\lambda^* = (\eta + \beta + \beta/d)/(\gamma/d + \eta + \gamma) \in (0,1)$, we have
%   \begin{align*}
%    \eta + \beta - \lambda^*(\eta + \gamma) = \frac{\lambda^* \gamma - \beta }{d} = \frac{(\gamma - \beta) \eta}{\gamma + d(\gamma + \eta)}
%   \end{align*}
% which is smaller than $\frac{1}{q}$ by assumption, and the claim follows from Lemma \ref{lem:cond_proc_unif_hoelder}. (ii) follows similarly.
\end{Dem}

\medskip

If $D$ is an open subset of $\RR^d$, and $\big(E,\| \cdot \|\big)$ is a normed space, $f \colon D \to E$ a function and $\rho \in (0,1]$, we define
\begin{align*}
 \| f \|_{\mathcal{C}^{\rho}}^{\sim} := \max \left\{ 2 \sup_{x \in D} \big\|f(x)\big\| , \sup_{\stackrel{x, y \in D}{0 < |x - y| \leq 1}} \frac{\big\|f(x) - f(y) \big\|}{|x - y|^{\rho}} \right\}.
\end{align*}
Let $f, g \colon D \to \RR^{m}$. Then we define the function $(f \odot g ) \colon D \to \RR^{m \times m}$ by setting $(f \odot g)^{ij}(x) = f^i(x) g^j(x)$. Note that  $\| \cdot \|_{\mathcal{C}^{\rho}}^{\sim}$ is equivalent to $\| \cdot \|_{\mathcal{C}^{\rho}}$ and that we have
\begin{align}\label{eqn:comp_norms}
 \| f \odot g \|_{\mathcal{C}^{\rho}}^{\sim} \leq \| f \|_{\mathcal{C}^{\rho}}^{\sim} \| g \|_{\mathcal{C}^{\rho}}^{\sim}
\end{align}
provided we equip $\RR^{m}$ and $\RR^{m \times m}$ with the sup norm. 

In the following, we will consider stochastic processes  $V \colon D \times [0,T] \to \RR^{m}$ and $W \colon D \times \{0 \leq s \leq t \leq T\} \to \RR^{m \times m}$ for which we assume that for every $s < u < t \in [0,T]$ and every $x \in D$,
   \begin{align}\label{eqn:algebraic_rel_first_second_level}
    W_{ts}(x) - W_{us}(x) - W_{tu}(x) = V_{us}(x) \odot \tilde{V}_{tu}(x)
   \end{align}
   holds almost surely. The next theorem is the main result of this section.

\medskip 

\begin{thm}[Kolmogorov criterion for rough drivers]\label{thm:kolmogorov_rough_drivers} {\sf
Let $D$ be an open subset of $\RR^d$, $\kappa > 0$ and $\gamma_1, \gamma_2 \in (0,1]$.

  \begin{itemize}
   \item[(i)] Let $V \colon D \times I \to \RR^{m}$ be a stochastic process and $q > 1$. Assume that for every $x,y \in D$ with $0 < |x-y| \leq 1$ and $s < t \in I$,
   \begin{align}\label{eqn:incr_bdd_first_level}
    \big\| V_{ts}(x) - V_{ts}(y) \big\|_{L^q} \leq \kappa |t - s|^{\gamma_1} |x-y|^{\gamma_2}
   \end{align}
   %Assume that there is a compact set in $\RR^d$ such that $V_{ts}$ is supported on this set almost surely for every $s < t \in I$, or 
   and that there is an $\eta \in (0, \infty)$ such and that for every $x \in D$ and $s<t \in I$,
   \begin{align}\label{eqn:growth_bdd_first_level}
    \big\| V_{ts}(x) \big\|_{L^q} \leq \frac{\kappa |t - s|^{\gamma_1}}{1 + |x|^{\eta}}. 
   \end{align}
   Let $\alpha \in (0, \gamma_1)$, $\beta \in (0,\gamma_2)$ and assume that
   \begin{align}\label{eqn:q_suff_large_kolmogorov}
    q > \max \left\{ \frac{\gamma_2}{\eta(\gamma_2 - \beta)} + d\left(\frac{\gamma_2}{\eta(\gamma_2 - \beta)} + \frac{1}{\gamma_2 - \beta} \right), \frac{1}{\gamma_1 - \alpha} \right\}.
   \end{align}
   % where we set $\eta = \infty$ in the case of $V$ having compact support (in the sense above). 
   Then there is a continuous modification of the process $V$. Moreover, 
   %if $V$ has compact support, 
   there is a constant $C = C(\gamma_1,\gamma_2,\alpha,\beta,\eta,d,T,q)$ such that
   \begin{align}\label{eqn:lq_estimates_first_level}
    \left\| \sup_{ s < t \in I} \frac{\big\| V_{ts} \big\|_{\mathcal{C}^{\beta}}}{|t-s|^{\alpha}} \right\|_{L^{q}} \leq C \kappa.
   \end{align}
%    If $V$ does not have compact support, for every $q' \in [1,q)$ there is a constant $C' = C'(\gamma_1,\gamma_2,\alpha,\beta,\eta,d,T,q,q')$ such that \eqref{eqn:lq_estimates_first_level} holds with $q$ and $C$ replaced by $q'$ resp. $C'$. \vspace{0.1cm}
      
   \item[(ii)] In addition, let $W \colon D \times \{0 \leq s \leq t \leq T\} \to \RR^{m \times m}$ be a stochastic process for which the relation \eqref{eqn:comp_norms} holds and $q > 2$. Assume that for every $x,y \in D$ with $0< |x - y| \leq 1$ and $s < t \in I$,
   \begin{align}\label{eqn:incr_bdd_second_level}
    \| W_{ts}(x) - W_{ts}(y) \|_{L^{\frac{q}{2}}} \leq \kappa^2 |t - s|^{2\gamma_1} |x-y|^{\gamma_2}
   \end{align}
   % Assume that $W$ has compact support (in the sense above) or 
   and that there is an $\eta \in (0,\infty)$ such that for every $x \in D$ and $s<t \in I$,
   \begin{align}\label{eqn:growth_bdd_second_level}
    \big\| W_{ts}(x) \big\|_{L^{\frac{q}{2}}} \leq \frac{\kappa^2 \, |t - s|^{2\gamma_1}}{1 + |x|^{2 \eta}}.
   \end{align}
    Let $\alpha \in (0, \gamma_1)$, $\beta \in (0,\gamma_2)$ and assume that $q$ is sufficiently large such that \eqref{eqn:q_suff_large_kolmogorov} holds. 
    %  (again with $\eta = \infty$ in the case of $W$ having compact support). Moreover, assume that for every $s < u < t \in I$ and every $x \in D$,
%    \begin{align}\label{eqn:algebraic_rel_first_second_level}
%     W_{ts}(x) - W_{us}(x) - W_{tu}(x) = V_{us}(x) \odot \tilde{V}_{tu}(x)
%    \end{align}
%    almost surely, where $V \colon D \times I \to \RR^{d_1}$ and $\tilde{V} \colon D \times I \to \RR^{d_2}$ are stochastic processes which satisfy the conditions of part (i), with either all of these processes having compact support or satisfying the respective growth conditions for the same $\eta > 0$. 
   Then there is a continuous modification of the process $W$. Moreover, 
   % if $W$ has compact support, 
   there is a constant $C = C(\gamma_1,\gamma_2,\alpha,\beta,\eta,d,T,q)$ such that
  \begin{align}\label{eqn:lq_estimates_second_level}
    \left\| \sup_{ s < t \in I} \frac{\big\| W_{ts} \big\|_{\mathcal{C}^{\beta}}}{|t-s|^{2 \alpha}} \right\|_{L^{\frac{q}{2}}} \leq C^2 \kappa^2.
   \end{align}
%    If $W$ does not have compact support, for every $q' \in [2,q)$ there is a constant $C' = C'(\gamma_1,\gamma_2,\alpha,\beta,\eta,d,T,q,q')$ such that \eqref{eqn:lq_estimates_second_level} holds with $q$ and $C$ replaced by $q'$ resp. $C'$. \vspace{0.1cm}
   
   \item[(iii)] Assume in addition that \eqref{eqn:incr_bdd_first_level}, \eqref{eqn:incr_bdd_second_level} and the growth conditions \eqref{eqn:growth_bdd_first_level} and \eqref{eqn:growth_bdd_second_level} hold for every $q \geq 2$ with constant $\kappa \leq \sqrt{q} \hat{\kappa}$. Then for every $\alpha \in (0,\gamma_1)$ and $\beta \in (0,\gamma_2)$, the random variables
   \begin{align*}
    \sup_{ s < t \in I} \frac{\big\| V_{ts} \big\|_{\mathcal{C}^{\beta}}}{|t-s|^{\alpha}} \qquad \text{and} \qquad \sqrt{\sup_{s < t \in I} \frac{\big\| W_{ts} \big\|_{\mathcal{C}^{\beta}}}{|t-s|^{2 \alpha}}}
   \end{align*}
   have Gaussian tails.
   % , and the $L^q$-estimates \eqref{eqn:lq_estimates_first_level} and \eqref{eqn:lq_estimates_second_level} hold for every $q \geq 1$ with $\kappa$ replaced by $\hat{\kappa}$, and the constant $C$ depends on $q$ in such a way that $C(q) = \mathcal{O}(\sqrt{q})$.
 \end{itemize}  } 
\end{thm}

\ssk

\begin{Dem}
Without loss of generality we may assume $\kappa = 1$, otherwise we can replace $V$ and $W$ by $V/\kappa$ resp. $W/ \kappa^2$. Furthermore, we will prove the result for the $\|\cdot \|_{\mathcal{C}^{\beta}}^{\sim}$ norm, claimed results follow by equivalence of norms.
 
% As in the previous proofs, we will only consider the case of $V$ and $W$ not having compact support, the other case is analogous. 
 We start with proving (i). Fix $s < t$. Using \eqref{eqn:incr_bdd_first_level} and the classical Kolmogorov theorem \cite[Theorem 3.23]{Kal02}, there is a continuous modification of the process $x \mapsto V_{ts}(x)$
 %on every bounded hypercube contained in $D$ with vertices in $D \cap \QQ^d$, and since $D$ is the countable union of such hypercubes also 
 on $D$. The estimates \eqref{eqn:incr_bdd_first_level} and \eqref{eqn:growth_bdd_first_level} and Corollary \ref{cor:cond_proc_unif_hoelder} imply that
 \begin{align*}
  \left\| \sup_{\stackrel{x,y \in D}{0 < |x - y| \leq 1}} \frac{\big| V_{ts}(x) - V_{ts}(y)\big|}{|x - y|^{\beta}} \right\|_{L^{q}} \leq C|t-s|^{\gamma_1}
 \end{align*}
 and Corollary \ref{cor:cond_proc_bdd} gives
 \begin{align*}
  \left\| \sup_{x \in D} \big| V_{ts}(x) \big| \right\|_{L^{q}} \leq C|t-s|^{\gamma_1}.
 \end{align*}
 Note in particular that the constant on the right hand side of both equations is independent of $s$ and $t$. We can repeat this procedure for every $s < t$ and obtain a process $t \mapsto V_t$ which, for every $t \in [0,T]$, takes values in $\mathcal{C}^{\beta}_b$ almost surely, and for which
 \begin{align}\label{eqn:incr_bdd_first_level_Cb_valued_proc}
    \Big\| \big\| V_t - V_s \big\|_{\mathcal{C}^{\beta}}^{\sim} \Big\|_{L^{q}} \leq C|t-s|^{\gamma_1} 
 \end{align}
 holds for every $s<t$. Applying again the Kolmogorov theorem for Banach space valued processes gives the claim.
 
\ssk 
 
We proceed with (ii). 
% Again, we may assume that $q' \in \big(\frac{1}{\gamma_1 - \alpha)}q\big)$. 
As in (i), for every $s<t$ there are modifications of the process $x \mapsto W_{ts}(x)$ such that
 \begin{align*}
    \Big\| \big\| W_{ts} \big\|_{\mathcal{C}^{\beta}}^{\sim} \Big\|_{L^{\frac{q}{2}}} \leq C|t-s|^{2\gamma_1}.
 \end{align*}
Using the algebraic relation \eqref{eqn:algebraic_rel_first_second_level}, the estimate \eqref{eqn:incr_bdd_first_level_Cb_valued_proc} for $V$ and the compatibility of the $\| \cdot \|_{\mathcal{C}^{\beta}}^{\sim}$ norms given in \eqref{eqn:comp_norms}, we can mimic the proof of the Kolmogorov criterion for rough paths (\cite[Theorem 3.1]{FH14}) to conlude.

\ssk
 
Assertion (iii) follows similarly by using part (ii) in the Corollaries \ref{cor:cond_proc_unif_hoelder} and \ref{cor:cond_proc_bdd}. 
\end{Dem}

\medskip

Finally, we give a Kolmogorov criterion for the distance between rough drivers, whose proof is very similar to the proof of Theorem \ref{thm:kolmogorov_rough_drivers}, using the Kolmogorov criterion for rough path distance \cite[Theorem 3.3]{FH14}; we leave it to the reader. 

\ssk

\begin{thm}[Kolmogorov criterion for rough driver distance]\label{thm:kolmogorov_rough_driver_distance} {\sf
  Let  $\kappa > 0$, $\gamma_1, \gamma_2 \in (0,1]$ and $(V,W)$, $(\hat{V},\hat{W})$ processes as in Theorem \ref{thm:kolmogorov_rough_drivers}. Set $\Delta V := V - \hat{V}$ and $\Delta W := W - \hat{W}$.
  % Let $I \subset \RR$ be a closed interval with $|I| = T$, $D$ an open subset of $\RR^d$ and $\kappa > 0$. Let $\gamma_1, \gamma_2 \in (0,1]$, %$\eta \in (0,\infty)$ 
%  and $q>2$. Let $(V,W)$ and $(\hat{V},\hat{W})$ be processes as in Theorem \ref{thm:kolmogorov_rough_drivers}. Set $\Delta V := V - \hat{V}$ and $\Delta W := W - \hat{W}$. 
  \begin{itemize}
   \item[(i)] Assume that $(V,W)$ and $(\hat{V},\hat{W})$ satisfy the same moment conditions as in Theorem \ref{thm:kolmogorov_rough_drivers} with $q$ sufficiently large as in \eqref{eqn:q_suff_large_kolmogorov}. Moreover, assume that there is an $\varepsilon > 0$ such that
   \begin{align}\label{eqn:incr_bdd_distance_first_level}
    \big\| \Delta V_{ts}(x) - \Delta V_{ts}(y) \big\|_{L^q} \leq \varepsilon \kappa |t - s|^{\gamma_1} |x-y|^{\gamma_2}
   \end{align}
   and
   \begin{align}\label{eqn:incr_bdd_distance_second_level}
    \big\|\Delta W_{ts}(x) - \Delta W_{ts}(y) \big\|_{L^{\frac{q}{2}}} \leq \varepsilon \kappa^2 |t - s|^{2\gamma_1} |x-y|^{\gamma_2}
   \end{align}
   for every $x,y \in D$ such that $0 < |x - y| \leq 1$ and every $s<t \in [0,T]$ and that
   \begin{align}\label{eqn:growth_bdd_distance_first_level}
    \big\| \Delta V_{ts}(x) \big\|_{L^q} \leq \varepsilon \frac{\kappa |t - s|^{\gamma_1}}{1 + |x|^{\eta}} 
   \end{align}
   and
   \begin{align}\label{eqn:growth_bdd_distance_second_level}
    \big\| \Delta W_{ts}(x) \big\|_{L^{\frac{q}{2}}} \leq \varepsilon \frac{\kappa^2 |t - s|^{2\gamma_1}}{1 + |x|^{2\eta}}
   \end{align}
   hold for every $x \in D$ and every $s < t \in [0,T]$. Then there are continuous modifications of the processes $(V,W)$ and $(\hat{V},\hat{W})$. Moreover, there is a constant $C = C(\gamma_1,\gamma_2,\alpha,\beta,\eta,d,T,q)$ such that
   \begin{align}\label{eqn:dist_est_first_level}
    \left\| \sup_{ s < t \in I} \frac{\big\| V_{ts} - \hat{V}_{ts}\big\|_{\mathcal{C}^{\beta}}}{|t-s|^{\alpha}} \right\|_{L^{q}} \leq \varepsilon C \kappa
   \end{align}
   and
   \begin{align}\label{eqn:dist_est_second_level}
    \left\| \sup_{ s < t \in I} \frac{\big\| W_{ts} - \hat{W}_{ts}\big\|_{\mathcal{C}^{\beta}}}{|t-s|^{2 \alpha}} \right\|_{L^{\frac{q}{2}}} \leq \varepsilon C^2 \kappa^2.
   \end{align}
%    If the processes do not have compact support, for every $q' \in [1,q)$ there is some constant $C' = C'(\gamma_1,\gamma_2,\alpha,\beta,\eta,d,T,q,q')$ such that \eqref{eqn:dist_est_first_level} and \eqref{eqn:dist_est_second_level} hold with $q$ and $C$ replaced by $q'$ resp. $C'$.  \vspace{0.1cm}
   
   \item[(ii)] Assume in addition that $(V,W)$ and $(\hat{V},\hat{W})$ satisfy the same moment conditions as in Theorem \ref{thm:kolmogorov_rough_drivers} for every $q \geq 2$ and that \eqref{eqn:incr_bdd_distance_first_level}, \eqref{eqn:incr_bdd_distance_second_level}, \eqref{eqn:growth_bdd_distance_first_level} and \eqref{eqn:growth_bdd_distance_second_level} also hold for every $q \geq 2$ with common constant  $\kappa \leq \sqrt{q} \hat{\kappa}$. Then there is a constant $C = C(\gamma_1,\gamma_2,\alpha,\beta,\eta,d,T)$ such that
   \begin{align*}
    \varepsilon^{-1} \left\| \sup_{ s < t \in I} \frac{\big\| V_{ts} - \hat{V}_{ts}\big\|_{\mathcal{C}^{\beta}}}{|t-s|^{\alpha}} \right\|_{L^{q}} + \varepsilon^{-\frac{1}{2}} \left\| \sqrt{ \sup_{ s < t \in I} \frac{\big\| W_{ts} - \hat{W}_{ts}\big\|_{\mathcal{C}^{\beta}}}{|t-s|^{2 \alpha}}} \right\|_{L^{q}} \leq \sqrt{q} \hat{\kappa} C
   \end{align*}
   for every $q \geq 2$.
\end{itemize}   }
\end{thm}

\bigskip

%-------------------------------------------%
\section{Stochastic and rough flows}
\label{SectionStochasticRoughFlows}
%-------------------------------------------%

The theory of stochastic flows grew out of the pioneering works of the Russian school \cite{BL61, GS72} on the dependence of solutions to stochastic differential equations with respect to parameters and the proof by Bismut \cite{B81} and Kunita \cite{Kun81} that stochastic differential equations generate continuous flows of diffeomorphisms under proper regularity conditions on the driving vector fields. The Brownian character of these random flows, that is the fact that they are continuous with stationary and independent increments, was inherited from the Brownian character of their driving noise. The next natural step consisted in the study of Brownian flows for themselves. After the works of Harris \cite{Ha81}, Baxendale \cite{Bax80} and Le Jan \cite{LeJan82}, they appeared to be generated by stochastic differential equations driven by infinitely many Brownian motions, or better, to be in one-to-one correspondence with vector field-valued Brownian motions. A probabilistic integration theory of such random 
time-varying velocity fields was developed to establish that correspondence, and it was extended by Le Jan and Watanabe \cite{LeJanW82} to a large class of continuous semimartingale flows and continuous semimartingale velocity fields. Kunita \cite{Kun86, KunTata} studied the problem of convergence of stochastic flows, with applications to averaging and homogenization results, and promoted the use of stochastic flows to implement a version of the characteristic method in the setting of first and second order stochastic partial differential equations, notably those coming from the nonlinear filtering theory.

\medskip

We shall show in this Section that the theory of semimartingale stochastic flows can be embedded into the theory of rough flows developed in Section \ref{SectionRoughFlows}. We review in Section \ref{SubsectionSemimartingaleFlows} the basics of the theory of stochastic flows and show in Section \ref{SubsectionMartingaleDrivers} that sufficiently regular (semi)martingale velocity fields can be lifted to rough drivers; this is done using our Kolmogorov-type criterion for rough drivers, Theorem \ref{thm:kolmogorov_rough_drivers}. The identification of (semi)martingale flows generated by (semi)martingale velocity fields to rough flows associated with the corresponding rough driver is done through the It\^o formula, on which one can read the local characteristics of a semimartingale flow.

\bigskip

%%----------------------------------------------%%
\subsection{Notations for function spaces}
\label{SubsectionNotations}
%%----------------------------------------------%%

The study of stochastic flows classically requires the introduction of a number of function spaces, that we recall here.

\medskip

Let $E$ and $F$ be Banach spaces. The derivative of a function $f$ from $E$ to $F$ is understood in the Fr\'echet sense. We shall equip tensor products of Banach spaces with a compatible tensor norm which makes the canonical embedding
\begin{align*}
 \mathcal{L}\big(E,\mathcal{L}(E,F)\big) \hookrightarrow \mathcal{L}\big(E \otimes E,F\big)
\end{align*}
continuous. The $n$-th derivative of $f$ can be seen as a function $D^n f \colon E \to \mathcal{L}\big(E^{\otimes n},F\big)$. For $n \in \NN_0$ and $\rho \in (0,1]$, we define
\begin{align*}
% \| f\|_{n} &:= \| f\|_{\mathcal{C}^n} := \sup_{x \in E} \big\| f(x) \big\| + \sum_{i = 1}^n \sup_{x \in E} \big\| D^i f(x) \big\| \quad \text{and} \\
 \| f\|_{n + \rho} &:= \| f\|_{\mathcal{C}^{n + \rho}} := \sum_{i = 0}^n \sup_{x \in E} \big\| D^i f(x) \big\| + \sup_{0 < \|x - y\| \leq 1} \frac{\big\|D^n f(x) - D^n f(y) \big\|}{\|x - y\|^{\rho}}.
\end{align*}
We define $ \mathcal{C}^{n, \rho}_b(E,F)$ to be the space of $n$-times continuously differentiable functions $f \colon E \to F$ such that $\| f\|_{\mathcal{C}^{n + \rho}} < \infty$.
% We will often just write $\mathcal{C}^{n, \rho}_b$ for $\mathcal{C}^{n, \rho}_b(E,F)$ when the domain and codomain of the function space is clear from the context.

\medskip

Next, we consider the finite dimensional case. Let be $D$ be a domain of $\RR^d$, $A \subseteq D$ a subset, $n \in \NN_0$ and $\rho \in (0,1]$. For a function $f \colon D \to \RR^k$, set
  \begin{align*}
%   \| f\|_{n; A} &:= \sup_{x \in A} \big|f(x)\big| + \sum_{1 \leq |\alpha| \leq n} \sup_{x \in A} \big|D^{\alpha}f(x)\big|, \\
   \| f\|_{n + \rho; A} &:= \sum_{0 \leq |\alpha| \leq n} \sup_{x \in A} \big|D^{\alpha}f(x)\big| + \sum_{|\alpha| = n} \sup_{\stackrel{x,y \in A}{0 < |x - y| \leq 1}} \frac{\big|D^{\alpha} f(x) - D^{\alpha} f(y)\big|}{|x-y|^{\rho}}.
  \end{align*}
We also set % $\| f\|_{n} := \| f\|_{n; D}$ resp. 
$\| f\|_{n + \rho} := \| f\|_{n + \rho; D}$. Note that this is consistent with the notation above when $D = \RR^d$. Let $\mathcal{C}^{n,\rho}(D,\RR^k)$ be the space of $n$-times continuously differentiable functions $f \colon D \to \RR^k$ such that $\| f\|_{n + \rho,K} < \infty$ for every compact subset $K \subset D$. Note that although the (semi-)norms we defined here differ slightly from those used by Kunita in his book \cite{Kun90}, they are actually equivalent on compact sets, hence the spaces coincide. We also define
  \begin{align*}
   \mathcal{C}^{n,\rho}_b(D,\RR^k) := \left\{ f \in \mathcal{C}^{n,\rho}(D,\RR^k) \, :\, \|f\|_{n + \rho} < \infty \right\}.
  \end{align*}
For a function $g \colon D \times D \to \RR^k$, we similarly define
% \begin{align*}
% \| g\|_{n; A}^{\wedge} := \sup_{x,y \in A} \big|g(x,y)\big| + \sum_{1 \leq |\alpha| \leq n} \sup_{x,y \in A} \big|D_x^{\alpha} D^{\alpha}_y g(x,y)\big|
% \end{align*}
% and 
\begin{align*}
&\| g\|_{n + \rho; A}^{\wedge} :=  \sum_{0 \leq |\alpha| \leq n} \sup_{x,y \in A} \big|D_x^{\alpha} D^{\alpha}_y g(x,y)\big|\\
&\quad + \sum_{|\alpha| = n} \sup_{\stackrel{x,y,x',y' \in A}{0 < |x - x'|,|y - y'| \leq 1}} \frac{\Big|D_x^{\alpha} D_y^{\alpha} g(x,y) - D_x^{\alpha} D_y^{\alpha} g(x',y) - D_x^{\alpha} D_y^{\alpha} g(x,y') + D_x^{\alpha} D_y^{\alpha} g(x',y')\Big| }{|x-x'|^{\rho} |y-y'|^{\rho} }.
\end{align*}
As above, set % $\| g\|_{n}^{\wedge} := \| g\|_{n; D}^{\wedge}$ resp. 
$\| g\|_{n + \rho}^{\wedge} := \| g\|_{n + \rho; D}^{\wedge}$.  We denote by $\widehat{\mathcal{C}}^{n,\rho}(D\times D,\RR^k)$ the space of functions $g \colon D \times D \to \RR^k$ which are $n$-times continuously differentiable with respect to each $x$ and $y$ and for which $\| g\|_{n + \rho; K}^{\wedge} < \infty$ for every compact subset $K \subset D$. Set % $\widehat{\mathcal{C}}^{n,\rho} := \widehat{\mathcal{C}}^{n,\rho} (D\times D,\RR)$ and  
\begin{align*}
\widehat{\mathcal{C}}^{n,\rho}_b (D\times D,\RR^k) := \left\{  g \in \widehat{\mathcal{C}}^{n,\rho}(D\times D,\RR^k) \, :\, \| g\|_{n + \rho}^{\wedge} < \infty \right\}.
\end{align*}
We will sometimes use the shorter notation $\mathcal{C}^{n, \rho}$, $\mathcal{C}^{n, \rho}_b$, $\widehat{\mathcal{C}}^{n,\rho}$ resp. $\widehat{\mathcal{C}}^{n,\rho}_b$ when domain and codomain of the function spaces are clear from the context.

\bigskip

%%-------------------------------------------------%%
\subsection{Semimartingale stochastic flows}
\label{SubsectionSemimartingaleFlows}
%%-------------------------------------------------%%

We describe in this Section the basics of the theory of semimartingale stochastic flows, and refer the reader to \cite{LeJanW82} or \cite{Kun90} for a complete account; we refer to Kunita's book for precise regularity and growth assumptions on the different objects involved. Readers familiar with this material can go directly to Section \ref{SubsectionMartingaleDrivers}.

\ssk

Let $\big(\Omega,\mcF,(\mcF_t)_{0\leq t\leq T},\PP\big)$ be a filtered probability space; denote by $\frak{Diff}$, resp. $\frak{F}$, the complete separable metric spaces of $\mcC^{k_0}$ diffeomorphisms of $\RR^d$, resp. $\mcC^{k_0}$ vector fields on $\RR^d$, for some integer $k_0\geq 2$.

\medskip

\begin{definition}   {\sf
A $\frak{Diff}$-valued continuous $\big(\mcF_t\big)_{0\leq t\leq T}$-adapted random process $(\phi_t)_{0\leq t\leq T}$ is called a \textbf{$\frak{Diff}$-valued semimartingale stochastic flow of maps} if the real-valued processes $f\big(\phi_\bullet(x)\big)$ are real-valued $\big(\mcF_t\big)_{0\leq t\leq T}$-semimartingales for all $x\in\RR^d$, and all $f\in\mcC^\infty_c(\RR^d)$. Such a $\frak{Diff}$-valued semimartingale is said to be \textbf{regular} if for every $x,y\in\RR^d$, and $f,g \in\mcC^\infty_c(\RR^d)$, the bounded variation part of $f\big(\phi_\bullet(x)\big)$ and the bracket $\big\langle f\big(\phi_\bullet(x)\big),g\big(\phi_\bullet(y)\big)\big\rangle$ are absolutely continuous with respect to Lebesgue measure $dt$.     }
\end{definition}

\medskip

Their densities $w^f_t(x)$ and $\{f,g\}_t(x,y)$ can be chosen to be jointly measurable and continuous in $f,g$ in the $\mcC^2$-norm \cite{LeJanW82}. Set
$$
\big(\mcL_tf\big)(x) := w^f_t\big(\phi_t^{-1}(x)\big), \qquad \langle f,g\rangle_t(x,y) := \{f,g\}_t\big(\phi_t^{-1}(x),\phi_t^{-1}(y)\big),
$$
so that the processes
$$
M_t^f(x) := f\big(\phi_t(x)\big) - f\big(\phi_0(x)\big) - \int_0^t \big(\mcL_sf\big)\big(\phi_s(x)\big)\,ds, \qquad x\in\RR^d,\,f\in\mcC^\infty_c(\RR^d)
$$
are continuous $\big(\mcF_t\big)_{0\leq t\leq T}$-local martingales with bracket
$$
\big\langle M^f(x), M^g(y) \big\rangle_t = \int_0^t \langle f,g\rangle_t\big(\phi_s(x),\phi_s(y)\big)\,ds.
$$
We have
$$
\big(\mcL_tf\big)(x) = \underset{h\downarrow 0}{\lim}\;\EE\left[\frac{f\big(\phi_{t+h,t}(x)\big) - f(x)}{h}\bigg| \mcF_t \right]
$$
and 
$$
\langle f,g\rangle (x,y) = \underset{h\downarrow 0}{\lim}\;\frac{1}{h}\,\EE\Big[\big\{f\big(\phi_{t+h,t}(x)\big)-f(x)\big\}\big\{g\big(\phi_{t+h,t}(y)\big)-g(y)\big\}\Big| \mcF_t\Big],
$$
with limits in $L^1$ whenever they exist. Under proper regularity conditions \cite{LeJanW82}, the operators $\langle f,g\rangle_s(x,y)$ can be seen to be random differential operators of the form
$$
\langle f,g\rangle_s(x,y) = A^{ij}_s(x,y)\,\partial^2_{x_iy_j},
$$
for some process $A_s(x,y)$ with values in the space ${\sf Symm}(d)$ of symmetric $d\times d$ matrices. The operators $\mcL_s$ can moreover be expressed in terms of $A_s$ and its differential with respect to the space variables, so that\textit{ the data of the processes $A_s$ and $\mcL_s$ is equivalent to the data of the process $A_s$ and an $\frak{F}$-valued process $b_\bullet$}. The family of random operators $\langle \cdot,\cdot\rangle_t$ and \textit{the drift} $b_t$ are called the \textbf{local characteristics} of the $\frak{Diff}$-valued semimartingale $\phi_\bullet$.  As an example, for the semimartingale flow generated by a stochastic differential equation of the form
$$
dx_t = V_i(x_t)\,{\circ d}B_t^i,
$$
driven by an $\ell$-dimensional Brownian motion, we have
$$
\mcL_tf = \frac{1}{2}\,\sum_{i=1}^\ell V_i^2f
$$
and 
$$
\langle f,g\rangle_s(x,y) = \big(V_if\big)(x)\big(V_jg\big)(y),
$$
and the drift $b_s$ in the local characteristic is given here by the time-independent vector field
$$
b_s(x) = \frac{1}{2}\big(V_iV_i\big)(x) = \frac{1}{2}\big(D_xV_i\big)V_i(x).
$$

\medskip

The infinitesimal counterpart of a $\frak{Diff}$-valued semimartingale is given by the following notion. 

\medskip

\begin{definition}  {\sf
A \textbf{semimartingale velocity field} is an $\frak{F}$-valued process $(V_t)_{0\leq t\leq T}$ such that the processes $\big(V_\bullet f\big)(x)$ are real-valued semimartingale for all $x\in\RR^d$ and all $f\in\mcC^\infty_c(\RR^d)$. It is called {\bf regular} if one can write 
$$
V_t = M_t + \int_0^t v_s\,ds
$$
for a vector field-valued adapted process $v_\bullet$, and an $\frak{F}$-valued local martingale $M_\bullet$ for which there exists a ${\sf Symm}(d)$-valued process $a_s(x,y)$ with 
$$
\big\langle \partial_x^\alpha M_\bullet(x), \partial_y^\beta M_\bullet(y)\big\rangle_t = \int_0^t \partial_x^\alpha\,\partial_y^\beta\, a_s(x,y)\,ds
$$
for a range of multiindices $\alpha,\beta$ depending on the regularity assumptions on $a_s$. The pair $(a_\bullet,v_\bullet)$ is called the \textbf{local characteristics} of the semimartingale velocity field $V_\bullet$.   }
\end{definition}

\medskip

A theory of Stratonovich integration can be constructed for making sense of integrals of the form 
$$
\int_0^t V_{\circ ds}(x_s),
$$
for some progressively measurable process $x_\bullet$ and some regular semimartingale velocity field $V_\bullet$, as a limit in probability of symmetric Riemann sums. This requires some almost sure regularity properties on the local characteristics $(a_\bullet,v_\bullet)$ of $V_\bullet$, and some almost sure bound on $\int_0^t \big|a_s(x_s,x_s)\big|\,ds$ and $\int_0^t\big|v_s(x_s)\big|\,ds$ -- see e.g. Section 2.3 of \cite{KunTata}. Under these conditions, the integral Stratonovich equation
\begin{equation}
\label{EqKunitaSDE}
x_t = x_0 + \int_0^t V_{\circ ds}(x_s)
\end{equation}
can be seen to have a unique solution started from any point $x_0\in\RR^d$. 

\medskip

\begin{thm}[\cite{LeJanW82}]
\label{ThmLeJanW82}
These solutions can be gathered into a semimartingale stochastic flow whose local characteristics are $(a_\bullet,\,v_\bullet+c_\bullet)$, where the time-dependent vector field $c_s$ has coordinates
$$
c^i_s(x) := \frac{1}{2} \sum_{j=1}^d\partial_{y_j}a^{ij}_s(x,y)_{\big| y=x}
$$
in the canonical basis of $\RR^d$. Conversely, one can associate to any regular stochastic flow of diffeomorphisms $\phi_\bullet$ a semimartingale velocity field $V_\bullet$, with the same local characteristics as $\phi_\bullet$, and such that $\phi_\bullet$ coincides with the stochastic flow generated by the Stratonovich equation 
$$
x_t = x_0 + \int_0^t V_{\circ ds}(x_s) - \int_0^t c_s(x_s)\,ds.
$$
\end{thm}

\medskip

The optimal regularity assumptions on the velocity fields and stochastic flows of maps are given in Theorems 4.4.1 and 4.5.1 of Kunita's book. We shall use the full strength of these two statements in Section \ref{SubsectionStochRoughFlows} to identify semimartingale stochastic flows of maps and the rough flows associated with the lift of the semimartingale velocity fields into a rough driver.

\ssk

The correspondence between semimartingale stochastic flows and semimartingale velocity fields via an \textit{It\^o equation} of the form 
$$
x_t = x_0 + \int_0^t V_{ds}(x_s)
$$
is exact, with no need to add the drift $\int_0^t c_s\,ds$. We state it here under the above form as we shall see below that rough flows are naturally associated with Stratonovich differential equations.

\ssk

The main difficulty in this business is to deal with the local martingale part of the dynamics, which is where probability theory is really needed. As a consequence, we shall concentrate our efforts on local martingale velocity fields in the sequel, the remaining changes to deal with regular semimartingale velocity fields being essentially cosmetic. As above, we shall freely identify in the sequel vector fields with first order differential operators. In order to keep consistent notations, we shall also denote by
$$
\int_0^t \alpha_s\, m_{ds} \qquad \textrm{and}\qquad \int_0^t \alpha_s\, m_{\circ ds}
$$
the It\^o and Stratonovich integrals of an adapted process $\alpha_s$ with respect to a local martingale $m_s$.

\bigskip

%%-----------------------------------------------%%
\subsection{Local martingale rough drivers}
\label{SubsectionMartingaleDrivers}
%%-----------------------------------------------%%

The aim of this section is to give conditions under which a local martingale velocity field can be lifted into a rough driver. Let $D$ be an open connected subset of $\RR^d$ and let $M$ stand for a local martingale velocity field. We prove in this Section that such a field can be lifted to a rough driver ${\bf M} = \big(M_{ts},\bbM_{ts}\big)_{0\leq s\leq t\leq T}$, with $M_{ts} := M_t - M_s$, under regularity and boundedness assumptions on the local characteristic of $M$. At a heuristic level, if $M$ is differentiable in space, the second level operator $\bbM_{ts}$ associated with $M_{ts} = M_t - M_s$ is given by the formula
\begin{equation*}
\begin{split}
\bbM_{ts} &= \int_s^t M_{us}M_{\circ du} \\
                 &= \left(\int_s^tM^i_{us}\partial_i M^k_{\circ du}\right)\partial_k + \left(\int_s^t M^j_{us}\,M^k_{\circ du}\right)\partial_{jk}^2,
\end{split}
\end{equation*}
with obvious notations for the operators $\partial_k$ and $\partial^2_{jk}$. In the following, we will use the notation
\begin{align*}
 \big(M_{ts}.M_{ts}\big) := \big(DM_{ts}\big)(M_{ts}).
\end{align*}
As the classical rules of Stratonovich integration give
$$
\left(\int_s^t M^j_{us}\,M^k_{\circ du}\right)\partial_{jk}^2 = \frac{1}{2}\,M^j_{ts}M^k_{ts}\partial_{jk}^2 = \frac{1}{2}\,M_{ts}M_{ts} - \frac{1}{2}\,\big(M_{ts}.M_{ts}\big),
$$
we see that the operators $\bbM_{ts}$ can be decomposed as
\begin{align}\label{eqn:def_second_order_process}
\bbM_{ts} = W_{ts} + \frac{1}{2}\,M_{ts}M_{ts},
\end{align}
where
\begin{equation*}
\begin{split}
W_{ts} &= \left(\int_s^tM^i_{us}\partial_i M^k_{\circ du}\right)\partial_k - \frac{1}{2}\,\big(M_{ts}.M_{ts}\big) \\
            &= \frac{1}{2}\,\left(\int_s^t M^i_{us}\partial_i M^k_{du} - M^i_{du}\partial_i M^k_{us} \right)\partial_k \\
            &= \frac{1}{2}\,\int_s^t \big(M_{us}.M_{du} - M_{du}.M_{us}\big) 
\end{split}
\end{equation*}
is a martingale velocity field defined pointwisely by an It\^o integral. The proof that the process $\mathbf{M}_{ts} := (M_{ts},\mathbb{M}_{ts})$ has a modification which is a $p$-rough driver for every $2<p<3$ will require two elementary intermediate results that heavily rest on a now classical modified version of Kolmogorov's regularity criterion that we recall here for the reader's convenience -- this is different from the content of the above Section \ref{SubsectionKolmogorovTheorem}; they can be found in Section 3 of Kunita's book \cite{Kun90}. 

\medskip

\begin{thm}   {\sf 
\label{ThmKunitaKolmogorov}
\begin{enumerate}
   \item Let $M_\bullet(x),\,x\in D$, be a family of continuous local martingales started from $0$, such that the joint quadratic variation $\big\langle M_\bullet(x),M_\bullet(y)\big\rangle$ has a continuous modification in $\widehat\mcC^{0,\delta}$. Then the process $M$ has a modification that is a continuous process with values in $\mcC^{0,\epsilon}$, for every $\epsilon<\delta$. \vspace{0.1cm}
   
   \item Let $M$ be a local martingale velocity field started from $0$, with local characteristic the random field $a_t(x,y)$. Assume that $a$ has a continuous modification that belongs almost surely to ${\widehat\mcC}^{m,\delta}$, for some integer $m\geq 1$, and $0<\delta\leq 1$. Then, for every $0<\epsilon<\delta$, the velocity field $M$ has a modification that is a continuous process with values in $\mcC^{m,\epsilon}$; we still denote it by $M$. Furthermore, for each multi-index $\alpha$, with $|\alpha|\leq m$, the time varying random field $\partial_x^\alpha M$ is a local martingale velocity field with quadratic variation 
$$
d\big\langle \partial_x^\alpha M_\bullet(x), \partial_x^\alpha M_\bullet(y) \big\rangle_t = \partial_x^\alpha\partial_y^\alpha a_t(x,y)\, dt. 
$$  

   \item Let here $M$ and $N$ be two local martingale velocity fields with values in $\mcC^{m,\delta}$. Then their joint quadratic variation 
$$
\big\langle M_\bullet(x),N_\bullet(y)\big\rangle_t
$$   
has a continuous modification taking values in $\widehat\mcC^{m,\epsilon}$ for every $\epsilon<\delta$. Furthermore, if $m\geq 1$, this modification satisfies the identity
$$
\partial_x^\alpha \partial_y^\beta \big\langle M_\bullet(x),N_\bullet(y)\big\rangle_t = \big\langle \partial_x^\alpha M_\bullet(x), \partial_y^\beta N_\bullet(y)\big\rangle_t,
$$
for all $|\alpha|,\,|\beta|\leq m$.
\end{enumerate}   }
\end{thm}

\ssk  

\begin{Dem}
 Cf. Theorem 3.1.1, Theorem 3.1.2 and Theorem 3.1.3 in \cite{Kun90}.
\end{Dem}

\ssk

These regularity results will be instrumental in the proof of the following intermediate result. Note that $N$ in equation \eqref{EqIntMN} below is seen as a vector field, not a differential operator, so $\big(MN\big)(x) = \big(D_xN\big)\big(M(x)\big)$.
  
\ssk  
  
\begin{prop}\label{prop:diff_and_prod_rule_stoch_integral}   {\sf
Let $M,N \colon D \times [0,T] \to \RR$ be continuous $\mathcal{C}^{m,\delta}$-valued local martingale fields, for $m \in \NN_0$ and $\delta \in (0,1]$. Assume $M$ is adapted to the filtration generated by $N$. Then the pointwisely defined It\^o integral
\begin{align} \label{EqIntMN}
  t \mapsto \int_0^t \big(M_s\, N_{ds}\big)(x)
\end{align}
has a continuous modification taking values in $\mathcal{C}^{m,\alpha}$ process for every $\alpha < \delta$. Moreover, if $m \geq 1$, the derivative is almost surely given by the formula
\begin{align} \label{eqn:prod_rule_ito_integral}
\partial_{x_i} \left( \int_0^t M_s\, N_{ds} \right) = \int_0^t \partial_{x_i} M_s\, N_{ds} + \int_0^t M_s\,\partial_{x_i} N_{ds}.
\end{align}
Both assertions also hold for the Stratonovich integral.   }
\end{prop}

\ssk

The proof of this result is somewhat lengthy but rests on classical considerations based on the regularization theorem \ref{ThmKunitaKolmogorov}.

\ssk

\begin{Dem}
Set 
$$
U_t(x) := \left(\int_0^t M_s\, N_{ds}\right)(x).
$$ 
This is a continuous local martingale field with joint quadratic variation given by
\begin{align}\label{eqn:joint_quadr_var_iter_ito}
\big\langle U_{\bullet}(x),U_{\bullet}(y) \big\rangle_t = \int_0^t M_s(x) \, M_s(y)\, d\big\langle N_\bullet(x), N_\bullet(y) \big\rangle_s.
\end{align}
Fix $t \in [0,T]$, some compact set $K \subset D$ and set 
$$
g(x,y) := \langle U_\bullet(x), U_\bullet(y) \rangle_t.
$$
We first consider the case $m = 0$. Choose $\alpha < \delta' < \delta$. Then, for $x,x',y,y' \in K$, we have
\begin{align*}
&\Big|g(x,y) - g(x',y) - g(x,y') + g(x',y')\Big| \\
=\ &\Big|\big\langle U_\bullet(x) - U_\bullet(x'), U_\bullet(y) - U_\bullet(y') \big\rangle_t\Big| \\
\leq\ &\left| \int_0^t \big(M_s(x) - M_s(x')\big)\,\big(M_s(y) - M_s(y')\big)\, d\big\langle N_\bullet(x,),N_\bullet(y) \big\rangle_s \right| \\
&+ \left| \int_0^t \big(M_s(x) - M_s(x')\big) M_s(x')\, d\big\langle N_\bullet(x),  N_\bullet(y) - N_\bullet(y') \big\rangle_s\right| \\
&+ \left| \int_0^t M_s(x')\big(M_s(y) - M_s(y')\big)\, d\big\langle N_\bullet(x) - N_\bullet(x'), N_\bullet(y) \big\rangle_s \right| \\
&+ \left| \int_0^t M_s(x')M_s(y')\, d\big\langle N_\bullet(x) - N_\bullet(x'), N_\bullet(y) - N_\bullet(y') \big\rangle_s \right|.
\end{align*}
For the first integral, we use Kunita's extended Cauchy-Schwarz inequality, as stated in \cite[Theorem 2.2.13]{Kun90}, to see that
\begin{align*}
&\left| \int_0^t \big(M_s(x) - M_s(x')\big)\,\big(M_s(y) - M_s(y')\big)\, d\big\langle N_\bullet(x), N_\bullet(y)\big\rangle_s \right| \\
\leq\ &\left( \int_0^t \big(M_s(x) - M_s(x')\big)^2\, d\big\langle N_\bullet(x) \big\rangle_s \right)^{\frac{1}{2}} \left( \int_0^t (M_s(y) - M_s(y'))^2\, d\big\langle N_\bullet(y) \big\rangle_s \right)^{\frac{1}{2}} \\
\leq\ &|x - x'|^{\delta} |y - y'|^{\delta} \sup_{s \in [0,T]} \|M_s \|_{\delta-\text{H\"ol}; K}^2 \sup_{z \in K} \big\langle N_\bullet(z) \big\rangle_T.
\end{align*}
Similarly, for the second integral,
\begin{align*}
&\left| \int_0^t \big(M_s(x) - M_s(x')\big) M_s(x')\, d\big\langle N_\bullet(x),  N_\bullet(y) - N_\bullet(y') \big\rangle_s\right| \\
\leq\ &\left( \int_0^t \big(M_s(x) - M_s(x')\big)^2\, d\big\langle N_\bullet(x)\big\rangle_s \right)^{\frac{1}{2}} \left( \int_0^t M_s(x')^2\, d\big\langle N_\bullet(y) - N_\bullet(y')\big\rangle_s \right)^{\frac{1}{2}} \\
\leq\ &|x - x'|^{\delta} |y - y'|^{\delta'} \sup_{s \in [0,T]} \big\|M_s \big\|_{\infty;K}  \sup_{s \in [0,T]} \big\|M_s \big\|_{\delta-\text{H\"ol};K}  \sup_{z \in K} \big\langle N_\bullet(z) \big\rangle_T^{\frac{1}{2}} \Big\| \big\langle N_\bullet, N_\bullet \big\rangle_T \Big\|_{\delta';K}^\frac{1}{2}.
\end{align*}
From point 3 in Theorem \ref{ThmKunitaKolmogorov} we know that there is a version of the joint quadratic variation of $N$ such that $\big\| \langle N_{\bullet}, N_{\bullet} \rangle_T \big\|_{\delta';K}^{\wedge} < \infty$. The other integrals are estimated similarly. This shows that
\begin{align*}
\sup_{\stackrel{x,x',y,y' \in K}{x \neq x', y\neq y'}} \frac{\big|g(x,y) - g(x',y) - g(x,y') + g(x',y')\big|}{|x - x'|^{\delta'}\,|y - y'|^{\delta'}} < \infty.
\end{align*}
Clearly $\|g\|_{\infty;K} < \infty$, thus we have shown that the joint quadratic variation of $U$ has a modification which is a continuous $\widehat{\mathcal{C}}^{0,\delta'}$-process. Point 1 in Theorem \ref{ThmKunitaKolmogorov} shows that $U$ has a modification which is a continuous $\mathcal{C}^{0,\alpha}$-process. Now let $m \geq 1$.  From  point 3 in Theorem \ref{ThmKunitaKolmogorov} we may deduce that the joint quadratic variation of $N$ has a modification which is a continuous $\tilde{\mathcal{C}}^{m,\delta'}$-process with
\begin{align*}
\partial_x^{\beta} \partial_y^{\gamma} \big\langle N_\bullet(x),N_\bullet(y) \big\rangle = \big\langle \partial_x^{\beta} N_\bullet(x), \partial_y^{\gamma} N_\bullet(y) \big\rangle
\end{align*}
for every $|\beta|, |\gamma| \leq m$. We may apply Proposition \ref{prop:diff_smooth_case} in Appendix iteratively in equation \eqref{eqn:joint_quadr_var_iter_ito} to show that $\big\langle U_\bullet(x),U_\bullet(y) \big\rangle_t$ has a modification which is $m$-times differentiable with respect to $x$ and $y$, and we can calculate the derivatives using the product rule stated in Proposition \ref{prop:diff_smooth_case}. As above, one can show that the $m$-th derivative has the claimed H\"older regularity, and we can conclude with point 2 of Theorem \ref{ThmKunitaKolmogorov} that $U$ has a modification which is a continuous $\mathcal{C}^{m,\alpha}$-process. The It\^o-Stratonovich conversion formula
\begin{align*}
\int_0^t M_s\, N_{\circ ds} = \int_0^t M_s\, N_{ds} + \frac{1}{2} \big\langle M_\bullet, N_\bullet \big\rangle_t
\end{align*}
and point 3 in Theorem \ref{ThmKunitaKolmogorov} show that the same is true for the Stratonovich integral. 

\ssk

We now come to equation \eqref{eqn:prod_rule_ito_integral}. In the following, we use $\| \cdot \|_{L^1}$ for the $\LL^1$-norm with respect to $\PP$. For $n \in \NN$, set
\begin{align*}
\tau_n = \inf \Big\{ t \in [0,T]\,:\, \big\|M_t \big\|_{\mathcal{C}^{1,\delta}} + \big\|N_t\big\|_{\mathcal{C}^{1,\delta}} \geq n \Big\}. 
\end{align*}
The random times $\tau_n$ define an increasing sequence of stopping times such that $\PP(\tau_n < T) \to 0$ for $n \to 0$. Let $x \in D$, $t > 0$ and choose $h$ such that $x + h e_i \in D$. Then
\begin{align*}
&\left\|\frac{1}{h} \left\{ \int_0^{t \wedge \tau_n} \big(M_s\, N_{ds}\big)(x + h e_i) - \int_0^{t \wedge \tau_n} \big(M_sN_{ds}\big)(x) \right\} \right. \\
&\qquad\qquad\qquad\qquad\qquad\qquad\qquad\qquad- \int_0^{t \wedge \tau_n} \big(\partial_i M_s\, N_{ds}\big)(x) - \int_0^{t \wedge \tau_n} \big(M_s\,\partial_i N_{ds}\big)(x) \Big\|_{L^1} \\
&= \left\| \int_0^{t \wedge \tau_n} \frac{M_s(x + h e_i) - M_s(x)}{h} \, N_{ds}(x + h e_i) + \int_0^{t \wedge \tau_n} M_s(x)\, \frac{N_{ds}(x + h e_i) - N_{ds}(x)}{h} \right. \\
&\quad-\int_0^{t \wedge \tau_n} \partial_i M_s(x)\, N_{ds}(x) + \int_0^{t \wedge \tau_n} M_s(x)\, \partial_i N_{ds}(x) \Big\|_{L^1} \\
&\leq\left\|  \int_0^{t \wedge \tau_n} \left(\frac{M_s(x + h e_i) - M_s(x)}{h}\right) \, \big(N_{ds}(x + h e_i) - N_{ds}(x)\big) \right\|_{L^1} \\
&\quad+\left\|  \int_0^{t \wedge \tau_n} \left(\frac{M_s(x + h e_i) - M_s(x)}{h} - \partial_i M_s(x)\right) \,  N_{ds}(x) \right\|_{L^1} \\
&\quad+\left\|  \int_0^{t \wedge \tau_n} M_s(x)\, \left(\frac{N_{ds}(x + h e_i) - N_{ds}(x)}{h} - \partial_i N_{ds}(x) \right) \right\|_{L^1}.
\end{align*}
We aim to show that the integrals on the right hand side vanish for $h \to 0$, using the Burkholder-Davis-Gundy inequality. We start with the first integral. Note that since $N$ is a $\mathcal{C}^{1,\delta}$ process, also the stopped process is a $\mathcal{C}^{1,\delta}$ process and its joint quadratic variation has a modification which is a $\tilde{\mathcal{C}}^{1,\delta'}$ process for any $\delta' < \delta$. In particular,
\begin{align*}
\big\langle N_\bullet(x + h e_i) - N_\bullet(x) \big\rangle_{t \wedge \tau_n} \to 0
\end{align*}
almost surely for $h \to 0$. From the Burkholder-Davis-Gundy inequality,
\begin{align*}
\EE \Big|\big\langle N_\bullet(x + h e_i) - N_\bullet(x) \big\rangle_{t \wedge \tau_n}\Big|^{p/2} \leq C_p \, \EE \Big| \sup_{s \in [0, t \wedge \tau_n]} \big|N_s(x + h e_i) - N_s(x)\big|\Big|^p \leq C \, 2^p n^p
\end{align*}
which shows that $\big\langle N_\bullet(x + h e_i) - N_\bullet(x)\big\rangle_{t \wedge \tau_n} \to 0$ in $L^p$ for any $p \geq 1$. Using the Burkholder-Davis-Gundy inequality for the first integral gives
\begin{align*}
&\EE \left| \int_0^{t \wedge \tau_n} \left(\frac{M_s(x + h e_i) - M_s(x)}{h}\right) \, \big(N_{ds}(x + h e_i) - N_{ds}(x)\big) \right| \\
\leq\ &C \, \EE \left| \int_0^{t \wedge \tau_n} \left(\frac{M_s(x + h e_i) - M_s(x)}{h}\right)^2 \, d\langle N_\bullet(x + h e_i) - N_\bullet(x) \rangle_s \right|^{\frac{1}{2}} \\
\leq\ &C \, \EE \left| \sup_{s \in [0, t \wedge \tau_n]} \big\|\partial_{x_i} M_s\big\|_{\infty}^2 \big\langle N_\bullet(x + h e_i) - N_\bullet(x) \big\rangle_{t \wedge \tau_n} \right|^{\frac{1}{2}} \\
\leq\ &C n \, \EE \left| \big\langle N_\bullet(x + h e_i) - N_\bullet(x) \big\rangle_{t \wedge \tau_n} \right|^{\frac{1}{2}} \to 0
\end{align*}
for $h \to 0$. For the second integral, the Burkholder-Davis-Gundy inequality gives
\begin{align*}
&\EE \left| \int_0^{t \wedge \tau_n} \left(\frac{M_s(x + h e_i) - M_s(x)}{h} - \partial_i M_s(x)\right) \,  N_{ds}(x) \right| \\
\leq\ &C \, \EE \left| \int_0^{t \wedge \tau_n} \left(\frac{M_s(x + h e_i) - M_s(x)}{h} - \partial_i M_s(x)\right)^2 \, d\big\langle N_\bullet(x)\big\rangle_s \right|^{\frac{1}{2}}
\end{align*}
For $h \to 0$, we can use dominated convergence twice for the expectation and the Lebesgue-Stieltjes integral to see that this term indeed converges to $0$. Now we come to the third integral. Set
\begin{align*}
N^h_t(x) := \frac{N_t(x + h e_i) - N_t(x)}{h}.
\end{align*}
As before, we can use the Burkholder-Davis-Gundy inequality to see that
\begin{align*}
\EE \left|  \int_0^{t \wedge \tau_n} M_s(x)\, \left( N^h_{ds}(x) - \partial_i N_{ds}(x) \right) \right| \leq C n \EE \left| \big\langle N^h_\bullet(x) - \partial_i N_\bullet(x) \big\rangle_{t \wedge \tau_n} \right|^{\frac{1}{2}}
\end{align*}
Taking a modification which is a $\tilde{\mathcal{C}}^{1,\delta'}$ process of the joint quadratic variation of $N$ gives, using point 3 of Theorem \ref{ThmKunitaKolmogorov}, that 
\begin{align*}
\big\langle N^h_\bullet(x) - \partial_i N_\bullet(x) \big\rangle_{t \wedge \tau_n} \to 0
\end{align*}
almost surely as $h \to 0$. Using the Burkholder-Davis-Gundy inequality as above shows that the convergence also holds in $L^p$ for any $p\geq 1$. To conclude, we have shown that
\begin{align*}
\frac{1}{h} \left( \int_0^{t \wedge \tau_n} M_s(x + h e_i) \, N_{ds}(x + h e_i) - \int_0^{t \wedge \tau_n} M_s(x)\, N_{ds}(x) \right)\\
\to \int_0^{t \wedge \tau_n} \partial_{x_i} M_s(x)\, N_{ds}(x) + \int_0^{t \wedge \tau_n} M_s(x)\, \partial_{x_i} N_{ds}(x)
\end{align*}
in $L^1$ for $h \to 0$. Since we already now that the convergence holds almost surely, the limits coincide and we have shown that
\begin{align}\label{eqn:prod_rule_stopping_times}
\partial_{x_i} \left( \int_0^{t \wedge \tau_n} M_s(x)\, N_{ds}(x) \right) = \int_0^{t \wedge \tau_n} \partial_{x_i} M_s(x)\, N_{ds}(x) + \int_0^{t \wedge \tau_n} M_s(x)\, \partial_{x_i} N_{ds}(x)
\end{align}
holds almost surely for every $n \in \NN$. Now,
\begin{align*}
\PP\left( \partial_{x_i} \left( \int_0^{t} M_s(x)\, N_{ds}(x) \right) \neq \partial_{x_i} \left( \int_0^{t \wedge \tau_n} M_s(x)\, N_{ds}(x) \right) \right) \leq \PP (\tau_n < T) \to 0
\end{align*}
as $n \to \infty$, therefore
\begin{align*}
\partial_{x_i} \left( \int_0^{t \wedge \tau_n} M_s(x)\, N_{ds}(x) \right) \to \partial_{x_i} \left( \int_0^{t} M_s(x)\, N_{ds}(x) \right)
\end{align*}
in probability as $n \to \infty$. The same is true for the other two integrals on the right hand side of \eqref{eqn:prod_rule_stopping_times} which shows \eqref{eqn:prod_rule_ito_integral} for the It\=o integral. For the Stratonovich integral, the assertion follows from the It\=o-Stratonovich conversion formula and the equality
\begin{align*}
\partial_{x_i} \big\langle M_\bullet(x) ,N_\bullet(x) \big\rangle_t = \big\langle \partial_{x_i} M_\bullet(x), N_\bullet(x)\big\rangle_t + \big\langle M_\bullet(x), \partial_{x_i} N_\bullet(x)\big\rangle_t.
\end{align*}
 \end{Dem}

\medskip

% Let $\delta \in (0,1]$. For a function $g \colon D \times D \to \RR^m$, set 
% \begin{align*}
%  \| g \|^{\wedge}_{\delta;\text{loc}} := \sup_{\stackrel{x,y,x',y' \in D}{0 < |x - x'|,|y - y'| \leq 1}} \frac{|g(x,y) - g(x',y) - g(x,y') + g(x',y')|}{|x - x'|^{\delta} |y - y'|^{\delta}} .
% \end{align*}
% 
% Note that $\| g \|^{\wedge}_{\delta;\text{loc}} < \infty$ implies that $g$ is continuous and that
% \begin{align*}
%   \sup_{\stackrel{x,y,x',y' \in K}{0 < |x - x'|,|y - y'| \leq 1}} \frac{|g(x,y) - g(x',y) - g(x,y') + g(x',y')|}{|x - x'|^{\delta} |y - y'|^{\delta}} < \infty
% \end{align*}
% for every compact set $K \subseteq \RR^d$, i.e. $g \in \hat{\mathcal{C}}^{0,\delta}$. Note also that if $g$ has compact support, $g \in \hat{\mathcal{C}}^{0,\delta}$ is equivalent to $\| g \|^{\wedge}_{\delta;\text{loc}} < \infty$.

Recall the definition of the space $\hat{\mathcal{C}}^{n,\delta}_b$.
% $\| \cdot \|^{\wedge}_{\delta}$ given in the appendix. For a function $a \colon D \times D \to \RR$, we will use the notation
% \begin{align*}
%  \gimel a(x,y) := a(x,x) - a(x,y) - a(y,x) + a(y,y).
% \end{align*}

 \ssk

\begin{prop}\label{prop:ex_mod_ito_int}   {\sf
Let $M,N \colon D \times [0,T] \to \RR$ be continuous local martingale fields adapted to the same filtration. Assume that the quadratic variation of the processes is given by
\begin{align*}
 d \langle M_{\bullet}(x), M_{\bullet}(y) \rangle_t &= a_t(x,y)\, dt\quad \text{resp.} \\
 d \langle N_{\bullet}(x), N_{\bullet}(y) \rangle_t &= b_t(x,y)\, dt
\end{align*}
for every $x,y \in D$ and every $t \in [0,T]$. Moreover, assume that there is a $\delta \in (0,1]$ such that $a$ and $b$ have continuous modifications in the space $\hat{\mathcal{C}}^{0,\delta}_b$.
%for which $\| a_t \|^{\wedge}_{\delta}$ and  $\| b_t \|^{\wedge}_{\delta}$ are finite for every $t \in [0,T]$. 
Let $p > 2$ and $\rho \in (0,\delta)$. 

\vspace{0.1cm}

\begin{itemize}
 \item[(i)] 
%   Assume that there is some compact set $K$ in $D$ such that $a_t$ is supported on $K \times K$ for every $t \in [0,T]$ almost surely, or that 
  Assume that there is a constant $\kappa > 0$ and some $\eta \in (0,\infty)$ such that
% %  \begin{align*}
% %   \sup_{\stackrel{x,y \in D}{0 < |x-y| \leq 1}} \frac{\| \sqrt{ | \gimel a(x,y)|} \|_{L^q}}{|x - y|^{\delta}} \leq \kappa
% %  \end{align*}
% %   for some $\delta \in (0,1]$ and
\begin{align}\label{eqn:growth_loc_char_a}
  \sup_{u \in [0,T]} \left\| \sqrt{a_u (x,x)} \right\|_{L^q(\Omega)} \leq \frac{\kappa}{1 + |x|^{\eta}}
\end{align}
  holds for every $x \in D$ where
\begin{align}\label{eqn:q_large_enough_level_1}
 q > \max \left\{ \frac{\delta}{\eta(\delta - \rho)} + d \left(\frac{\delta}{\eta(\delta - \rho)} + \frac{1}{\delta - \rho} \right), \frac{1}{\frac{1}{2} - \frac{1}{p}} \right\}.
\end{align}
Then the process $M$ has a modification which satisfies
\begin{align*}
 \sup_{0 \leq s < t \leq T} \frac{\| M_{ts} \|_{\mathcal{C}^{\rho}}}{|t-s|^{\frac{1}{p}}} < \infty
\end{align*}
almost surely. \vspace{0.1cm}

  \item[(ii)] Set
  \begin{align*}
    U_{ts}(x) := \int_s^t (M_{us} N_{du})(x).
  \end{align*}
% Assume either that 
%  there is a constant $\kappa > 0$ such that
%  \begin{align}\label{eqn:quadr_var_reg}
%   \sup_{\stackrel{x,y \in D}{0 < |x-y| \leq 1}} \frac{\| \sqrt{ | \gimel a(x,y)|} \|_{L^q}}{|x - y|^{\delta}} + \sup_{\stackrel{x,y \in D}{0 < |x-y| \leq 1}} \frac{\| \sqrt{ | \gimel b(x,y)|} \|_{L^q}}{|x - y|^{\delta}} \leq \kappa
%  \end{align}
%   for some $\delta \in (0,1]$ and
% $a$ and $b$ are compactly supported in the sense above, or 
    Assume that there is a constant $\kappa > 0$ and some $\eta \in (0,\infty)$ such that
\begin{align}\label{eqn:quadr_var_growth}
  \sup_{u \in [0,T]} \left\| \sqrt{a_u (x,x)} \right\|_{L^q(\Omega)} + \sup_{u \in [0,T]} \left\| \sqrt{b_u (x,x)} \right\|_{L^q(\Omega)} \leq \frac{\kappa}{1 + |x|^{\eta}}
\end{align}
 for every $x \in D$ with $q$ sufficiently large as in \eqref{eqn:q_large_enough_level_1}. 

Then the process $U$ 
has a modification which satisfies
\begin{align*}
 \sup_{0 \leq s < t \leq T} \frac{\big\| U_{ts} \big\|_{\mathcal{C}^{\rho}}}{|t-s|^\frac{2}{p}} < \infty
\end{align*}
almost surely. \vspace{0.1cm} 
 
  \item[(iii)] Assume that $\sup_{t \in [0,T]} \|a_t\|_{\delta}^{\wedge}$ and $\sup_{t \in [0,T]} \|b_t\|_{\delta}^{\wedge}$ are almost surely bounded random variables
  and that \eqref{eqn:quadr_var_growth} holds for some $\kappa > 0$ and $\eta \in (0,\infty)$ uniformly for every $q \geq 1$. 
  
  Then the random variables
  \begin{align*}
   \sup_{0 \leq s < t \leq T} \frac{\big\| M_{ts} \big\|_{\mathcal{C}^{\rho}}}{|t-s|^{\frac{1}{p}}} \qquad \text{and} \qquad \sqrt{\sup_{0 \leq s < t \leq T} \frac{\| U_{ts} \|_{\mathcal{C}^{\rho}}}{|t-s|^\frac{2}{p}}}
  \end{align*}
  have Gaussian tails.
\end{itemize}   }
\end{prop}

\ssk

\begin{Dem}
 We start with (i). %Assume that $a$ is supported on some compact set $K$ as stated. 
 In a first step, we assume that $\sup_{t \in [0,T]} \|a_t\|_{\delta}^{\wedge}$ is an almost surely bounded random variable and that \eqref{eqn:growth_loc_char_a} holds uniformly for every $q \geq 1$. Under these assumptions, we aim to show that $M$ has a modification such that 
  \begin{align}\label{eqn:1st_level_fin_and_gauss_tails}
   \sup_{0 \leq s < t \leq T} \frac{\| M_{ts} \|_{\mathcal{C}^{\rho}}}{|t-s|^{\frac{1}{p}}}
  \end{align}
  is finite almost surely and has Gaussian tails. Set
  \begin{align}\label{eqn:norm_quadr_var_finite}
   \alpha_1 := \left\| \sup_{u \in [0,T] } \| a_u \|^{\wedge}_{\delta} \right\|_{L^{\infty}}. % = \sup_{u \in [0,T] } \| a_u \|^{\wedge}_{\delta;K}
  \end{align}
  %which is an almost surely constant and, by continuity, finite random variable. Abusing notation, we denote its value by the same symbol. 
  Let $s <t$ and $x,y \in D$ with $0 < |x - y| \leq 1$. By the Burkholder Davis Gundy inequality with optimal constants (cf. \cite[Theorem A]{CK91}), for every $q \geq 2$,
  \begin{align*}
   \EE | M_{ts}(x) - M_{ts}(y) |^q &\leq C^q q^{q/2} \EE \left| \int_s^t a_u(x,x) - a_u(x,y) - a_u(y,x) + a_u(y,y)\, du \right|^{q/2}  \\
   &\leq C^q q^{q/2} |t-s|^{q/2} \alpha_1^{q/2} |x - y|^{q \delta}
  \end{align*}
  where we used the estimate
  \begin{align*}
   \sup_{u \in [0,T]} |a_u(x,x) - a_u(x,y) - a_u(y,x) + a_u(y,y)| \leq \alpha_1 |x - y|^{2\delta}
  \end{align*}
  in the last step.  
%   If $a_t$ is supported on some compact set $K$ for every $t \in [0,T]$, we can use the Burkholder Davis Gundy inequality to show that $M_t(x) = 0$ for every $t \in [0,T]$ and every $x \in D \setminus K$ on a set of full measure which, a priori, depends on $t$ and $x$. By a standard density argument, using continuity of $M$, we can deduce that there is a set of full measure on which $M_t$ is supported on $K$ for every $t \in [0,T]$. If $a$ does not have compact support, 
  Similarly, for $s < t$ and $x \in D$,
  \begin{align*}
   \EE | M_{ts}(x)  |^q \leq C q^{q/2} \EE \left| \int_s^t a_u(x,x) \, du \right|^{q/2} \leq C q^{q/2} |t-s|^{q/2} \left( \frac{\kappa}{1 + |x|^{\eta} }\right)^q
  \end{align*}
  for every $q \geq 2$. Hence we can can conclude with Theorem \ref{thm:kolmogorov_rough_drivers} (iii) that $M$ has a modification such that \eqref{eqn:1st_level_fin_and_gauss_tails} is finite and has Gaussian tails.

  Now we drop the assumption that $\alpha_1 < \infty$. Consider the stopped process $M_t^n := M_{t\wedge \tau_n}$ where
\begin{align*}
  \tau_n =  \inf \left\{ t \in [0,T]\,:\,  \sup_{u \in [0,t]}  \| a_u \|^{\wedge}_{\delta} \geq n \right\},
\end{align*}
   hence
  \begin{align*}
   \langle M^n_{\bullet}(x) , M^n_{\bullet}(y) \rangle_t = \int_0^{t \wedge \tau_n} a_u(x,y)\, du.
  \end{align*}
  Fix $n \in \NN$. As before, we have the estimates
  \begin{align*}
   \left\| M^n_{ts}(x) - M^n_{ts}(y) \right\|_{L^q} \leq C \sqrt{n} \sqrt{q} |t-s|^{\frac{1}{2}} |x - y|^{\delta}
  \end{align*}
  for every $x,y \in D$, $s < t$ and $q \geq 1$, and
  % . If $a$ is compactly supported, also $M$ and $M^n$ are compactly supported. Otherwise we have the estimate 
  \begin{align*}
   \left\| M^n_{ts}(x)  \right\|_{L^q} \leq C \sqrt{q} |t-s|^{1/2} \left( \frac{\kappa}{1 + |x|^{\eta} }\right)
  \end{align*}
  for every $x \in D$, $s < t$ and $q$ as in \eqref{eqn:q_large_enough_level_1}. Thus Theorem \ref{thm:kolmogorov_rough_drivers} implies that there is a modification of $M^n$ such that
  \begin{align*}
    \sup_{0 \leq s < t \leq T} \frac{\| M_{ts}^n \|_{\mathcal{C}^{\rho}}}{|t-s|^{\frac{1}{p}}}
  \end{align*}
  is finite almost surely for every $n \in \NN$. Now,
  \begin{align*}
    \PP \left( \sup_{0 \leq s < t \leq T} \frac{\big\|M_{ts}\big\|_{\mathcal{C}^{\rho}}}{|t-s|^{\frac{1}{p}}} = \infty \right) &\leq \PP \left( \sup_{0 \leq s < t \leq T} \frac{\big\|M_{ts}\big\|_{\mathcal{C}^{\rho}}}{|t-s|^{\frac{1}{p}}} \neq \sup_{0 \leq s < t \leq T} \frac{\big\|M^n_{ts}\big\|_{\mathcal{C}^{\rho}}}{|t-s|^{\frac{1}{p}}}  \right) \\
    &\leq \PP( \tau_n < T) \to 0
  \end{align*}
as $n \to \infty$, which shows that indeed
  \begin{align*}
   \sup_{0 \leq s < t \leq T} \frac{\big\|M_{ts}\big\|_{\mathcal{C}^{\rho}}}{|t-s|^{\frac{1}{p}}} < \infty
  \end{align*}
  almost surely. This shows (i) and the first part of (iii). 
  
  We proceed with (ii). As above, we first assume that
  \begin{align*}
   \alpha_2 :=  \left\| \sup_{u \in [0,T] } \| a_u \|^{\wedge}_{\delta} \right\|_{L^{\infty}} + \left\| \sup_{u \in [0,T] } \| b_u \|^{\wedge}_{\delta} \right\|_{L^{\infty}} < \infty
  \end{align*}
  and that \eqref{eqn:quadr_var_growth} holds uniformly for every $q \geq 1$. Let $x,y \in D$ such that $0 < |x - y| \leq 1$ and $s < t$. By the triangle inequality, for every $q \geq 2$,
\begin{equation*}
\begin{split}
   \big\| U_{ts}(x) - U_{ts}(y) \big\|_{L^\frac{q}{2}} \leq &\left\| \int_s^t (M_{us}(x) - M_{us}(y))N_{du}(x) \right\|_{L^{q/2}} \\ 
   &+ \left\| \int_s^t M_{us}(y) (N_{du}(x) - N_{du}(y)) \right\|_{L^{q/2}}.
\end{split}
\end{equation*}
  Note that
  \begin{align*}
   \langle \int_0^{\bullet} (M_u(x) - M_u(y))N_{du}(x) \rangle_t &= \int_0^t |M_u(x) - M_u(y)|^2\, \langle N_{\bullet}(x) \rangle_{du} \\ 
   &\leq \sup_{u \in [0,t]} |M_u(x) - M_u(y)|^2 \langle N_{\bullet}(x) \rangle_t.
  \end{align*}
  Therefore, applying twice the Burkholder Davis Gundy inequality,
  \begin{align*}
   \left\| \int_s^t (M_{us}(x) - M_{us}(y))N_{du}(x) \right\|_{L^{q/2}} &\leq C \sqrt{q} \left\| \sqrt{ \langle N_{\bullet}(x) \rangle_{ts} } \sup_{u \in [s,t]} |M_{us}(x) - M_{us}(y)| \right\|_{L^{q/2}} \\
   &\leq C \sqrt{q} \| \sqrt{ \langle N_{\bullet}(x) \rangle_{ts} } \|_{L^{\infty}} \big\| \sup_{u \in [s,t]} |M_{us}(x) - M_{us}(y)| \big\|_{L^{q/2}} \\
   &\leq C q \| \sqrt{ \langle N_{\bullet}(x) \rangle_{ts} } \|_{L^{\infty}} \| \sqrt{ \langle M_{\bullet}(x) - M_{\bullet}(y) \rangle_{ts}} \|_{L^{q/2}}.
  \end{align*}
  Now we have the estimate
  \begin{align*}
   \| \sqrt{ \langle N_{\bullet}(x) \rangle_{ts} } \|_{L^{\infty}} = \left\| \sqrt{ \int_s^t b_u(x,x)\, du } \right\|_{L^{\infty}} \leq \sqrt{\alpha_2} |t - s|^{\frac{1}{2}}.
  \end{align*}
  %and $\sup_{x \in D} \sup_{u \in [0,T]} \|\sqrt{b_u(x,x)} \|_{L^{\infty}}$ is finite both in the case of $b$ having compact support or satisfying the stated growth condition. 
  Furthermore, as seen above,
  \begin{align*}
   \| \sqrt{ \langle M_{\bullet}(x) - M_{\bullet}(y) \rangle_{ts}} \|_{L^{q/2}} \leq \sqrt{\alpha_2} |t - s|^{\frac{1}{2}} |x - y|^{\delta}
  \end{align*}
  which implies
  \begin{align*}
   \left\| \int_s^t (M_{us}(x) - M_{us}(y)) N_{du}(x) \right\|_{L^{q/2}} \leq C \alpha_2 q |t-s| |x - y|^{\delta}. 
  \end{align*}
  Similarly,
  \begin{align*}
   \left\| \int_s^t M_{us}(y) (N_{du}(x) - N_{du}(y)) \right\|_{L^{q/2}} \leq C \alpha_2 q |t-s| |x - y|^{\delta},
  \end{align*} 
  hence we have shown that
  \begin{align*}
   \| U_{ts}(x) - U_{ts}(y) \|_{L^{q/2}} \leq C \alpha_2 q |t-s| |x - y|^{\delta}
  \end{align*}
  holds for every $q \geq 2$, every $s < t$ and every $x,y \in D$ such that $0 < |x - y| \leq 1$. 
%   If $a$ and $b$ are compactly supported, it follows (e.g. by the Burkholder Davis Gundy inequality) that also $U$ is compactly supported. Assume now that $a$ and $b$ are not compactly supported, but that the stated growth conditions hold. 
  Similarly, for every $q \geq 2$, every $x \in D$ and every $s<t$,
  \begin{align*}
   \| U_{ts}(x) \|_{L^{q/2}} &\leq C q \| \sqrt{ \langle M_{\bullet}(x) \rangle_{ts} } \|_{L^{q}} \| \sqrt{ \langle N_{\bullet}(x) \rangle_{ts} } \|_{L^{q}} \\
   &\leq  \frac{C q |t-s|}{1 + |x|^{2\eta}}.
  \end{align*}
  Hence we can apply Theorem \ref{thm:kolmogorov_rough_drivers} to show that the process $U$ has a modification such that
  \begin{align*}
    \sqrt{\sup_{0 \leq s < t \leq T} \frac{\| U_{ts} \|_{\mathcal{C}^{\rho}}}{|t-s|^{\frac{2}{p}}}}
  \end{align*}
  is finite almost surely and has Gaussian tails. It remains to prove (ii) when $\alpha_2 = \infty$. In this case, we can consider the stopped processes  $M_t^n := M_{t\wedge \sigma_n}$ and $N_t^n := N_{t\wedge \sigma_n}$ where
\begin{align*}
  \sigma_n =  \inf \left\{ t \in [0,T]\,:\,  \sup_{u \in [0,t]}  \| a_u \|^{\wedge}_{\delta} + \sup_{u \in [0,t] } \| b_u \|^{\wedge}_{\delta} \geq n \right\}
\end{align*}
and proceed as above. We leave the details to the reader.
%   and set 
%   \begin{align*}
%    U_{ts}^n(x) := \int_s^t (M_{us}^n N^n_{du})(x).
%   \end{align*}
%   As before,
%   \begin{align*}
%    \| U^n_{ts}(x) - U^n_{ts}(y) \|_{L^{q/2}} \leq C n q |t-s| |x - y|^{\delta'}
%   \end{align*}
%   and
%    \begin{align*}
%    \| U^n_{ts}(x) \|_{L^{q/2}} % &\leq C \| \sqrt{ \langle M_{\bullet}^1(x) \rangle_{ts} } \|_{L^{q}} \| \sqrt{ \langle M_{\bullet}^2(x) \rangle_{ts} } \|_{L^{q}} \\
%    &\leq  \frac{C q |t-s|}{1 + |x|^{2\eta}}
%   \end{align*}
%   where the constant $C$ now depends on $n$ and $q$. Theorem \ref{thm:kolmogorov_rough_drivers} implies that for every $n \in \NN$, $U^n$ has a modification such that
%   \begin{align*}
%     \sup_{0 \leq s < t \leq T} \frac{\| U^n_{ts} \|_{\mathcal{C}^{\rho}}}{|t-s|^{\frac{2}{p}}}
%   \end{align*}  
%   is finite almost surely. We conclude as already seen for the first level process.
\end{Dem}
	
\medskip	

\begin{thm}\label{ThmLiftMartingaleVelocityFields}	{\sf
  Let $M$ be a continuous local martingale velocity field in $\mathcal{C}^{2,\delta}(D,\RR^d)$ with continuous local characteristic $a$ in $\hat{\mathcal{C}}^{2,\delta}_b$ for some $\delta \in (0,1]$.
  \begin{itemize}
   \item[(i)] Let $\rho \in (0,\delta)$ and $p \in (2,3)$. Assume that % either $a$ has compact support, or that 
   there is an $\eta \in (0, \infty)$ and a constant $\kappa > 0$ such that
  \begin{equation}
  \label{EqCondifionLifting}
   \sum_{0 \leq |\alpha| \leq 2} \sup_{u \in [0,T]} \Big\| \sqrt{ \partial_x^{\alpha} \partial_y^{\alpha} a_u(z,z) } \Big\|_{L^q} \leq \frac{\kappa}{1 + |z|^{ \eta}}
  \end{equation}
  for every $z \in D$ where $q$ satisfies
  \begin{align*}
    q > \max \left\{ \frac{\delta}{\eta(\delta - \rho)} + d \left(\frac{\delta}{\eta(\delta - \rho)} + \frac{1}{\delta - \rho} \right), \frac{1}{\frac{1}{2} - \frac{1}{p}} \right\}.
  \end{align*}
  
  Then $\mathbf{M} = (M,\mathbb{M})$, $\mathbb{M}$ being defined as in \eqref{eqn:def_second_order_process}, has a modification which is a weak geometric $(p,\rho)$-rough driver. We call $\mathbf{M}$ the \emph{natural lift} of $M$. \vspace{0.2cm}
  
  \item[(ii)] Assume that $\sup_{t \in [0,T]} \|a_t\|_{2 + \delta}^{\wedge}$ is an almost surely bounded random variable and that there is an $\eta \in (0,\infty)$ and a constant $\kappa > 0$ such that \eqref{EqCondifionLifting} holds uniformly for every $z \in D$ and every $q \geq 1$. 
%   \begin{align}\label{eqn:decay_loc_char}
%    \sum_{0 \leq |\alpha| \leq 2} \sup_{u \in [0,T]} \big\| \sqrt{ \partial_x^{\alpha} \partial_y^{\alpha} a_u(z,z) } \big\|_{L^{\infty}} \leq \frac{\kappa}{1 + |z|^{ \eta}}
%   \end{align}
%   holds for every $z \in D$. 
  
  Then for every $p \in (2,3)$ and $\rho \in (0,\delta)$, $\mathbf{M} = (M,\mathbb{M})$ has a modification which is a weak geometric $(p,\rho)$-rough driver, and the random variables
  \begin{align*}
   \sup_{0 \leq s < t \leq T} \frac{ \| M_{ts} \|_{\mathcal{C}^{2 + \rho}}}{|t - s|^{\frac{1}{p}}}\qquad \text{and} \qquad \sqrt{\sup_{0 \leq s < t \leq T} \frac{ \| W_{ts} \|_{\mathcal{C}^{1 + \rho}}}{|t - s|^{\frac{2}{p}}}},
  \end{align*}
  $W$ being defined as in \eqref{eqn:def_second_order_process}, have Gaussian tails.

  \end{itemize}  } 
\end{thm}

\ssk

\begin{Dem}
 The claim for $M$ follows by applying Proposition \ref{prop:ex_mod_ito_int} to $M$ and its derivatives. For $W$, we use the product rule in Proposition \ref{prop:diff_and_prod_rule_stoch_integral} for calculating the derivatives  and apply Proposition \ref{prop:ex_mod_ito_int} afterwards. The estimates for $W$ together with $M$ yield the claimed estimates for $\mathbb{M}$. We leave the details to the reader.
\end{Dem}

\ssk

\begin{rem}
 In the special case where $M$ or $a$ have compact support, i.e. when there exists a deterministic compact set $K \subset D$ such that $M$ resp. $a$ are supported on $K$ almost surely, assertion (i) of Theorem \ref{ThmLiftMartingaleVelocityFields} holds without any moment conditions on $a$. Indeed, this follows from the fact that the stopped processes in Proposition \ref{prop:ex_mod_ito_int} trivially satisfy the growth condition stated in the Kolmogorov theorem \ref{thm:kolmogorov_rough_drivers}, and this was the only point where these assumptions were needed.
\end{rem}

\bigskip
\bigskip

%%-------------------------------------------%%
\subsection{Stochastic and rough flows}
\label{SubsectionStochRoughFlows}
%%-------------------------------------------%% 

We keep in this Section the notations of the previous sections, and denote in particular by $(\mcF_t)_{0\leq t\leq T}$ a filtration to which the semimartingale velocity field $M$ is adapted. Assume that the local characteristic $a$ of $M$ satisfies the boundedness assumptions of point {\sf (i)} in Theorem \ref{ThmLiftMartingaleVelocityFields}. Then we can use Theorem \ref{ThmLiftMartingaleVelocityFields} to define the natural lift $\bf M$ of $M$ into a rough driver, and one can make sense of the rough flow $\varphi$ as pathwise solution to the equation
$$
d\varphi = {\bf M}(\varphi\,;dt)
$$
using Theorem \ref{ThmMain}. It follows from equation \eqref{EqApproxVarphiMu}, giving $\varphi_{ts}$ as a limit of compositions  of $\mu_{ba}$'s, that $\varphi$ is a semimartingale stochastic flow of homeomorphisms. One can read its local characteristics on the It\^o formula that it satisfies. Given $x,y\in\RR^d$ and $f,g\in\mcC^3_b$, we have
\begin{equation*}
\begin{split}
f\big(\varphi_{ts}(x)\big) = f(x) + \big(M_{ts}f\big)(x) + \frac{1}{2} \left\{\Big(\int_s^t M_{us}.M_{du} - M_{du}.M_{us}\Big)f\right\}(x) &+ \frac{1}{2}\,\big(M_{ts}^2f\big)(x) \\
&+ O\Big(|t-s|^\frac{3}{p}\Big),
\end{split}
\end{equation*}
with an $O(\cdot)$ term depending only on $\|\bf M\|$ and $\|f\|_{\mcC^3}$, with a similar formula for $g\big(\varphi_{ts}(y)\big)$. We read on this identity that
$$
\underset{h\downarrow 0}{\lim}\;\EE\left[\frac{f\big(\varphi_{t+h,t}(x)\big) - \varphi(x)}{h}\bigg| \mcF_t \right] \;\overset{L^1}{=}\; \underset{h\downarrow 0}{\lim}\;\EE\left[ \frac{1}{2}\,\Big(M_{t+h,t}^2f\Big)(x) \Big| \mcF_t \right], 
$$
and 
\begin{equation*}
\begin{split}
&\underset{h\downarrow 0}{\lim}\;\frac{1}{h}\,\EE\Big[\big\{f\big(\varphi_{t+h,t}(x)\big)-f(x)\big\}\big\{g\big(\varphi_{t+h,t}(y)\big)-g(y)\big\}\Big| \mcF_t\Big] \\
&= \underset{h\downarrow 0}{\lim}\;\frac{1}{h}\,\EE\Big[\big(M_{t+h,t}f\big)(x)\big(M_{t+h,t}g\big)(y)\Big| \mcF_t\Big] = \langle f,g\rangle_t(x,y).
\end{split}
\end{equation*}
So the semimartingale stochastic flow $\varphi$ has the same local characteristics as the semimartingale stochastic flow generated by the Stratonovich differential equation
\begin{equation}
\label{EqStratoSDE}
dx_t = M_{\circ dt}(x_t);
\end{equation} 
they coincide by Theorem \ref{ThmLeJanW82}, such as stated in Theorems 4.4.1 and 4.5.1 in Kunita's book, as assumption \eqref{EqCondifionLifting} on the local characteristic $a$ of $M$ is clearly stronger than the optimal assumptions of Kunita.

\ssk

\begin{thm} 
\label{ThmIdentificationStochasticRough}
Let $M$ be a continuous local martingale velocity field in $\mathcal{C}^{2,\delta}(\RR^d,\RR^d)$, for some $\delta \in \big(\frac{2}{3},1\big]$, with continuous local characteristic $a$ in $\widehat{\mathcal{C}}^{2,\delta}_b$. Let $\bf M$ be the rough driver associated with $M$ by Theorem \ref{ThmLiftMartingaleVelocityFields}. Under the condition that $M$ or $a$ have compact support or that the growth condition \eqref{EqCondifionLifting} holds, the rough flow solution to the differential equation 
$$
d\varphi = {\bf M}(\varphi\,;dt)  
$$
coincides with the stochastic flow generated by the Stratonovich differential equation \eqref{EqStratoSDE}.
\end{thm}

\bigskip

%%-------------------------------------%%
\subsection{Strong approximations}
\label{SubsectionStrongApproximations}
%%-------------------------------------%%

We give in this Section an example of use of the continuity of the It\^o map, in the setting of rough drivers and rough flows, by proving a Wong-Zaka\"i type theorem for semimartingale stochastic flows of maps. That is, we prove that such flows are limits in probability of flows generated by ordinary differential equations. Granted the continuity of the It\^o map, the core of the proof consists in showing that a rough lift of a continuous piecewise linear time interpolation of a semimartingale velocity field $M$ converges in probability to $\bf M$ in the topology of rough drivers.

\medskip

As in the last section, let
\begin{align*}
	M \colon [0,T] \to \mathcal{C}^{2,\delta}\big(D,\RR^d\big)
\end{align*}
be a continuous local martingale velocity field with quadratic variation  
\begin{align*}
  \big\langle M^i_\bullet(x), M^j_\bullet(y) \big\rangle_t = \int_0^t a^{ij}_s(x,y)\, ds
\end{align*}
and $\delta \in (0,1]$. Let $\mathcal{D} = \big\{0=t_0 < t_1 < \ldots < t_n = T\big\}$ be a partition of the interval $[0,T]$ and define the piecewise linear approximation of $M$ with respect to $\mathcal{D}$ as
\begin{align*}
	M_t^{\mathcal{D}} := M_{t_i} + (t-t_i) \, \frac{M_{t_{i+1}} - M_{t_i}}{t_{i+1} - t_i}\quad \text{if } t \in [t_i,t_{i+1}].
\end{align*}
Note that $\mathcal{D} \mapsto M^{\mathcal{D}}$ commutes with the spatial derivate, i.e. 
\begin{align*}
	\partial_{x_i}\big(M^{\mathcal{D}}\big) = \big(\partial_{x_i}M\big)^{\mathcal{D}} =: \partial_{x_i} M^{\mathcal{D}}.
\end{align*}
Define the mesh size of the partition by the formula $|\mathcal{D}| := \max_i \big| t_{i+1} - t_i \big|$. We define the first order differential operator
   \begin{align*}
   W^{\mathcal{D}}_{ts} &:= \frac{1}{2} \left( \int_s^t M^{\mathcal{D};i}_{us} \partial_i M^{\mathcal{D};k}_{du} - \int_s^t M^{\mathcal{D};i}_{du} \partial_i M^{\mathcal{D};k}_{us} \right) \partial_k  \\
   &= \frac{1}{2} \int_s^t (M^{\mathcal{D}}_{us}.M^{\mathcal{D}}_{du} - M^{\mathcal{D}}_{du}.M_{us}^{\mathcal{D}} )
  \end{align*}
via usual Riemann-Stieltjes integration. Then we set
\begin{align*}
 \mathbb{M}^{\mathcal{D}}_{ts} := W^{\mathcal{D}}_{ts} + \frac{1}{2} \, M^{\mathcal{D}}_{ts} M^{\mathcal{D}}_{ts}
\end{align*}
and 
\begin{align}\label{eqn:rough_driver_approx}
  \mathbf{M}^{\mathcal{D}} := \big(M^{\mathcal{D}},\mathbb{M}^{\mathcal{D}}\big).
\end{align}
 Our aim is to prove that $\mathbf{M}^{\mathcal{D}}$ converges towards the natural lift $\mathbf{M}$ of $M$ when $|\mathcal{D}| \to 0$ in probability (or even in $L^p(\Omega)$) in the topology of rough drivers. Note that
\begin{align*}
 W_{ts} = \frac{1}{2} \int_s^t (M_{us}.M_{du} - M_{du}.M_{us} ) = \frac{1}{2} \int_s^t (M_{us}.M_{\circ du} - M_{\circ du}.M_{us} ),
\end{align*}
hence it is enough to prove that the Riemann-Stieltjes integrals of the approximated processes converge towards the Stratonovich integrals (in the right topology), and this is what we are going to do.

\ssk 

  \begin{lem}\label{lemma:est_mart_and_approx}   {\sf
   	Let $M = \big(M^1,\ldots,M^d\big) \colon [0,T] \to \RR^d$ be a continuous local martingale and assume that
   	\begin{align*}
   		 \big\|\langle M \rangle_T\big\|_{L^\frac{q}{2}} \leq K < \infty
   	\end{align*} 
   	for some $q \geq 2$ and $K > 0$. Let $\mathbf{M}$ and $\mathbf{M^{\mathcal{D}}}$ be the associated rough paths lifts to $M$ and $M^{\mathcal{D}}$, i.e. $\mathbf{M} = (M, \mathbb{M})$ and $\mathbf{M}^{\mathcal{D}} = \big(M^{\mathcal{D}},\mathbb{M}^{\mathcal{D}}\big)$ where
   	\begin{align*}
%   	 \mathbb{M} = \int_{\Delta^2} \circ dM \otimes \circ dM \quad \text{and} \quad \mathbb{M}^{\mathcal{D}} = \int_{\Delta^2} dM^{\mathcal{D}} \otimes dM^{\mathcal{D}}
	  \mathbb{M} = \int_s^t M_{us} M_{\circ du} \quad \text{and} \quad \mathbb{M}^{\mathcal{D}} = \int_s^t M^{\mathcal{D}}_{us} M^{\mathcal{D}}_{du}
   	\end{align*}
   	are iterated Stratonovich, resp. Riemann-Stieltjes, integrals. Set

   \begin{align*}
    \varepsilon := \left\| \ \sup_{0 < v - u \leq |\mathcal{D}|} \| \mathbf{{M}}_{vu} \| \right\|_{L^q} 
   \end{align*}
   and assume that $\varepsilon \leq 1$.
   
  Then for every $\delta \in (0,1/3)$ there is a constant $C = C(\delta,q,K)$ such that
    \begin{align*}
     \left\| M_{ts} - M_{ts}^{\mathcal{D}} \right\|_{L^q} &\leq C \varepsilon^\frac{\delta}{8}  \,\big\| \langle M \rangle_{ts}\big\|_{L^\frac{q}{2}}^{\frac{1-\delta}{2}} \quad \text{and} \\
     \left\| \mathbb{M}_{ts} - \mathbb{M}_{ts}^{\mathcal{D}} \right\|_{L^\frac{q}{2}} &\leq C  \varepsilon^\frac{\delta}{4}\, \Big\| \langle M \rangle_{ts}\Big\|_{L^\frac{q}{2}}^{1-\delta}
    \end{align*}
    for every $s,t \in [0,T]$.   }
   \end{lem}
   
\begin{Dem}
 Let $d$ denote the Carnot-Caratheodoy metric on the step-two free nilpotent Lie group $G^2_d$ over $\RR^d$ [cf. \cite[Chapter 7]{FV10}). By interpolation and \cite[Proposition 8.15]{FV10}, for fixed $s<t$ and $2 < p' < p < 3$,
   \begin{align*}
    d(\mathbf{{M}}_{ts},\mathbf{M}^{\mathcal{D}}_{ts}) &\leq \left( d_{0-\text{H\"ol}}(\mathbf{{M}},\mathbf{{M}^{\mathcal{D}}})^{1-\frac{p'}{p}} d_{p'-\text{var};[s,t]}(\mathbf{{M}},\mathbf{{M}^{\mathcal{D}}})^{\frac{p'}{p}}  \right) \\
    &\lesssim  \left( d_{\infty}(\mathbf{{M}},\mathbf{{M}^{\mathcal{D}}})^{1-\frac{p'}{p}} + d_{\infty}(\mathbf{{M}},\mathbf{{M}^{\mathcal{D}}})^{\frac{1}{2} -\frac{p'}{2p}}(\|\mathbf{{M}} \|^{\frac{1}{2} -\frac{p'}{2p}}_{\infty} + \|\mathbf{{M}^{\mathcal{D}}}\|^{\frac{1}{2} -\frac{p'}{2p}}_{\infty} ) \right) \\
    &\qquad \times \left( \|\mathbf{{M}} \|^{p'/p}_{p'-\text{var};[s,t]} + \|\mathbf{{M}^{\mathcal{D}}}\|^{p'/p}_{p'-\text{var};[s,t]} \right).
   \end{align*}
   Taking the $q$-th moment, using H\"older's inequality and Cauchy-Schwarz, we obtain
   \begin{align*}
    &\EE\Big[ d\big(\mathbf{M}_{ts},\mathbf{M}^{\mathcal{D}}_{ts}\big)^q\Big] \\
    & \lesssim \left( \EE\Big[d_{\infty}\big(\mathbf{{M}},\mathbf{M}^{\mathcal{D}}\big)^q\Big] + \sqrt{ \EE\Big[ d_{\infty}\big(\mathbf{{M}},\mathbf{{M}^{\mathcal{D}}}\big)^q\Big] \, \Big( \EE\big[ \|\mathbf{{M}} \|^q_{\infty}\big] + \EE\big[ \|\mathbf{{M}^{\mathcal{D}}}\|^q_{\infty}\big] \Big) } \right)^{1- \frac{p'}{p}} \\
    &\qquad \times \left( \EE\Big[ \|\mathbf{{M}} \|^q_{p'-\text{var};[s,t]}\Big] + \EE\Big[ \|\mathbf{{M}^{\mathcal{D}}}\|^{q}_{p'-\text{var};[s,t]} \Big]\right)^\frac{p'}{p}.
   \end{align*}
   By \cite[Theorem 14.8 and Theorem 14.15]{FV10},
   \begin{align*}
    \EE\big[ \| \mathbf{M} \|_{\infty}^q\big] &\lesssim \EE\Big[ \langle M \rangle_T |^\frac{q}{2}\Big] \qquad \text{and} \\
    \EE\big[ \|\mathbf{{M}^{\mathcal{D}}}\|^q_{\infty}\big] &\leq \EE\Big[ \|\mathbf{M}^{\mathcal{D}}\|^q_{p-\text{var}}\Big] \lesssim \EE\Big[\langle M \rangle_T |^\frac{q}{2}\Big].
   \end{align*}
   Using the same theorems, we also have
   \begin{align*}
    \EE\Big[ \|\mathbf{{M}} \|^q_{p'-\text{var};[s,t]}\Big] \lesssim \EE\Big[ \langle M \rangle_{ts}^\frac{q}{2}\Big] \qquad \text{and} \qquad \EE\Big[ \|\mathbf{{M}^{\mathcal{D}}}\|^{q}_{p'-\text{var};[s,t]}\Big] \lesssim \EE\Big[\langle M \rangle_{ts}^\frac{q}{2}\Big].
   \end{align*}
   The estimate \cite[Equation (14.6) on p. 400]{FV10} gives
   \begin{align*}
    \EE\Big[d_{\infty}\big(\mathbf{{M}},\mathbf{{M}^{\mathcal{D}}}\big)^q\Big] &\lesssim \left(\EE\Big[  \sup_{0 < v - u \leq |\mathcal{D}|} \big\| \mathbf{{M}}_{vu} \big\|^q\Big]\right)^\frac{1}{4} \left( \EE\Big[\langle M \rangle_T |^\frac{q}{2}\Big] \right)^\frac{3}{4} + \EE\Big[  \sup_{0 < v - u \leq |\mathcal{D}|} \| \mathbf{{M}}_{vu} \|^q \Big] \\
    &\lesssim \varepsilon^{\frac{q}{4}}(K^{\frac{3q}{4}} + 1).
   \end{align*}
   Setting $\delta = 1 - \frac{p'}{p}$, we have shown that
   \begin{align*}
    \Big\| d\big(\mathbf{{M}}_{ts},\mathbf{M}^{\mathcal{D}}_{ts}\big)\Big\|_{L^q} \lesssim \varepsilon^\frac{\delta}{8} | \langle M \rangle_{ts}|_{L^\frac{q}{2}}^\frac{1-\delta}{2}.
   \end{align*}
   The result follows by equivalence of homogeneous norms on $G^2$ (\cite[Theorem 7.44]{FV10}).
\end{Dem}

\begin{prop}\label{prop:conv_strat_integral}   {\sf
Let $M,N \colon D \times [0,T] \to \RR$ be continuous local martingale fields adapted to the same filtration.  Assume that the quadratic variation of the processes is given by
\begin{align*}
 d \langle M_{\bullet}(x), M_{\bullet}(y) \rangle_t &= a_t(x,y)\quad \text{resp.} \\
 d \langle N_{\bullet}(x), N_{\bullet}(y) \rangle_t &= b_t(x,y)
\end{align*}
for every $x,y \in D$ and every $t \in [0,T]$. Moreover, assume that there is a $\delta \in (0,1]$ such that $a$ and $b$ have continuous modifications in the space $\hat{\mathcal{C}}^{0,\delta}_b$. Let $p > 2$ and $\rho \in (0,\delta)$.

\begin{itemize}
 \item[(i)] Assume
 %either that there is some compact set $K$ in $D$ such that $a_t$ is supported on $K \times K$ for every $t \in [0,T]$ almost surely, or 
 that there is an $\eta \in (0,\infty)$ and a constant $\kappa > 0$ such that
\begin{align*}
  \sup_{u \in [0,T]} \left\| \sqrt{a_u (x,x)} \right\|_{L^q} \leq \frac{\kappa}{1 + |x|^{\eta}}
\end{align*}
for every $x \in D$ where
\begin{align}\label{eqn:q_large_enough_level_1_conv}
 q > \max \left\{ \frac{\delta}{\eta(\delta - \rho)} + d \left(\frac{\delta}{\eta(\delta - \rho)} + \frac{1}{\delta - \rho} \right), \frac{1}{\frac{1}{2} - \frac{1}{p}} \right\}.
\end{align}
Let $M$ be the modification of the process given in Proposition \ref{prop:ex_mod_ito_int}. 

Then 
\begin{align}\label{eqn_conv_first_level}
 \sup_{0 \leq s < t \leq T} \frac{\big\| M_{ts} - M_{ts}^{\mathcal{D}} \big\|_{\mathcal{C}^{\rho}}}{|t-s|^\frac{1}{p}} \to 0
\end{align}
in probability when $|\mathcal{D}| \to 0$.

  \item[(ii)] Define
  \begin{align*}
    U_{ts}(x) := \int_s^t (M_{us} N_{\circ du})(x) \quad \text{and} \quad U^{\mathcal{D}}_{ts}(x) := \int_s^t (M^{\mathcal{D}}_{us} N^{\mathcal{D}}_{du})(x)
  \end{align*}
  as Stratonovich resp. Riemann-Stieltjes integral. Assume % either that $a$ and $b$ are compactly supported in the sense above, or 
  that there is an $\eta > 0$ and a constant $\kappa > 0$ such that
\begin{align}\label{eqn:quadr_var_growth_conv}
  \sup_{u \in [0,T]} \left\| \sqrt{a_u (x,x)} \right\|_{L^q} + \sup_{u \in [0,T]} \left\| \sqrt{b_u (x,x)} \right\|_{L^q} \leq \frac{\kappa}{1 + |x|^{\eta}}
\end{align}
for every $x \in D$ with $q$ sufficiently large as in \eqref{eqn:q_large_enough_level_1_conv}. Let $U$ be the modification of the process given in Proposition \ref{prop:ex_mod_ito_int}. 

Then
\begin{align}\label{eqn_conv_second_level}
 \sup_{0 \leq s < t \leq T} \frac{\big\| U_{ts} - U_{ts}^{\mathcal{D}} \big\|_{\mathcal{C}^{\rho}}}{|t-s|^\frac{2}{p}} \to 0
\end{align}
in probability when $|\mathcal{D}| \to 0$.

  \item[(iii)] Assume that $\sup_{u \in [0,T]} \| a_u \|^{\wedge}_{\delta}$ and $\sup_{u \in [0,T]} \| b_u \|^{\wedge}_{\delta}$  are almost surely bounded random variables % and that either both processes are compactly supported or 
  and that \eqref{eqn:quadr_var_growth_conv} holds for some $\eta \in (0,\infty)$ and $\kappa> 0$ uniformly for all $q \geq 1$. 
  
  Then the convergences in \eqref{eqn_conv_first_level} and \eqref{eqn_conv_second_level} hold in $L^q$ for every $q \geq 1$.
\end{itemize}   }
\end{prop}

\ssk

\begin{Dem}
  We will only give a proof in the case of $a$ and $b$ having compact support. The more general case works analogous using the stated growth conditions, as seen in the proof of Proposition \ref{prop:ex_mod_ito_int}. 
  We start with (i). Assume first that $\sup_{u \in [0,T]} \| a_u \|^{\wedge}_{\delta}$ is an almost surely bounded random variable. Fix $x,y \in D$ such that $0 < |x - y| \leq 1$ and let $s < t \in [0,T]$. Define the martingale $\hat{M} := M(x) - M(y)$ and let $\hat{\mathbf{M}}$ denote its canonical rough path lift (given by $\hat{M}$ and its iterated Stratonovich integrals) . From the Burkholder Davis Gundy inequality for enhanced martingales (cf. \cite[Theorem 14.8]{FV10}), for every $q \geq 2$,
\begin{align*}
 \left\| \sup_{0 < v - u \leq |\mathcal{D} |} \big\| \hat{\mathbf{M}}_{vu} \big\| \right\|_{L^q} \leq C \left\| \sqrt{ \langle \hat{M} \rangle_{|\mathcal{D}|} } \right\|_{L^q} \leq C|\mathcal{D}|^\frac{1}{2} |x-y|^{\delta}
\end{align*}
where the constant $C$ depends on the essential supremum of $\sup_{u \in [0,T]} \| a_u \|^{\wedge}_{\delta}$. We may assume that $|\mathcal{D}|$ is sufficiently small such the the right hand side of the equation is smaller than $1$. By Lemma \ref{lemma:est_mart_and_approx}, for every $q \geq 2$ and every $\beta \in (0,1/3)$,
  \begin{align*}
   &\Big\| M_{ts}(x) - M_{ts}^{\mathcal{D}}(x) - M_{ts}(y) + M_{ts}^{\mathcal{D}}(y) \Big\|_{L^q} \\
   \leq\ &C |\mathcal{D}|^\frac{\beta}{16} |x-y|^{\frac{\beta \delta}{8}} \, \big\| \langle M(x) - M(y) \rangle_{ts} \rangle \big\|_{L^\frac{q}{2}}^\frac{1-\beta}{2}  \\
   \leq\ &C |\mathcal{D}|^\frac{\beta}{16} |x-y|^{ \delta \big( 1 - \frac{7 \beta}{8} \big) } |t-s|^\frac{1-\beta}{2}. 
  \end{align*}
  We have already seen that 
  \begin{align*}
   \big\| M_{ts}(x) - M_{ts}(y)  \big\|_{L^q} \leq C |x-y|^{ \delta} |t-s|^\frac{1}{2},
  \end{align*}
  and by the triangle inequality,
  \begin{align*}
   \Big\| M_{ts}^{\mathcal{D}}(x) - M_{ts}^{\mathcal{D}}(y)  \Big\|_{L^q} \leq C |x-y|^{ \delta} |t-s|^\frac{1-\beta}{2}.
  \end{align*}
  Choosing $\beta$ small enough and $q$ large enough, we can apply Theorem \ref{thm:kolmogorov_rough_driver_distance} to see that
  \begin{align*}
   \left\| \sup_{s < t} \frac{\big\|M_{ts} - M_{ts}^{\mathcal{D}} \big\|_{\mathcal{C}^{\rho}}}{|t - s|^\frac{1}{p}} \right\|_{L^q} \leq C |\mathcal{D}|^\frac{\beta}{16}
  \end{align*}
  which shows the claim. If the essential supremum of $\sup_{u \in [0,T]} \| a_u \|^{\wedge}_{\delta}$ is not bounded, define the stopping times
  \begin{align*}
   \tau_n = \inf \left\{ t \in [0,T]\, :\, \sup_{u \in [0,t]} \| a_u \|^{\wedge}_{\delta} \geq n \right\}
  \end{align*}
  and set $M^n := M_{t \wedge \tau_n}$. We can repeat the argument above and conclude that 
  \begin{align*}
   \sup_{s < t} \frac{\Big\|M_{ts}^n - M_{ts}^{n;\mathcal{D}} \Big\|_{\mathcal{C}^{\rho}}}{|t - s|^\frac{1}{p}} \to 0
  \end{align*}
  in probability as $|\mathcal{D}|$ tends to $0$, and every $n \in \NN$. Fix some $\varepsilon > 0$ and some $n \in \NN$. Then
	\begin{align*}
	 &\PP\left(\left| \sup_{s < t} \frac{\Big\|M_{ts} - M_{ts}^{\mathcal{D}} \Big\|_{\mathcal{C}^{\rho}}}{|t - s|^\frac{1}{p}} \right| \geq \varepsilon \right) \\
	 &\quad\leq \PP\left(\left| \sup_{s < t} \frac{\Big\|M_{ts} - M_{ts}^{\mathcal{D}} \Big\|_{\mathcal{C}^{\rho}}}{|t - s|^\frac{1}{p}} -   \sup_{s < t} \frac{\Big\|M_{ts}^n - M_{ts}^{n;\mathcal{D}} \Big\|_{\mathcal{C}^{\rho}}}{|t - s|^\frac{1}{p}} \right| \geq \frac{\varepsilon}{2} \right) \\
	 &\quad\quad  + \PP\left(\left|   \sup_{s < t} \frac{\Big\|M_{ts}^n - M_{ts}^{n;\mathcal{D}} \Big\|_{\mathcal{C}^{\rho}}}{|t - s|^\frac{1}{p}} \right| \geq \frac{\varepsilon}{2} \right) \\
	 &\quad \leq \PP\left( \tau_n < T\right) + \PP\left(\left|  \sup_{s < t} \frac{\Big\|M_{ts}^n - M_{ts}^{n;\mathcal{D}} \Big\|_{\mathcal{C}^{\rho}}}{|t - s|^\frac{1}{p}} \right| \geq \frac{\varepsilon}{2} \right),
	\end{align*} 
	and the first term converges to $0$ for $n \to \infty$. This shows that indeed
  \begin{align*}
    \sup_{s < t} \frac{\Big\|M_{ts} - M_{ts}^{\mathcal{D}} \Big\|_{\mathcal{C}^{\rho}}}{|t - s|^\frac{1}{p}} \to 0
  \end{align*}
  in probability when $|\mathcal{D}| \to 0$.

\medskip
   
Now we prove (ii). As above, we first assume that $\sup_{u \in [0,T]} \| a_u \|^{\wedge}_{\delta}$ and $\sup_{u \in [0,T]} \| b_u \|^{\wedge}_{\delta}$ are almost surely bounded random variables. Let $x,y \in D$ with $0 < |x - y| \leq 1$ and $s < t$. We have already seen that in case of the It\^o integral, for every $q \geq 2$,
  \begin{align*}
   \left\| \int_s^t (M_{us} N_{ du})(x) - \int_s^t (M_{us} N_{ du})(y)\right\|_{L^\frac{q}{2}} \leq C |t-s| |x-y|^{\delta}.
  \end{align*}
  Moreover, we have
  \begin{align*}
   &\Big\| \big\langle M(x),N(x) \big\rangle_{ts} - \big\langle M(y),N(y) \big\rangle_{ts} \Big \|_{L^{q/2}} \\
   \leq\ &\Big\| \big\langle M(x) - M(y),N(x) \big\rangle_{ts} \Big\|_{L^{\frac{q}{2}}} + \Big\| \big\langle M(y),N(x)  - N(y) \big\rangle_{ts} \Big\|_{L^{\frac{q}{2}}} \\
   \leq\ & \Big\| \sqrt{ \big\langle M(x) - M(y) \big\rangle_{ts} } \Big\|_{L^q} \, \Big\| \sqrt{ \big\langle N(x) \big\rangle_{ts} } \Big\|_{L^q} + \Big\| \sqrt{ \big\langle N(x) - N(y) \big\rangle_{ts} } \Big\|_{L^q} \, \Big\| \sqrt{ \big\langle M(x) \big\rangle_{ts} } \Big\|_{L^q} \\
   \leq\ & C |t-s| |x - y|^{\delta}. 
  \end{align*}
  This implies that the same estimate holds for the Stratonovich integral, thus
  \begin{align*}
   \big\| U_{ts}(x) - U_{ts}(y) \big\|_{L^{\frac{q}{2}}} \leq C |t-s|\, |x - y|^{\delta}.
  \end{align*}
  By the triangle inequality,
  \begin{align*}
   &\Big\| U_{ts}^{\mathcal{D}}(x) - U_{ts}^{\mathcal{D}}(y) \Big\|_{L^{\frac{q}{2}}} \\
   \leq\ &\left\| \int_s^t (M_{us}^{\mathcal{D}}(x) - M_{us}^{\mathcal{D}}(y)) N_{du}^{\mathcal{D}}(x) \right\|_{L^{\frac{q}{2}}} + \left\| \int_s^t M_{us}^{\mathcal{D}}(y) (N_{du}^{\mathcal{D}}(x) - N_{du}^{\mathcal{D}}(y)) \right\|_{L^{\frac{q}{2}}}
  \end{align*}
  Now we define the martingale $\hat{M} := \left(\frac{M(x) - M(y)}{|x-y|^{\delta}},N(x)\right)$. We can check that for its quadratic variation, we have
  \begin{align*}
   \Big\| \big\langle \hat{M} \big\rangle_{ts} \Big\|_{L^{\frac{q}{2}}} \leq C|t-s|
  \end{align*}
  where the constant $C$ does not depend on $x$ or $y$. From \cite[Theorem 14.15]{FV10}, it follows that
  \begin{align*}
   \left\| \int_s^t \Big(M_{us}^{\mathcal{D}}(x) - M_{us}^{\mathcal{D}}(y)\Big) N_{du}^{\mathcal{D}}(x) \right\|_{L^{\frac{q}{2}}} \leq C |t-s| \, |x - y|^{\delta}.
  \end{align*}
  Similarly,
  \begin{align*}
   \left\| \int_s^t M_{us}^{\mathcal{D}}(y) \Big(N_{du}^{\mathcal{D}}(x) - N_{du}^{\mathcal{D}}(y)\Big) \right\|_{L^{\frac{q}{2}}} \leq C |t-s|\, |x - y|^{\delta}
  \end{align*}
and hence
  \begin{align*} 
   \big\| U_{ts}^{\mathcal{D}}(x) - U_{ts}^{\mathcal{D}}(y) \big\|_{L^{\frac{q}{2}}} \leq C |t-s| |x - y|^{\delta}.
  \end{align*}
Therefore, by the triangle inequality,
  \begin{align}\label{eqn:approx_sec_level_est1}
   \Big\| U_{ts}(x) - U_{ts}^{\mathcal{D}}(x) - U_{ts}(y) + U_{ts}^{\mathcal{D}}(y) \Big\|_{L^{\frac{q}{2}}} \leq C |t-s|\, |x - y|^{\delta}.
  \end{align}
On the other hand, we can apply Lemma \ref{lemma:est_mart_and_approx} to the martingale $M(x)$, $x$ fixed, to see that for every $q \geq 2$ and every $\beta \in \big(0,\frac{1}{3}\big)$,
\begin{align*}
\Big\| U_{ts}(x) - U_{ts}^{\mathcal{D}}(x) \Big\|_{L^{\frac{q}{2}}} \leq C |\mathcal{D}|^\frac{\beta}{8} \Big\| \big\langle M(x) \big\rangle_{ts} \Big\|_{L^{\frac{q}{2}}}^{1-\beta} \leq C |\mathcal{D}|^\frac{\beta}{8}\, |t-s|^{1 - \beta}.
\end{align*}
The same estimate holds if we replace $x$ by $y$, and by the triangle inequality, we also get
\begin{align}\label{eqn:approx_sec_level_est2}
\Big\| U_{ts}(x) - U_{ts}^{\mathcal{D}}(x) - U_{ts}(y) + U_{ts}^{\mathcal{D}}(y) \Big\|_{L^{\frac{q}{2}}} \leq C |\mathcal{D}|^\frac{\beta}{8} |t-s|^{1 - \beta}.
\end{align}
Interpolating the two inequalities \eqref{eqn:approx_sec_level_est1} and \eqref{eqn:approx_sec_level_est2}, we see that for every $q \geq 2$, every $\beta \in (0,1/3)$ and every $\lambda \in [0,1]$, we have the estimate
\begin{align*}
\Big\| U_{ts}(x) - U_{ts}^{\mathcal{D}}(x) - U_{ts}(y) + U_{ts}^{\mathcal{D}}(y) \Big\|_{L^{\frac{q}{2}}} \leq C |\mathcal{D}|^\frac{\lambda \beta}{8} |x - y|^{(1-\lambda)\delta} |t-s|^{1 - \lambda \beta}.
\end{align*}
Choosing $\lambda > 0$ and $\beta > 0$ small enough and $q$ large enough, we can again use Theorem \ref{thm:kolmogorov_rough_driver_distance} to see that
\begin{align*}
\left\| \sup_{s < t} \frac{\big\| U_{ts} - U_{ts}^{\mathcal{D}} \big\|_{\mathcal{C}^{\rho}}}{|t - s|^\frac{2}{p}} \right\|_{L^{\frac{q}{2}}} \leq C |\mathcal{D}|^\frac{\lambda \beta}{8}
\end{align*}
which proves the claim if $\sup_{u \in [0,T]} \| a_u \|^{\wedge}_{\delta}$ and $\sup_{u \in [0,T]} \| b_u \|^{\wedge}_{\delta}$ are bounded random variables. The general case follows by the same stopping argument as above.
\end{Dem}

\medskip

\begin{thm}\label{ThmConvergenceMartingaleVelocityFields}	{\sf
  Let $M$ be a continuous local martingale velocity field in $\mathcal{C}^{2,\delta}(D,\RR^d)$, for some $\delta \in (0,1]$, 
  with continuous local characteristic $a$ in $\hat{\mathcal{C}}^{2,\delta}_b$.
  \begin{itemize}
   \item[(i)] Let $\rho \in (0,\delta)$ and $p \in (2,3)$. Assume that % either $a$ has compact support, or that 
   there is an $\eta \in (0, \infty)$ and a constant $\kappa > 0$ such that
  \begin{equation}\label{eqn:growth_cond_conv_thm}
   \sum_{0 \leq |\alpha| \leq 2} \sup_{u \in [0,T]} \Big\| \sqrt{ \partial_x^{\alpha} \partial_y^{\alpha} a_u(z,z) } \Big\|_{L^q} \leq \frac{\kappa}{1 + |z|^{ \eta}}
  \end{equation}
  for every $z \in D$ where $q$ satisfies
  \begin{align*}
    q > \max \left\{ \frac{\delta}{\eta(\delta - \rho)} + d \left(\frac{\delta}{\eta(\delta - \rho)} + \frac{1}{\delta - \rho} \right), \frac{1}{\frac{1}{2} - \frac{1}{p}} \right\}.
  \end{align*}
  Let $\mathbf{M} = (M,\mathbb{M})$  be the weak geometric $(p,\rho)$-rough driver given in Theorem \ref{ThmLiftMartingaleVelocityFields}. 
  
  Then
  \begin{align}\label{eqn:conv_rough_driver}
   \mathbf{M}^{\mathcal{D}} \to \mathbf{M}
  \end{align}
  in probability for $|\mathcal{D}| \to 0$, with $\mathbf{M}^{\mathcal{D}}$ defined as in \eqref{eqn:rough_driver_approx}.
\vspace{0.2cm}
  
  \item[(ii)] Assume that $\sup_{u \in [0,T]} \| a_u \|^{\wedge}_{\delta}$ is almost surely bounded, and that % either $a$ has compact support or that 
  the growth condition \eqref{eqn:growth_cond_conv_thm} holds for some $\eta \in (0,\infty)$ and a constant $\kappa > 0$ uniformly for all $q \geq 1$. 
%   such that
%   \begin{align*}
%    \sum_{0 \leq |\alpha| \leq 2} \sup_{u \in [0,T]} \left\| \sqrt{ \partial_x^{\alpha} \partial_y^{\alpha} a_u(z,z) } \right\|_{L^{\infty}} \leq \frac{\kappa}{1 + |z|^{ \eta}}
%   \end{align*}
%   holds for every $z \in \RR^d$. 
  
  Then for every $p \in (2,3)$ and $\rho \in (0,\delta)$, $\mathbf{M} = (M,\mathbb{M})$ has a modification which is a weak geometric $(p,\rho)$-rough driver, and the convergence in \eqref{eqn:conv_rough_driver} holds in $L^q$ for every $q \geq 1$.

  \end{itemize}  } 
\end{thm}

\ssk

\begin{Dem}
 This follows by applying Proposition \ref{prop:conv_strat_integral} to $M$, $W$,and $M^{\mathcal{D}}$ and $W^{\mathcal{D}}$ and its derivatives. We use the product rule in Proposition \ref{prop:diff_and_prod_rule_stoch_integral} for calculating the derivatives of $W$ and Proposition \ref{prop:diff_smooth_case} for the derivatives of $W^{\mathcal{D}}$. The details are left to the reader.
\end{Dem}

\ssk

\begin{rem}
 As seen in the proof of Proposition \ref{prop:conv_strat_integral}, assertion (i) in Theorem \ref{ThmConvergenceMartingaleVelocityFields} holds without any moment conditions on $a$ in the case where $M$ or $a$ have compact support.
\end{rem}

\bigskip 

It follows then directly from this statement and the continuity of the It\^o solution map, Theorem \ref{ThmMain}, that the solution flow to the equation 
$$
d\varphi = {\bf M}(\varphi\,;dt)
$$
satisfies a Wong-Zaka\" i theorem. Using that the flow coincides with the one generated by the corresponding Stratonovich SDE, we obtain the following corollary.

\begin{cor}
 Let $M$ be a continuous local martingale velocity field in $\mathcal{C}^{2,\delta}(\RR^d,\RR^d)$, for some $\delta \in \big(\frac{2}{3},1\big]$, with local characteristic $a$ in $\widehat{\mathcal{C}}^{2,\delta}_b$. Assume that $M$ or $a$ have compact support or that the growth assumption \eqref{EqCondifionLifting} holds. Let $\varphi$ be the flow generated by the Stratonovich solution to
 \begin{align*}
  d \varphi = M(\varphi\,; \circ dt)
 \end{align*}
 and $\varphi^{\mathcal{D}}$ be the pathwise solution to
 \begin{align*}
  d \varphi^{\mathcal{D}} = M^{\mathcal{D}}\big(\varphi^{\mathcal{D}}\,; dt\big).
 \end{align*} 
 Then $\varphi^{D} \to \varphi$ in the space of $\mathcal{C}^{\rho}$ homeomorphisms uniformly in probability when $|\mathcal{D}| \to 0$ for all $\rho \in (0,\delta)$.
\end{cor}

\bigskip

%%----------------------------------------%%
\section{Application: Large deviations}
\label{SectionLDPSupp}
%%----------------------------------------%%
 
We provide in this Section another illustration of use of the continuity of the It\^o map by proving a large deviation theorem for Brownian flows. Relatively few works were dedicated to these topices, and we mention \cite{BDM10} and \cite{DD12}. In \cite{BDM10}, Dupuis and his co-authors use Dupuis' weak convergence approach to large deviation principles to prove such a result for Brownian flows of maps, building on a general large deviation criterion proved earlier in \cite{BDM08}. Dereich and Dimitroff's approach in \cite{DD12} is more in the line of the present work. They consider Brownian flows of maps as solutions to rough differential equations driven by a Banach space valued Brownian rough path, whose construction in a vector field setting was made possible by the previous work \cite{Der10} of Dereich. The support and large deviation theorems for Brownian flows are then inherited from the corresponding results proved in \cite{Der10} for the above mentioned vector field-valued Brownian rough 
% path, via the classical It\^o map. Our approach here is similar to the latter, and perhaps more direct.  

\bigskip

Let $(E,\mathcal{H}_1,\gamma)$ be a Gaussian Banach space with norm $\| \cdot \|$, i.e. $(E, \| \cdot \|)$ is a separable Banach space, $\gamma$ is a Gaussian measure defined on the Borel $\sigma$ algebra and $\mathcal{H}_1$ the Cameron-Martin Hilbert space (cf. \cite{Bog98} or \cite[Chapter 4]{Led96} for the precise definitions and further properties). Recall that $\mathcal{H}_1$ is continuously embedded in $E$, and for every $h \in \mathcal{H}_1$,
\begin{align*}
\|h \| \leq \sigma_{\gamma} \| h \|_{\mathcal{H}_1} = \sigma_{\gamma} \sqrt{\langle h, h \rangle_{\mathcal{H}_1}}
\end{align*}
where
\begin{align*}
\sigma_{\gamma}^2 := \int_E \|x \|^2\, \gamma(dx).
\end{align*}
A process $X \colon [0,T] \to E$, defined on some probability space, is called a \emph{$E$-valued Wiener process} if it has almost surely continuous sample paths starting from $0$, has independent increments, and for every $\xi \in E^*$, the distribution of $\langle X_t - X_s, \xi \rangle$ is a centered Gaussian random variable with variance  $|t-s| |\xi|_{\mathcal{H}_1}^2$ (cf. \cite{LLQ02} and \cite{Der10} for more properties about $E$-valued Wiener processes). The law of $X$ on the space $\mathcal{C}\big([0,T],E\big)$ is again Gaussian, and one can see that the corresponding Cameron-Martin space $\mathcal{H}$ is given by
\begin{align*}
\mathcal{H} = \left\{ \int_0^\bullet f_s \, ds\,:\, f \in L^2\big([0,T],\mathcal{H}_1\big) \right\}
\end{align*}
where the integral is a Bochner-integral.  Moreover, if $h^i_t = \int_0^t \dot{h}^i_s\, ds$, $i = 1,2$, the scalar product is given by
\begin{align*}
\big\langle h^1, h^2 \big\rangle_{\mathcal{H}} = \int_0^T \big\langle \dot{h}^1_s, \dot{h}^2_s \big\rangle_{\mathcal{H}_1}\, ds.
\end{align*}

Let $p \in [1,\infty)$. In the following, we will use the notion of \emph{$p$-variation} of a path $h \colon [0,T] \to E$ which is defined as
\begin{align*}
\|h\|_{p-\text{var};[s,t]} := \sup_{(t_i) \subset [s,t]} \left( \sum_{t_i} \big\|h_{t_{i+1}} - h_{t_i}\big\|^p \right)^{\frac{1}{p}}
\end{align*}
where the supremum is taken over all finite partitions $(t_i)$ of the interval $[s,t]$. If $\|h\|_{1-\text{var};[0,T]} < \infty$, we say that $h$ has finite variation.

\begin{lem}\label{lemma:smoothness_CM_paths}
For every $h \in \mathcal{H}$ we have
\begin{align}\label{eqn:CM_paths_12_Hoelder}
\sup_{0 \leq s < t \leq T} \frac{\big\| h_t - h_s \big\|}{|t-s|^{\frac{1}{2}}} \leq \sigma_{\gamma} \sqrt{ \langle h, h \rangle_{\mathcal{H}} }
\end{align}
and for $[s,t] \subseteq [0,T]$,
\begin{align}\label{eqn:CM_paths_bdd_var}
\|h\|_{1-\text{var};[s,t]} \leq \sigma_{\gamma} |t - s|^{\frac{1}{2}} \sqrt{ \langle h, h \rangle_{\mathcal{H}} }.
\end{align}
In particular, every $h \in \mathcal{H}$ is $\frac{1}{2}$-H\"older continuous and has finite variation on $[0,T]$.
\end{lem}
  
\ssk  
  
\begin{Dem}
Clearly, \eqref{eqn:CM_paths_12_Hoelder} follows from \eqref{eqn:CM_paths_bdd_var}, hence we only prove the second estimate. Let $h \in \mathcal{H}$ with $h(t) = \int_0^t \dot{h}_s \, ds$ and let $(t_i)$ be a partition of some interval $[s,t] \subseteq [0,T]$. Then
\begin{align*}
\sum_{i} \| h_{t_{i+1}} - h_{t_i} \| &\leq \sum_i \int_{t_i}^{t_{i+1}} \| \dot{h}_u \|\, du = \int_s^t  \| \dot{h}_u \|\, du \\
&\leq |t - s|^{\frac{1}{2}} \left( \int_0^T \| \dot{h}_u \|^2\, du \right)^{\frac{1}{2}} \leq \sigma_{\gamma} |t - s|^{\frac{1}{2}} \sqrt{ \langle h, h \rangle_{\mathcal{H}} }.
\end{align*}
Taking the supremum over all partitions shows the claim.
\end{Dem}

\ssk

Let $D$ be a an open, relatively compact, connected subset in $\RR^d$, $m \in \NN_0$ and $\delta \in (0,1]$. In the following, we would like to take the space $\mathcal{C}_b^{m,\delta}(D,\RR^d)$ as $E$ and consider a Gaussian measure on this space. However, $\mathcal{C}_b^{m,\delta}(D,\RR^d)$ is not separable (which is usually the case for H\"older-type spaces). Instead, we define the space $\mathcal{C}_b^{m,0,\delta}(D,\RR^d)$ as the closure of smooth paths from $D$ to $\RR^d$ with respect to the norm $\| \cdot \|_{m + \delta}$. As for H\"older spaces, using boundedness of $D$, one can show that these spaces are separable. From now on, let $E = \mathcal{C}_b^{m,0,\delta}(D,\RR^d)$ and assume that there is a Gaussian Banach space $(E,\mathcal{H}_1,\gamma)$. 

If $v$ is a $\mathcal{C}_b^{m,\delta}(D,\RR^d)$ valued path with finite variation and if $m \geq 1$, we define the pair $S(v)_{ts} = (v_{ts},\mathbbm{v}_{ts})$ by setting
  \begin{align*}
   v_{ts}(x) = v_t(x) - v_s(x)\qquad \text{and} \qquad \mathbbm{v}_{ts} = w_{ts} + \frac{1}{2} v_{ts} v_{ts}
  \end{align*}
  where $w_{ts}$ is the first order differential operator
  \begin{align*}
   w_{ts} = \frac{1}{2} \left( \int_s^t v_{us} . v_{du} - v_{du} . v_{us} \right)
  \end{align*}
  and the integral is a Riemann-Stieltjes integral. Note that if $X$ is a Wiener process in $\mathcal{C}^{m,0,\delta}_b$, $S(h)$ is always defined for every Cameron-Martin path $h$ since these paths are continuous and have bounded variation by Lemma \ref{lemma:smoothness_CM_paths}. Moreover, the following holds:
  
\begin{lem}\label{lemma:estim_Sh}
Let $(E,\mathcal{H}_1,\gamma)$ be a Gaussian Banach space with $E = \mathcal{C}^{2,0,\delta}_b(D,\RR^d)$ and $\delta \in (0,1]$. Then, for every $h \in \mathcal{H}$, $S(h)$ is a geometric $(2,\delta)$-rough driver, and there is a constant $C$ such that
\begin{align*}
\sup_{s<t} \frac{\|h_t - h_s \|_{\mathcal{C}^{2 + \delta}}}{|t-s|^{\frac{1}{2}}} \leq \sigma_{\gamma} \sqrt{ \langle h,h \rangle } \quad \text{and} \quad \sup_{s<t} \frac{\| w_{ts} \|_{\mathcal{C}^{1+\delta}}}{|t-s|} \leq C \sigma_{\gamma}^2 \, \langle h,h \rangle
\end{align*}
where
\begin{align*}
w_{ts} = \frac{1}{2} \left( \int_s^t h_{us} . h_{du} - h_{du} . h_{us} \right).
\end{align*}
\end{lem}

\ssk  
  
\begin{Dem}
Let $h \in \mathcal{H}$ and $S(h) = (h, \mathbbm{h})$ be defined as above. The claim for $h$ follows directly from Lemma \ref{lemma:smoothness_CM_paths}, and the algebraic condition for $(h, \mathbbm{h})$ follows from well-known identities for Riemann-Stieltjes integrals. Let $i,k \in \{1,\ldots, d\}$ , $x \in D$ and $s < t$. Then, by Riemann-Stieltjes estimates and Lemma \ref{lemma:smoothness_CM_paths},
\begin{align*}
    \left| \int_s^t h^i_{us}(x)\, \partial_i h^k_{du}(x) \right| &\leq  \sup_{u \in [s,t]} \big| h^i_{us}(x)\big| \sup_{(t_j) \subset [s,t]} \sum_{j} \big|\partial_i h^k_{t_{j+1}}(x) - \partial_i h^k_{t_j}(x)\big| \\
    &\leq \sup_{u \in [s,t]} \| h_{us} \|_{2 + \delta}  \sup_{(t_j) \subset [s,t]} \sum_{j} \big\| h_{t_{j+1}} -  h_{t_j} \big\|_{2 + \delta} \\
    &\leq \sigma_{\gamma}^2|t-s| \langle h,h \rangle.
\end{align*}
One can perform the same estimate for the second term in $w_{ts}$. By the triangle inequality, this shows that
\begin{align*}
\sup_{s<t} \frac{\| w_{ts} \|_{\mathcal{C}^0}}{|t-s|} \leq C \langle h,h \rangle.
\end{align*}
We can calculate the derivative of $w$ using Proposition \ref{prop:diff_smooth_case} and perform similar estimates, also for the H\"older norm, to conclude. 
  \end{Dem}

\bigskip

%%-----------------------------%%
\subsection{Schilder's theorem and Freidlin-Ventzel large deviations for stochastic flows}
\label{SubsectionLDP}
%%-----------------------------%%

In this section, we will prove a large deviation result for a $(p,\rho)$-rough driver $\mathbf{X}$ in the case where the underlying vector field $X$ is a Wiener process. We will do this by using the \emph{extended contraction principle} (see e.g. \cite[Theorem 4.2.23]{DZ98}), a strategy which has proven to be useful in rough paths theory (\cite{LQZ02}, \cite{MS06}, \cite{FV07}, \cite{FV10}). As a corollary, we obtain a Freidlin-Ventzel-type large deviation result for the flow generated by this driver. The key step is to prove that $\mathbf{X}^n := \mathbf{X}^{\mathcal{D}_n}$ with $\mathcal{D}_n = \left\{ \frac{kT}{n}\, :\, k = 0,\ldots, n \right\}$ is an \emph{exponentially good approximation} to $\mathbf{X}$. This is done in the Lemmata \ref{lemma:exp_good_approx_1} and \ref{lemma:exp_good_approx_2}. 

\medskip

For $\varepsilon > 0$, set $\delta_{\varepsilon} \mathbf{X} := (\varepsilon X,\delta_{\varepsilon}\mathbb{X})$ where
\begin{align*}
\delta_{\varepsilon}\mathbb{X}_{ts} &:= \delta_{\varepsilon} W_{ts} + \frac{\varepsilon^2}{2} X_{ts} X_{ts} \\
    &:= \frac{1}{2} \left(\int_s^t (\varepsilon X)_{us} .  (\varepsilon X)_{\circ du} (\varepsilon X)_{ \circ du} .  (\varepsilon X)_{us} \right)  + \frac{\varepsilon^2}{2} X_{ts} X_{ts} .
\end{align*}
We similarly define $\delta_{\varepsilon}\mathbf{X}^{n}$ with Riemann-Stieltjes integrals. Note that the homogeneous metric $\mathfrak{d}$ defined in \eqref{eqn:def_homog_metric} enjoys the property that $\mathfrak{d}(\delta_{\varepsilon} \mathbf{X}, \delta_{\varepsilon}\mathbf{X}^{n}) = \varepsilon \mathfrak{d}(\mathbf{X},\mathbf{X}^{n})$ for every $\varepsilon > 0$ (which is the reason why we call it \emph{homogeneous} metric). Note that since $X$ is Gaussian, the local characteristic $a$ is almost surely deterministic. Moreover, since $D$ is bounded, the growth condition \eqref{EqCondifionLifting} in Theorem \ref{ThmLiftMartingaleVelocityFields} is trivially satisfied for $a$ and possible derivatives.

\ssk
  
\begin{lem}\label{lemma:exp_good_approx_1}  {\sf 
Let $D$ be a relatively compact domain in $\RR^d$ and $X$ be a Wiener process in $\mathcal{C}^{2,0,\delta}_b$ for some $\delta \in (0,1]$. Let $\mathbf{X} = (X,\mathbb{X})$ denote its natural lift to a $(p,\rho)$-rough driver for some $\rho \in (0,\delta)$ and $p \in (2,3)$.
% with $p-2 < \rho$. 
Let $\eta > 0$ be fixed. Then the following holds:
\begin{align*}
\lim_{n \to \infty} \limsup_{\varepsilon \to 0} \varepsilon^2 \log \PP \Big( \mathfrak{d}_{p,\rho}(\delta_{\varepsilon} \mathbf{X}, \delta_{\varepsilon} \mathbf{X}^n)  > \eta \Big) = - \infty. 
\end{align*}   }
\end{lem}

\ssk
  
\begin{Dem}
Let $\varepsilon > 0$ and $n \in \NN$ be fixed. Since $X$ is Gaussian, the quadratic variation process $a$ is deterministic, and all estimates for $X$, $X^n$, $X - X^{n}$ and its iterated integrals in the proof of Proposition \ref{prop:conv_strat_integral} hold for $q = 2$. Moreover, $|X|_{L^q} \lesssim \sqrt{q} |X|_{L^2}$. The iterated Stratonovich and Riemann-Stieltjes integrals are both elements in the second inhomogeneous Wiener chaos, therefore
   \begin{align*}
    \left\| \int \circ dX \otimes \circ dX \right\|_{L^q} \lesssim q \left\| \int \circ dX \otimes \circ dX \right\|_{L^2} \quad \text{and} \quad \left\| \int dX^n \otimes dX^n \right\|_{L^q} \lesssim q \left\| \int dX^n \otimes dX^n \right\|_{L^2},
   \end{align*} 
    for all $q \geq 2$,  see e.g. \cite[Theorem D.8]{FV10}, and similar estimates hold for the other quantities. Therefore, we may apply Theorem \ref{thm:kolmogorov_rough_driver_distance} with $\kappa$ equal to a constant times $\sqrt{q}$ which shows that 
    \begin{align*}
     \big\| \mathfrak{d}_{p,\rho}( \mathbf{X},  \mathbf{X}^n) \big\|_{L^q} =  \alpha_n \sqrt{q}
    \end{align*}
    holds for all $q \geq 2$ and some constant $\alpha_n$. Repeating this argument for every $n \in \NN$, we obtain a sequence $(\alpha_n)$ converging to $0$ for $n \to \infty$. Thus
    \begin{align*}
     \PP \big( \mathfrak{d}_{p,\rho}(\delta_{\varepsilon} \mathbf{X}, \delta_{\varepsilon} \mathbf{X}^n) > \eta \big) &= \PP \left( \mathfrak{d}_{p,\rho} (\mathbf{X},  \mathbf{X}^n) > \frac{\eta}{\varepsilon} \right) \\
     &\leq \left( \frac{\varepsilon}{\eta} \right)^q q^{\frac{q}{2}} \alpha_n^q \\
     &\leq \exp\left[ q \log\left( \frac{\varepsilon \alpha_n \sqrt{q}}{\eta} \right) \right].
    \end{align*}
    Choosing $q = \varepsilon^{-2}$ we obtain the inequality
    \begin{align*}
     \varepsilon^2  \log \PP \big( \mathfrak{d}_{p,\rho} (\delta_{\varepsilon} \mathbf{X}, \delta_{\varepsilon} \mathbf{X}^n) > \eta \big) \leq \log ( \alpha_n / \eta )
    \end{align*}
from which the claim follows.
\end{Dem}

\ssk  
  
If $\mathcal{H}$ is the Cameron-Martin space of a $\mathcal{C}^{2,0,\delta}_b$-valued Wiener process and $v$ is a path with values in $\mathcal{C}^{2,0,\delta}_b$, set 
\begin{align*}
   I(v) := \begin{cases}
            \frac{1}{2} \langle v,v \rangle_{\mathcal{H}}, &\text{if } v \in \mathcal{H} \\
            + \infty &\text{otherwise.}
           \end{cases}
\end{align*}

\begin{lem}\label{lemma:exp_good_approx_2}  {\sf 
Let $\mathcal{H}$ be the Cameron-Martin space for some $\mathcal{C}_b^{2,0,\delta}$-valued Wiener process. Choose $\Lambda > 0$. Then 
    \begin{align*}
     \lim_{|\mathcal{D}| \to 0} \sup_{\{h \in \mathcal{H}\,:\, I(h) \leq \Lambda \}} \mathfrak{d}_{p,\delta}\Big(S(h^{\mathcal{D}}),S(h)\Big) = 0 
    \end{align*}
for every $p > 2$.  }
\end{lem}
  
\begin{Dem}
It is easy to check (cf. \cite[Proposition 5.20]{FV10} and Lemma \ref{lemma:estim_Sh}) that
   \begin{align*}
    \sup_{0 \leq s < t \leq T} \frac{\big\|h^{\mathcal{D}}_t - h^{\mathcal{D}}_s \big\|_{2 + \delta}}{|t-s|^{\frac{1}{2}}} \leq \sqrt{3} \sup_{0 \leq s < t \leq T} \frac{\|h_t - h_s \|_{2 + \delta}}{|t-s|^{\frac{1}{2}}} \leq \sqrt{3} \sigma_{\gamma} \sqrt{ \langle h,h \rangle}.
   \end{align*}
From Lemma \ref{lemma:estim_Sh}, we know that
   \begin{align*}
    \sup_{s<t} \frac{\| w_{ts} \|_{\mathcal{C}^{1+\delta}}}{|t-s|} \leq C \sigma_{\gamma}^2 \langle h,h \rangle
   \end{align*}
Now fix $j,k \in \{1,\ldots,d\}$, $s<t$ and $x \in D$. Then
   \begin{align*}
     \left| \int_s^t h_{su}^{\mathcal{D};j}(x)\, \partial_j h^{\mathcal{D};k}_{du}(x) \right| \leq \big\|h^{\mathcal{D};j}(x)\big\|_{1-\text{var};[s,t]}\, \big\|\partial_j h^{\mathcal{D};k}(x) \big\|_{1-\text{var};[s,t]}.
   \end{align*}
For $u \in [0,T]$, define $u_{\mathcal{D}} := \sup\big\{ t_i \in \mathcal{D}\,:\, t_i \leq u \big\}$ and $u^{\mathcal{D}} := \inf\big\{ t_i \in \mathcal{D}\,:\, t_i \geq u \big\}$. With this notation, using the estimates from Lemma \ref{lemma:smoothness_CM_paths},
\begin{align*}
    \Big\|h^{\mathcal{D};j}(x)\Big\|_{1-\text{var};[s,t]} &\leq \Big\|h^{\mathcal{D};j}(x)\Big\|_{1-\text{var};[s,s^{\mathcal{D}}]} + \Big\|h^{\mathcal{D};j}(x)\Big\|_{1-\text{var};[s^{\mathcal{D}},t_{\mathcal{D}}]} + \Big\|h^{\mathcal{D};j}(x)\Big\|_{1-\text{var};[t_{\mathcal{D}},t]} \\
    &\leq \frac{|s^{\mathcal{D}} - s|}{|s^{\mathcal{D}} - s_{\mathcal{D}}|} \, \big|h^{j}_{s^{\mathcal{D}}}(x) - h^{j}_{s_{\mathcal{D}}}(x)\big| + \big\|h^{j}(x)\big\|_{1-\text{var};[s^{\mathcal{D}},t_{\mathcal{D}}]} \\
    &\qquad + \frac{|t - t_{\mathcal{D}}|}{|t^{\mathcal{D}} - t_{\mathcal{D}}|}\, \big|h^{j}_{t^{\mathcal{D}}}(x) - h^{j}_{t_{\mathcal{D}}}(x)\big| \\
    &\leq \frac{\big|s^{\mathcal{D}} - s\big|}{|s^{\mathcal{D}} - s_{\mathcal{D}}|^{\frac{1}{2}}} \sup_{0 \leq u < v \leq T} \frac{\big\|h_v - h_u\big\|_0}{|v-u|^{\frac{1}{2}}} + \sigma_{\gamma} \big|t_{\mathcal{D}} - s^{\mathcal{D}}\big|^{\frac{1}{2}} \sqrt{ \langle h,h \rangle } \\
    &\qquad + \frac{\big|t - t_{\mathcal{D}}\big|}{|t^{\mathcal{D}} - t_{\mathcal{D}}|^{\frac{1}{2}}} \sup_{0 \leq u < v \leq T} \frac{\big\|h_v - h_u\big\|_0}{|v-u|^{\frac{1}{2}}} \\
    &\leq 3 \sigma_{\gamma} \sqrt{ \langle h,h \rangle } |t-s|^{\frac{1}{2}}.
   \end{align*}
A similar estimate holds for $\partial_j h^{\mathcal{D};k}$. Therefore,
   \begin{align*}
    \sup_{0\leq s < t \leq T}\frac{\left\| \int_s^t h_{su}^{\mathcal{D};j}(\cdot)\, \partial_j h_{du}^{\mathcal{D};k}(\cdot) \right\|_0}{|t-s|} \leq 9 \sigma_{\gamma}^2 \langle h,h\rangle.
   \end{align*}
Let $s<t$ and $x,y \in D$. We have
   \begin{align*}
    \left| \int_s^t h_{su}^{\mathcal{D};j}(x)\, \partial_j h^{\mathcal{D};k}_{du}(x) - \int_s^t h_{su}^{\mathcal{D};j}(y)\, \partial_j h^{\mathcal{D};k}_{du}(y) \right|
    &\leq \Big\|h^{\mathcal{D};j}(x) - h^{\mathcal{D};j}(y)\Big\|_{1-\text{var};[s,t]} \Big\|\partial_j h^{\mathcal{D};k} \Big\|_{1-\text{var};[s,t]} \\ & + \Big\|h^{\mathcal{D};j}\Big\|_{1-\text{var};[s,t]} \Big\|\partial_j h^{\mathcal{D};k}(x) - \partial_j h^{\mathcal{D};k}(y) \Big\|_{1-\text{var};[s,t]}
   \end{align*}
Similar to $h^{\mathcal{D};j}(x)$, one can estimate
   \begin{align*}
    \Big\|h^{\mathcal{D};j}(x) - h^{\mathcal{D};j}(y)\Big\|_{1-\text{var};[s,t]} &\leq 3 \sigma_{\gamma} \sqrt{ \langle h,h \rangle } \, |t-s|^{\frac{1}{2}} \, |x-y|^{\delta}
   \end{align*}
and similarly for $\partial_j h^{\mathcal{D};k}(x) - \partial_j h^{\mathcal{D};k}(y)$. Thus
   \begin{align*}
    \sup_{0\leq s < t \leq T}\frac{\left\| \int_s^t h_{su}^{\mathcal{D};j}(\cdot)\, \partial_j h^{\mathcal{D};k}_{du}(\cdot) \right\|_{\delta}}{|t-s|} \leq 9 \sigma_{\gamma}^2 \langle h,h\rangle.
   \end{align*}
Using the product rule (Proposition \ref{prop:diff_smooth_case}),
   \begin{align*}
    \partial_i \left(\int_s^t h_{su}^{\mathcal{D};j}(x)\, \partial_j h^{\mathcal{D};k}_{du}(x)\right) = \int_s^t \partial_i h_{su}^{\mathcal{D};j}(x)\, \partial_j h^{\mathcal{D};k}_{du}(x) + \int_s^t h_{su}^{\mathcal{D};j}(x)\, \partial_{i,j}^2 h^{\mathcal{D};k}_{du}(x)
   \end{align*}
and similar estimates as above, we can show that
   \begin{align*}
    \sup_{0\leq s < t \leq T}\frac{\Big\| w^{\mathcal{D}}_{ts} \Big\|_{1 + \delta}}{|t-s|} \leq C \sigma_{\gamma}^2 \langle h,h\rangle
   \end{align*}
where
   \begin{align*}
    w^{\mathcal{D}}_{ts} =  \frac{1}{2} \left( \int_s^t h^{\mathcal{D}}_{us} . h_{du}^{\mathcal{D}} - h_{du}^{\mathcal{D}} . h_{us}^{\mathcal{D}} \right).
   \end{align*}
This implies
   \begin{align*}
    \sup_{\mathcal{D}} \sup_{h\, :\, I(h) \leq \Lambda} \Big\| S\big(h^{\mathcal{D}}\big) \Big\|_{2,\delta} < \infty.
   \end{align*} 
Let $p > 2$. Then
   \begin{align*}
    \sup_{0 \leq s < t \leq T} \frac{\Big\|h^{\mathcal{D}}_{ts} - h_{ts}\Big\|_{2 + \delta}}{|t-s|^{\frac{1}{p}}} 
    &\leq \left( \sup_{0 \leq s < t \leq T} \frac{ \Big\|h^{\mathcal{D}}_{ts} - h_{ts}\Big\|_{2 + \delta}}{|t-s|^{\frac{1}{2}}} \right)^{\frac{2}{p}} 
    \left( \sup_{0 \leq s < t \leq T} \Big\|h^{\mathcal{D}}_{ts} - h_{ts}\Big\|_{2 + \delta} \right)^{1 - \frac{2}{p}} \\
    &\leq \left( \sup_{0 \leq s < t \leq T} \frac{ \Big\|h^{\mathcal{D}}_{ts}\Big\|_{2 + \delta}}{|t-s|^{\frac{1}{2}}} + \sup_{0 \leq s < t \leq T} \frac{\|h_{ts}\|_{2 + \delta}}{|t-s|^{\frac{1}{2}}} \right)^{\frac{2}{p}} \\
    &\qquad \times \left( \sup_{0 \leq s < t \leq T} \Big\|h^{\mathcal{D}}_{ts} - h_{ts}\Big\|_{2 + \delta} \right)^{1 - \frac{2}{p}}
   \end{align*}
and 
   \begin{align*}
    &\sup_{0 \leq s < t \leq T} \frac{\left\| \int_s^t h_{su}^{\mathcal{D};j}(\cdot)\, \partial_j h^{\mathcal{D};k}_{du}(\cdot) - \int_s^t h_{su}^{j}(\cdot)\, \partial_j h^{k}_{du}(\cdot) \right\|_{1 + \delta}}{|t-s|^{\frac{2}{p}}} \\
    \leq\ &\left( \sup_{0 \leq s < t \leq T} \frac{\left\| \int_s^t h_{su}^{\mathcal{D};j}(\cdot)\, \partial_j h^{\mathcal{D};k}_{du}(\cdot) \right\|_{1 + \delta}}{|t-s|} + \sup_{0 \leq s < t \leq T} \frac{\left\| \int_s^t h_{su}^{j}(\cdot)\, \partial_j h^{k}_{du}(\cdot) \right\|_{1 + \delta}}{|t-s|} \right)^{\frac{2}{p}} \\
    &\qquad \times \left( \sup_{0 \leq s < t \leq T} \left\| \int_s^t h_{su}^{\mathcal{D};j}(\cdot)\, \partial_j h^{\mathcal{D};k}_{du}(\cdot) - \int_s^t h_{su}^{j}(\cdot)\, \partial_j h^{k}_{du}(\cdot) \right\|_{1 + \delta} \right)^{1 - \frac{2}{p}}.
   \end{align*}
Therefore, the claim follows if we can prove
   \begin{align}\label{eqn:unif_conv_CM_1st}
    \lim_{|\mathcal{D}| \to 0} \sup_{h\, :\, I(h) \leq \Lambda}  \sup_{0 \leq s < t \leq T} \Big\|h^{\mathcal{D}}_{ts} - h_{ts}\Big\|_{2 + \delta} = 0
   \end{align}
and
   \begin{align}\label{eqn:unif_conv_CM_2nd}
    \lim_{|\mathcal{D}| \to 0} \sup_{h\, :\, I(h) \leq \Lambda} \sup_{0 \leq s < t \leq T} \left\| \int_s^t h_{su}^{\mathcal{D};j}(\cdot)\, \partial_j h^{\mathcal{D};k}_{du}(\cdot) - \int_s^t h_{su}^{j}(\cdot)\, \partial_j h^{k}_{du} (\cdot) \right\|_{1 + \delta} = 0.
   \end{align}
Concerning \eqref{eqn:unif_conv_CM_1st}, note that
  \begin{align*}
     \sup_{0 \leq s < t \leq T} \Big\|h^{\mathcal{D}}_{ts} - h_{ts}\Big\|_{2 + \delta} \leq 2 \sup_{t \in [0,T]} \Big\|h^{\mathcal{D}}_{t} - h_{t}\Big\|_{2 + \delta}.
  \end{align*}
If $t \in [0,T]$ is fixed, using Lemma \ref{lemma:smoothness_CM_paths}, we have
  \begin{align*}
   \Big\|h^{\mathcal{D}}_{t} - h_{t}\Big\|_{2 + \delta} &= \left\| (t - t_{\mathcal{D}})\frac{h_{t^{\mathcal{D}}} - h_{t_{\mathcal{D}}}}{t^{\mathcal{D}} - t_{\mathcal{D}}}  + h_{t_{\mathcal{D}}} - h_t \right\|_{2 + \delta} \\
   &\leq 2 \sup_{|v - u| \leq |\mathcal{D}|} \|h_v - h_u \|_{2 + \delta} \\
   &\leq 2 \sqrt{|\mathcal{D}|} \sup_{|v - u| \leq |\mathcal{D}|} \frac{\|h_v - h_u \|_{2 + \delta}}{|v-u|^{\frac{1}{2}}} \\
   &\leq 2 \sigma_{\gamma} \sqrt{ \langle h,h\rangle_{\mathcal{H}}} \, \sqrt{|\mathcal{D}|},
  \end{align*}
and \eqref{eqn:unif_conv_CM_1st} follows. Fix $x \in D$. By Young integration estimates (cf. e.g. \cite[Theorem 6.8]{FV10}),
  \begin{align*}
   &\left| \int_s^t h_{su}^{\mathcal{D};j}(x)\, \partial_j h^{\mathcal{D};k}_{du}(x) - \int_s^t h_{su}^{j}(x)\, \partial_j h^{k}_{du}(x) \right| \\
   \lesssim\ &\Big\| h^{\mathcal{D};j} - h^j \Big\|_{2-\text{var};[0,T]} \Big\|\partial_j h^{\mathcal{D};k}\Big\|_{1-\text{var};[0,T]} + \big\|h^j\big\|_{1-\text{var};[0,T]} \Big\| \partial_j h^{\mathcal{D};k} - \partial_j h^k \Big\|_{2-\text{var};[0,T]}.
  \end{align*}
From interpolation for the $p$-variation \cite[Proposition 5.5]{FV10} and our former estimates,
  \begin{align*}
   \Big\| h^{\mathcal{D};j} - h^j \Big\|_{2-\text{var};[0,T]} &\leq \left( \sup_{0 \leq u < v \leq T} \Big\| h^{\mathcal{D};j}_{uv} - h^j_{uv} \Big\|_{1 + \delta} \right)^{\frac{1}{2}} \left( \Big\| h^{\mathcal{D};j} \Big\|_{1-\text{var};[0,T]} + \big\|  h^j \big\|_{1-\text{var};[0,T]} \right)^{\frac{1}{2}} \\
   &\leq C \sigma_{\gamma} \sqrt{ \langle h,h \rangle_{\mathcal{H}} } \, |\mathcal{D}|^\frac{1}{4}.
   \end{align*}
A similar estimate holds for $\| \partial_j h^{\mathcal{D};k} - \partial_j h^k \|_{2-\text{var};[0,T]}$ and we obtain
   \begin{align*}
    \left\| \int_s^t h_{su}^{\mathcal{D};j}(\cdot)\, \partial_j h^{\mathcal{D};k}_{du}(\cdot) - \int_s^t h_{su}^{j}(\cdot)\, \partial_j h^{k}_{du}(\cdot) \right\|_0 \leq C \sigma_{\gamma}^2 \langle h,h \rangle_{\mathcal{H}} \, |\mathcal{D}|^\frac{1}{4}.
   \end{align*}
As above, one can use the product rule and obtain similar estimate for the H\"older norm of the derivative. This shows that
   \begin{align*}
   \sup_{0 \leq s < t \leq T} \left\| \int_s^t h_{su}^{\mathcal{D};j}(\cdot)\, \partial_j h^{\mathcal{D};k}_{du}(\cdot) - \int_s^t h_{su}^{j}(\cdot)\, \partial_j h^{k}_{du}(\cdot) \right\|_{1 + \delta} \leq C \sigma_{\gamma}^2 \langle h,h \rangle_{\mathcal{H}} \, |\mathcal{D}|^\frac{1}{4}
    \end{align*}
and \eqref{eqn:unif_conv_CM_2nd} follows. 
\end{Dem}

\bigskip

\begin{thm}[Schilder's theorem for Wiener rough drivers]\label{theorem:schilder_rough_driver} 
Let $D$ be a relatively compact domain and let $X$ be a Wiener process in $\mathcal{C}_b^{2,0,\delta}(D,\RR^d)$ with Cameron-Martin space $\mathcal{H}$. %Assume that the local characteristic satisfies the conditions stated in Theorem \ref{ThmLiftMartingaleVelocityFields}. 
Denote by $\mathbf{X}$ its natural lift to a $(p,\rho)$-rough driver. For $\varepsilon > 0$, set $\mathbf{P}_{\varepsilon} := \PP \circ (\delta_{\varepsilon} \mathbf{X})^{-1}$. Then the family $\{\mathbf{P}_{\varepsilon} \,:\, \varepsilon > 0\}$ of probability measures satisfies a large deviation principle on the space of rough drivers with speed $\varepsilon^{-2}$ and good rate function
\begin{align*}
 J(\mathbf{v}) = \begin{cases}
                  \frac{1}{2} \langle v,v \rangle &\text{if } \mathbf{v} = (v,\mathbbm{v}) \text{ and } v \in \mathcal{H} \\
                  + \infty &\text{otherwise.}
                 \end{cases}
\end{align*}

\end{thm}

\begin{Dem}
 The proof is standard, using the large deviation principle for Gaussian measure \cite[Section 3.4]{DS89}, the extended contraction principle \cite[Theorem 4.2.23]{DZ98} and the results in the Lemmas \ref{lemma:exp_good_approx_1} and \ref{lemma:exp_good_approx_2} (cf. e.g. \cite[Theorem 13.42]{FV10}). 
\end{Dem}

\ssk

As an immediate corollary, we obtain Freidlin-Ventzel large devations for a class of stochastic flows.

\ssk

\begin{thm}
Let $X$ be a Wiener process in $\mathcal{C}_b^{2,0,\delta}(D,\RR^d)$, for some $\delta \in \big(\frac{2}{3},1\big]$, %with local characteristic satisfying the conditions of Theorem \ref{ThmLiftMartingaleVelocityFields}. Let
and let $\varphi^{\varepsilon}$ be the flow generated by the Stratonovich solution to
 \begin{align*}
  d \varphi^{\varepsilon} = \varepsilon X(\varphi^{\varepsilon}\,; \circ dt).
 \end{align*}
 Let $\nu^{\varepsilon}$ denote the law of $\varphi^{\varepsilon}$ in the space of $\mathcal{C}^{\rho}$ homeomorphisms, $\rho \in \big(\frac{2}{3},\delta\big)$. Then the family $\{ \nu^{\varepsilon}\, :\, \varepsilon > 0 \}$ of probability measures satisfies a large deviation principle with speed $\varepsilon^{-2}$ and good rate function
 \begin{align*}
  L(\psi) = \inf \Big\{ J(\mathbf{v})\,:\, d \psi = \mathbf{v}(\psi\,; dt) \Big\}.
 \end{align*}
\end{thm}

\ssk

\begin{Dem}
 The Stratonovich solution equals the solution generated by the $(p,\rho)$-rough driver $\mathbf{X}$. Using Theorem \ref{theorem:schilder_rough_driver} and the pathwise continuity $\mathbf{X} \mapsto \varphi$, we can use the usual contraction principle in large deviation theory \cite[Theorem 4.2.1]{DZ98} to conclude.
\end{Dem}

\bigskip
\bigskip

%--------------------%
\section{Appendix}
\label{SectionAppendix}
%--------------------%

We provide in this Appendix an elementary regularity result for integrals depending on a parameter. 

\medskip

Let $\delta_{\varepsilon}$ be a standard Dirac sequence. If $I$ is a closed interval and $f \colon I \to \RR$ is a continuous function, let $\bar{f} \colon \RR \to \RR$ denote the unique continuous extension which coincides with $f$ on $I$ and which is constant outside this interval. Set $f^{\varepsilon} := \delta_{\varepsilon} \ast \bar{f}$. If $D$ is some subset of $\RR^d$ and if $f \colon D \times I \to \RR$ is a continuous function in time for every $x \in D$, set 
$$
f^{\varepsilon}(x,t) :=  \big(\delta_{\varepsilon} \ast \bar{f}(x,\cdot)\big)(t).
$$

\ssk

\begin{prop}\label{prop:diff_smooth_case}  {\sf
Let $D \subset \RR^d$ be an open set and let $f \colon D \times [0,T] \to \RR$ and $g \colon D \times [0,T] \to \RR$ be continuous. Assume that $f$ and $g$ are continuously differentiable on $D$ and that $f(x,0) = \partial_{x_i} f(x,0) = 0$ for every $x \in D$ and every $i = 1,\ldots,d$. Moreover, assume that there are $p, q \in [1, \infty)$ with $\frac{1}{p} + \frac{1}{q} > 1$ such that
 \begin{align*}
 	\sup_{(t_i)} \sum_{t_i} \big\| f (\cdot,t_{i+1}) - f (\cdot,t_{i}) \big\|_{\mathcal{C}^1}^p \quad \text{and} \quad \sup_{(t_i)} \sum_{t_i} \big\| g (\cdot,t_{i+1}) - g (\cdot,t_{i}) \big\|_{\mathcal{C}^1}^q
 \end{align*}
 are finite, where the suprema are taken over all finite partitions of the interval $[0,T]$. Then the Young integral (cf. e.g. \cite[Chapter 6]{FV10} for the precise definition) $\int_0^T f(x,t)\, g(x,dt)$ exists, is continuously differentiable for all $x \in D$ and the derivative is given by
 \begin{align*}
  \partial_{x_i} \left( \int_0^T f(x,t)\, g(x,dt) \right) = \int_0^T \partial_{x_i}  f(x,t)\, g(x,dt) + \int_0^T f(x,t)\, \partial_{x_i} g(x,dt)
 \end{align*}
 for all $i = 1,\ldots,d$.  } 
\end{prop}

\ssk

\begin{Dem}
One can suppose without loss of generality that $i = 1$. Fix $x \in D$. % and let $K = [a,b]$, $a<b$ be a closed interval such that $x \in K \subset D$. Let $y \in K$. Taking the spatial derivate of the approximated integral, an application of the product rule gives
 \begin{align*}
  \partial_{x_1} \left( \int_0^T f(x,t)\, g^{\varepsilon}(x,dt) \right) &= \partial_{x_1} \left( \int_0^T f(x,t) \partial_t( g^{\varepsilon}(x,t))\, dt \right) \\
  &= \int_0^T \partial_{x_1} f(x,t) \partial_t( g^{\varepsilon}(x,t))\, dt + \int_0^T f(x,t) \partial_t( (\partial_{x_1} g(x,\cdot))^{\varepsilon}(t))\, dt \\
  &= \int_0^T \partial_{x_1} f(x,t) \,g^{\varepsilon}(x,dt) + \int_0^T f(x,t) (\partial_{x_1} g(x,\cdot))^{\varepsilon}(dt). \\
 % &\to \int_I \partial_x f(x,t) \,g(x,dt) + \int_I f(x,t)\, \partial_x g(x,dt)
 \end{align*}
Let $q' > q$ such that $\frac{1}{p} + \frac{1}{q'} > 1$ . Let $U$ be some neighbourhood of $x$ and let $y \in U$. From Young estimates and interpolation, we obtain
 \begin{align*}
  &\left| \int_0^T \partial_y f(y,t) \,g^{\varepsilon}(y,dt) - \int_0^T \partial_y f(y,t) \,g(y,dt) \right| \\
  \leq\ &C \sup_{(t_i) \subset [0,T]} \left( \sum_{t_i} \big|\partial_{x_1} f(y,t_{t_{i+1}}) - \partial_{x_1} f(y,t_i)\big|^p \right)^{\frac{1}{p}} \\
  & \times \sup_{(t_i) \subset [0,T]} \left( \sum_{t_i} \big| g^{\varepsilon}(y,t_{i+1}) -  g(y,t_{i+1}) - g^{\varepsilon}(y,t_i) + g(y,t_i) \big|^{q'} \right)^{\frac{1}{q'}} \\
  \leq\ &C \sup_{(t_i) \subset [0,T]} \left(\sum_{t_i} \big\| f(\cdot,t_{i+1}) -  f(\cdot,t_i) \big\|_{\mathcal{C}^1}^p \right)^{\frac{1}{p}} \, \\
  & \times  \left\{ \sup_{(t_i) \subset [0,T]} \left( \sum_{t_i} \big\| g^{\varepsilon}(\cdot,t_{i+1}) - g^{\varepsilon}(\cdot,t_i) \big\|_{\mathcal{C}^1}^q \right)^{\frac{1}{q}} +  \sup_{(t_i)  \subset [0,T]} \left( \sum_{t_i} \big\| g(\cdot,t_{i+1}) - g(\cdot,t_i) \big\|_{\mathcal{C}^1}^q \right)^{\frac{1}{q}}\right\}^{\frac{q}{q'}} \\
  & \times 2^{1 - \frac{q}{q'}} \sup_{0 \leq t \leq T} \big\| g^{\varepsilon}(\cdot,t) - g(\cdot,t) \big\|_{\mathcal{C}^0}^{1 - \frac{q}{q'}}.
 \end{align*}
It is easy to check that 
 \begin{align*}
 	\sup_{(t_i) \subset [0,T]} \left( \sum_{t_i} \big\| g^{\varepsilon}(\cdot,t_{i+1}) - g^{\varepsilon}(\cdot,t_i) \big\|_{\mathcal{C}^1}^q \right)^{\frac{1}{q}} \leq \sup_{(t_i) \subset [0,T]} \left( \sum_{t_i} \big\| g(\cdot,t_{i+1}) - g(\cdot,t_i) \big\|_{\mathcal{C}^1}^q \right)^{\frac{1}{q}}.
 \end{align*}
Therefore, we obtain a bound of the form
 \begin{align*}
  \left| \int_0^T \partial_y f(y,t) \,g^{\varepsilon}(y,dt) - \int_0^T \partial_y f(y,t) \,g(y,dt) \right| \leq C \sup_{0 \leq t \leq T} \big\| g^{\varepsilon}(\cdot,t) - g(\cdot,t) \big\|_{\mathcal{C}^0}^{1 - \beta' / \beta}
 \end{align*}
where $C$ is independent of $y$ and $\varepsilon$. Thus,
 \begin{align*}
  \int_0^T \partial_{x_1} f(y,t) \,g^{\varepsilon}(y,dt) \to \int_0^T \partial_{x_1} f(y,t) \,g(y,dt)
 \end{align*}
uniformly in a neighbourhood around $x$ when $\varepsilon \to 0$. Similarly,
 \begin{align*}
  \int_0^T f(y,t) (\partial_{x_1} g(y,\cdot))^{\varepsilon}(dt) \to \int_0^T f(y,t) \partial_{y_1} g(y,dt)
 \end{align*}
uniformly in a neighbourhood around $x$ when $\varepsilon \to 0$. This shows differentiability in $x$ of the integral 
 $$
 \int_0^T f(x,t)\, g(x,dt)
 $$ 
and the claimed identity.
\end{Dem}

\bigskip
\bigskip

% -------------------------------------------------------------------------------------------
%                         bibliography
% -------------------------------------------------------------------------------------------

\bibliographystyle{alpha}
\bibliography{refs}

\end{document}